%% file: CLIFF101.tex
\documentclass{article}

\usepackage{maple2e,amsfonts}   
\MaplePromptSecondary = {\raise 1pt \hbox{$\phantom{\scriptstyle>}$\space}} 

\usepackage{amssymb}

\input Makro98                  
\newcommand{\MBS}{\raise 1pt \hbox{$\phantom{\scriptstyle>}$\space}}

\newcommand{\JJB}{\mathop{\, \JJ \,}\limits_B} 
\newcommand{\JJA}{\mathop{\, \JJ \,}\limits_A} 
\newcommand{\JJg}{\mathop{\, \JJ \,}\limits_g} 

\newcommand{\JJgg}{\mathop{\, \JJ \,}\displaylimits_g} 

\newcommand{\LLg}{\mathop{\, \LL \,}\limits_g} 

\newcommand{\LLgg}{\mathop{\, \LL \,}\displaylimits_g} 

\newcommand{\Bprod}[2]{\mathop{#1 #2}\limits_B} 
\newcommand{\gprod}[2]{\mathop{#1 #2}\limits_g} 

\newcommand{\wexp}[1]{e^{\w #1}} 

\newcommand{\bjinv}{\bj^{-1}} 

\newcommand{\ed}{\end{document}} 
\newcommand{\edc}{\end{document}} 

\newcommand{\Y}[2]{Y^{#1}_{#2}} 
\newcommand{\G}[2]{G^{#1}_{#2}} 
\newcommand{\Id}{{\bf 1}}

\newcommand{\alphaq}{\alpha_q}

\newcommand{\jtilde}{\tilde \bj}

\newcommand{\revside}[1]{#1 \, \tilde{}}
\newcommand{\reversion}[1]{#1 \, \tilde{}}
\newcommand{\revtop}[1]{\tilde{#1}}

\newcommand{\Star}[1]{#1^{\ast}}

\newcommand{\BK}{\mathbb{K}}
\newcommand{\BF}{\mathbb{F}}
\newcommand{\BC}{\mathbb{C}}
\newcommand{\BR}{\mathbb{R}}
\newcommand{\BH}{\mathbb{H}}
\newcommand{\diag}{\mbox{\rm diag}}

\newcommand{\AAT}{AA^T}
\newcommand{\ATA}{A^T A}
\newcommand{\SST}{\Sigma \Sigma^T}
\newcommand{\STS}{\Sigma^T \Sigma}
\newcommand{\VT}{V^T}
\newcommand{\UT}{U^T}
\newcommand{\AT}{A^T}

\newcommand{\phtab}[1]{%
                         \mbox{\phantom{ #1 }}}

\newcommand{\example}[2]{%
            \vgap
            {\bf #1 #2 }
            \vhgap}

\def\dim{\hbox{\rm dim\,}}

\def\span{\hbox{\rm span\,}}

\def\pin3{\Pin(3)}
\def\spin3{\Spin(3)}
\def\axis{{\bf axis}}
\def\qaxis{{\bf qaxis}}

\def\vhgap{\vskip 5pt}
\def\vgap{\vskip 10pt}

\def\beq{\begin{equation}}
\def\eeq{\end{equation}} 

\begin{document}

\title{Hecke Algebras, SVD, and Other Computational Examples with {\sc CLIFFORD}\thanks{Extended version of a talk presented at ``ACACSE'99:Applied Clifford Algebra in Cybernetics, Robotics, Image Processing and Engineering'', International Workshop as a Special Parallel Session of the 5th International Conference on Clifford Algebras and their Applications in Mathematical Physics,                              June 27-July 4, 1999, Ixtapa, Zihuatanejo, Mexico}}
\author{Rafa\l \ Ab\l amowicz \\
        Department of Mathematics, \\
        Tennessee Technological University\\
        Cookeville, TN 38505\\
        rablamowicz@tntech.edu}
\date{September 30, 1999}
\maketitle

\noindent
{\bf Abstract:}
{\sc CLIFFORD} is a Maple package for computations in Clifford algebras $\cl (B)$ of an arbitrary 
symbolic or numeric bilinear form $B.$ In particular, $B$ may have a non-trivial antisymmetric part. It is well known that the symmetric part $g$ of $B$ determines a unique (up to an isomorphism) Clifford structure on $\cl(B)$ while the antisymmetric part of $B$ changes the multilinear structure of $\cl(B).$ As an example, we verify Helmstetter's formula which relates Clifford product in $\cl(g)$ to the Clifford product in $\cl(B).$ Experimentation with Clifford algebras $\cl(B)$ of a general form~$B$ is highly desirable for physical reasons and can be easily done with {\sc CLIFFORD}. One such application includes a derivation of a representation of Hecke algebras in ideals generated by $q$-Young operators. Any element (multivector) of $\cl(B)$ is represented in Maple as a multivariate Clifford polynomial in the Grassmann basis monomials although other bases, such as the Clifford basis,  may also be used. Using the well-known isomorphism between simple Clifford algebras $\cl(Q)$ of a quadratic form $Q$ and matrix algebras through a faithful spinor representation, one can translate standard matrix algebra problems into the Clifford algebra language. We show how the Singular Value Decomposition of a matrix can be performed in a Clifford algebra. Clifford algebras of a degenerate quadratic form provide a convenient tool with which to study groups of rigid motions in robotics. With the help from {\sc CLIFFORD} we can actually describe all elements of $\Pin(3)$ and $\Spin(3).$ Rotations in $\BR^3$ can then be generated by unit quaternions realized as even elements in $\cl^{+}_{0,3}.$ Throughout this work all symbolic computations are performed with {\sc CLIFFORD} and its extensions.\\


\pagebreak

\tableofcontents


\section{Introduction}
\label{sec:intro}

A first working version of a Maple package {\sc CLIFFORD} was presented in Banff in 1995 \cite{Abl96}. From a modest program capable of symbolic computations in Clifford algebras of an arbitrary bilinear form, {\sc CLIFFORD} has grown to include $96$ main procedures, $21$ new Maple types, close to $4,000$ lines of code written in the Maple programming language, and an extensive on-line documentation. There is a number of special-purpose extensions available to {\sc CLIFFORD} such as {\sc suppl} and {\sc asvd} used in this paper \cite{CLIFFORD}. In fact, anyone who uses Maple can easily write additional procedures to tackle specific problems.

There are major advantages in using {\sc CLIFFORD} on a Computer Algebra System. One is an ability to solve equations and find the most general elements in the Clifford algebra satisfying given conditions. This approach has been presented in Sections \ref{sec:hecke} and \ref{sec:robotics}. In Section \ref{sec:hecke} Young operators in the Hecke algebra $H_{\BF}(3,q)$ are eventually found by systematically solving three equations that define them. Computations in this section were first reported in \cite{AblFauser99a} along with a physical motivation. There, experimentation with {\sc CLIFFORD} led to finding Young operators realized as idempotents in the Hecke algebra which in turn had been embedded into the even part of a Clifford algebra $\cl(B)$ of a suitable bilinear form $B.$ Such embedding was first given in \cite{Fauser-hecke} where the bilinear form $B$ was found so that the defining relations on Hecke generators were satisfied. Furthermore, it was shown in \cite{AblFauser99a} that Young operators corresponding to conjugate tableaux were related through the operation of reversion in the Clifford algebra, an idea that was first proposed in  \cite{Fauser-hecke}. Here, we only show the mechanics of the search for such operators, Garnir elements, and bases in the representation spaces as they were performed with {\sc CLIFFORD}. 

In the same spirit, in Section \ref{sec:robotics}, we describe a search for the elements in $\Pin(3)$ considered as a subgroup of the group of units in the Clifford algebra $\cl_{0,3}.$ Seven general types, not entirely exclusive, are eventually found through a systematic search and analysis. Then, the elements of $\Spin(3)$ are computed and related to unit quaternions. Rotations in coordinate planes and in a plane orthogonal to an arbitrary non-zero axis vector are described using quaternions realized as elements of $\cl^{+}_{0,3}.$ A symbolic formula describing the most general rotation is derived. Finally, using the ability of {\sc CLIFFORD} to compute in Clifford algebras of a degenerate quadratic form, the semi-direct product $\Spin(3) \rtimes \BR^3$ is shown to generate rigid motions on a suitable subspace of the Clifford algebra $\cl_{0,3,1}.$

The second advantage of using symbolic program like {\sc CLIFFORD} is its ability to compute with expressions containing totally undefined symbolic coefficients. It is possible, of course, like in Section~\ref{sec:robotics} to impose additional conditions on these coefficients when needed (by defining aliases for roots of polynomial equations). In Section \ref{sec:helm} we verify one of Helmstetter's formulas \cite{Helm87} that relates Clifford product in $\cl(B),$ the Clifford algebra of an arbitrary bilinear form $B,$ to the Clifford product in $\cl(g)$ where $g=g^T$ is the symmetric part of $B.$ A re-wording of the Helmstetter formula presented to the Author by Pertti Lounesto \cite{Lounesto_private} proved to be suitable for symbolic verification with {\sc CLIFFORD}. In fact, while this problem turned up to be a challenge for {\sc CLIFFORD} in view of its complexity, it also has helped to fine-tune the program to make such computations feasible. We will only illustrate computations in dimension $3;$ however, computations in dimension up to $9$ have been successfully completed.

The third major advantage of using {\sc CLIFFORD} shows up in Section \ref{sec:svd}: here, we perform within the same workspace not only symbolic computations with a Clifford algebra, but also with a linear algebra package built into Maple. This way we can illustrate in two low-dimensional examples two parallel approaches to the Singular Value Decomposition (SVD) of a matrix: one through the matrix algebra, and one through the Clifford algebra. In Section \ref{sec:comments} we comment on SVD based entirely on the Clifford algebra approach.

\section{Verification of the Helmstetter formula}
\label{sec:helm}

In his paper \cite{Helm87} Helmstetter studies canonical isomorphisms between Clifford algebras $\cl(Q)$ and $\cl(Q')$ of two quadratic forms $Q$ and $Q'$ defined on the same (real or complex) vector space $V.$ The forms are related via the identity $Q'(\bx) = Q(\bx) + B(\bx,\bx)$ for every $\bx \in V$ and some bilinear form $B$ on $V.$ Helmstetter constructs a {\it deformed\/} Clifford product $\ast$ on $\cl(Q)$ by extending the Clifford product $\bx\by$ of two elements $\bx$ and $\by$ in $V \hookrightarrow \cl(Q)$
$$
\bx \ast \by = \bx\by + B(\bx,\by)
$$
to all elements in $\cl(Q).$ Together with the new product $\ast,$ the Clifford algebra $\cl(Q)$ becomes a deformed Clifford algebra $\cl(Q,B).$ Given now two different bilinear forms $B$ and $B'$ on the quadratic space $(V,Q)$ such that $B(\bx,\bx) = B'(\bx,\bx)$ for every $\bx \in V,$ Helmstetter proves that there exists $F \in \bigw^{2} V$ such that 
$$
B'(\bx,\by)-B(\bx,\by) = <F,\,\bx \w \by>
$$
and the mapping
\beq
\phi: \cl(Q,B) \rightarrow \cl(Q,B'), 
\quad \quad 
u \mapsto e^{\w F} \JJ u
\label{eq:phiiso}
\eeq 
gives an isomorphism from $\cl(Q,B)$ to $\cl(Q,B')$ which acts as an identity on $V.$ In the above, $e^{\w F} \JJ u$ denotes the {\it left contraction\/} of $u$ by the exterior exponential of $F$ (see \cite{AblLoun96}, \cite{Lounesto}). A special case of (\ref{eq:phiiso}) occurs when $B$ is symmetric, that is, $B=g=g^T$ and $B'=g+A$ for some antisymmetric form $A.$

Thus, with a slight change of notation, let $B=g+A,\, g^T=g,\, A^T=-A$ and let us consider two Clifford algebras $\cl(g)$ and $\cl(B)$ on the same vector space $V.$  We have therefore three contractions: $\bx \JJB \by = B(\bx,\by),$ $\bx \JJA \by = A(\bx,\by),$ and 
$\bx \JJg \by = g(\bx,\by).$\footnote{The symbols $\bx \JJB \by,\, \bx \JJA \by,$ and $\bx \JJg \by$ denote the left contraction of $\by$ by $\bx$ with respect to $B,\,A,$ and $g$ respectively.} Then, the $B$-dependent Clifford product $\Bprod{u}{v}$ of any two elements $u$ and $v$ in $\cl(B)$ can be written \cite{Lounesto_private} solely in terms of the operations in $\cl(g)$ as
\beq
\Bprod{u}{v} = ((u \LLgg \wexp{F})(v \LLgg \wexp{F})) \LLgg \wexp{(-F)}
\label{eq:Bprod}
\eeq
where $u \LLg \wexp{F}$ denotes the {\it right contraction\/} of $u$ by $\wexp{F}$ with respect to $g.$ The product of $u \LLg \wexp{F}$ and $v \LLg \wexp{F}$ in (\ref{eq:Bprod}) is taken in $\cl(g).$ Let $\dim_\BR(V) = n.$ Then, the element 
$F \in \bigw^2 V$ is defined as
\beq
F = \sum_{K,L} (-1)^{|\pi(K,L)|}A_K \e_L \bjinv
\label{eq:defF}
\eeq
where $\bj^{-1}$ denotes the inverse of the unit pseudoscalar 
$\bj = \e_1 \w \e_2 \w \cdots \w \e_n$ in $\cl(g),$ the product $\e_L \bjinv$ is taken in $\cl(g),$ and the summation is taken over all multi indices $K=[k_1,k_2]$ and 
$L=[l_1,l_2,\ldots,l_s]$ satisfying the following relations:
$$
K \cap L = \emptyset, \; K \cup L = \{1,2,\ldots,n\}, \; n=2+s,\; 
K \mbox{ and } L \mbox{ are ordered by } <.
$$
$\pi(K,L)$ denotes a permutation which puts the list $[k_1,k_2,l_1,l_2,\ldots,l_s]$ in the standard order $[1,2,\ldots,n]$ and $|\pi(K,L)|$ equals $0$ or $1$ depending whether $\pi(K,L)$ is an even or odd element of~$S_n.$ In (\ref{eq:defF}) we have also adopted notation
$\e_L = \e_{l_1 l_2 \ldots l_s} = \e_{l_1} \w \e_{l_2} \w \cdots \w \e_{l_s}.$

Before we proceed to verify formula (\ref{eq:Bprod}) with {\sc CLIFFORD} \cite{CLIFFORD}, let's observe the following properties of the left and right contraction:
\beq
u \LLgg v = \bjinv((\bj u) \w v),\quad
u \JJgg v = (u \w (v \bj)) \bjinv, \quad \bjinv = \frac{\jtilde}{\det(g)},
\eeq
where $\tilde{~}$ is the $g$-dependent reversion in $\cl(g).$\footnote{From now on in this section we assume that $\det(g) \not= 0.$} Observe also that since $F \in \bigw^2 V,$ we have
\beq
\wexp{F} = \sum_{k=0}^{N} \displaystyle{\frac{F^{\w k}}{k!}},
\quad
N = \lfloor n/2 \rfloor
\eeq
where $F^{\w k} = F \w F \w \cdots \w F$ is the exterior product of $F$ computed $k$-times and $\lfloor \cdot \rfloor$ denotes the floor function. For example, for different values of $n,$ $F$ has the following form:
\beq
n = 2, \; F  =  -A_{12} \, \e_{12}^{-1},
\eeq
\beq
n = 3, \; 
    F  =  -(A_{12}\, \e_3 - A_{13} \, \e_2 + A_{23} \, \e_1) \, \e_{123}^{-1},
\eeq
\beq
n = 4, \; 
    F  =  -(A_{12}\, \e_{34} - A_{13} \, \e_{24} + A_{14} \, \e_{23} + 
            A_{23} \, \e_{14} + A_{34} \, \e_{12} - A_{24} \, \e_{13}) \, \e_{1234}^{-1}, 
\eeq
and so on. We will verify the validity of (\ref{eq:Bprod}) in a numeric and a symbolic case. In the Maple symbolic language, formula (\ref{eq:Bprod}) becomes:\footnote{In {\sc CLIFFORD}, the Clifford product $uv$ of two elements $u$ and $v$ can be entered as {\tt u \&c v} (the infix form) or as {\tt cmul(u,v)}.}

\vskip 3pt
\begin{maplelatex}
\mapleinline{inert}{2d}{cmul(u,v) =
RCg(cmulg(RCg(u,wexp(F,N)),RCg(v,wexp(F,N)),wexp(-F,N)));}{%
\[
\mathrm{cmul}(u, \,v)=\mathrm{RCg}(\mathrm{cmulg}(\mathrm{RCg}(u
, \,\mathrm{wexp}(F, \,N)), \,\mathrm{RCg}(v, \,\mathrm{wexp}(F, 
\,N)), \,\mathrm{wexp}( - F, \,N)))
\]
}
\end{maplelatex}
\vskip 3pt 
\noindent
with the {\sc CLIFFORD} procedures {\tt cmulg} and {\tt RCg} representing the Clifford product and the right contraction in $\cl(g)$ and {\tt wexp} giving the exterior exponential in $\bigw V.$ We limit our two examples to dimension $n=3.$ Computations presented in the following two sections can be extended with {\sc CLIFFORD} to higher dimensions.

\subsection{Numeric example when $n=3$}
\label{sec:num3}

Let's first assign an arbitrary matrix to $B,$ split $B$ into its symmetric and antisymmetric parts $g$ and $A,$ and compute the bivector $F$ with a procedure {\tt makeF}:
\begin{maplegroup}
\begin{mapleinput}
\mapleinline{active}{1d}{
dim:=3:eval(makealiases(dim,'ordered')):B:=matrix(dim,dim,[4,8,3,0,9,5,-2,1,7]);}{%
}
\end{mapleinput}
\mapleresult
\begin{maplelatex}
\mapleinline{inert}{2d}{B := matrix([[4, 8, 3], [0, 9, 5], [-2, 1, 7]]);}{%
\[
B :=  \left[ 
{\begin{array}{rrr}
4 & 8 & 3 \\
0 & 9 & 5 \\
-2 & 1 & 7
\end{array}}
 \right] 
\]
}
\end{maplelatex}
\end{maplegroup}
\begin{maplegroup}
\begin{mapleinput}
\mapleinline{active}{1d}{g,A:=splitB(B);}{%
}
\end{mapleinput}
\mapleresult
\begin{maplelatex}
\mapleinline{inert}{2d}{g, A := matrix([[4, 4, 1/2], [4, 9, 3], [1/2, 3, 7]]), matrix([[0, 4,
5/2], [-4, 0, 2], [-5/2, -2, 0]]);}{%
\[
g, \,A :=  \left[ 
{\begin{array}{crc}
4 & 4 & {\displaystyle \frac {1}{2}}  \\ [2ex]
4 & 9 & 3 \\
{\displaystyle \frac {1}{2}}  & 3 & 7
\end{array}}
 \right] , \, \left[ 
{\begin{array}{crc}
0 & 4 & {\displaystyle \frac {5}{2}}  \\ [2ex]
-4 & 0 & 2 \\
{\displaystyle \frac {-5}{2}}  & -2 & 0
\end{array}}
 \right] 
\]
}
\end{maplelatex}
\end{maplegroup}
\begin{maplegroup}
\begin{mapleinput}
\mapleinline{active}{1d}{F:=makeF(dim);}{%
}
\end{mapleinput}
\mapleresult
\begin{maplelatex}
\mapleinline{inert}{2d}{F := 2/91*e13+86/455*e12;}{%
\[
F := {\displaystyle \frac {2}{91}} \,\mathit{e13} + 
{\displaystyle \frac {86}{455}} \,\mathit{e12}
\]
}
\end{maplelatex}
\end{maplegroup}
\noindent
Next, we find the exterior exponentials of $F$ and $-F$ which we assign to Maple variables $F_1$ and $F_2$ respectively.
\begin{maplegroup}
\begin{mapleinput}
\mapleinline{active}{1d}{N:=floor(dim/2):F1:=wexp(F,N);F2:=wexp(-F,N);}{%
}
\end{mapleinput}
\mapleresult
\begin{maplelatex}
\mapleinline{inert}{2d}{F1 := 1+2/91*e13+86/455*e12;F2 := 1-2/91*e13-86/455*e12;}{%
\[
\mathit{F1} := 1 + {\displaystyle \frac {2}{91}} \,\mathit{e13}
 + {\displaystyle \frac {86}{455}} \,\mathit{e12}, \quad
\mathit{F2} := 1 - {\displaystyle \frac {2}{91}} \,\mathit{e13}
 - {\displaystyle \frac {86}{455}} \,\mathit{e12}.
\]
}
\end{maplelatex}
\end{maplegroup}
\noindent
Let $u$ and $v$ be two arbitrary elements in $\cl(B):$
\begin{maplegroup}
\begin{mapleinput}
\mapleinline{active}{1d}{u:=2+e1-e23+e123;v:=3-e3+e12+e23;}{%
}
\end{mapleinput}

\mapleresult
\begin{maplelatex}
\mapleinline{inert}{2d}{u := 2+e1-e23+e123;v := 3-e3+e12+e23;}{%
\[
u := 2 + \mathit{e1} - \mathit{e23} + \mathit{e123}, \quad
v := 3 - \mathit{e3} + \mathit{e12} + \mathit{e23}.
\]
}
\end{maplelatex}
\end{maplegroup}
\noindent
The Clifford product $\Bprod{u}{v}$ of $u$ and $v$ in $\cl(B),$ the left hand side of (\ref{eq:Bprod}), is then equal to
\begin{maplegroup}
\begin{mapleinput}
\mapleinline{active}{1d}{cmul(u,v);}{%
}
\end{mapleinput}
\mapleresult
\begin{maplelatex}
\mapleinline{inert}{2d}{-48*e3+79*Id+13*e13-6*e12+81*e2-8*e123-81*e1;}{%
\[
 - 48\,\mathit{e3} + 79\,\mathit{Id} + 13\,\mathit{e13} - 6\,
\mathit{e12} + 81\,\mathit{e2} - 8\,\mathit{e123} - 81\,\mathit{
e1}
\]
}
\end{maplelatex}
\end{maplegroup}
\noindent
while the right hand side of (\ref{eq:Bprod}) gives the same result:
\begin{maplegroup}
\begin{mapleinput}
\mapleinline{active}{1d}{RCg(cmulg(RCg(u,F1),RCg(v,F1)),F2);}{%
}
\end{mapleinput}
\mapleresult
\begin{maplelatex}
\mapleinline{inert}{2d}{-48*e3+79*Id+13*e13-6*e12+81*e2-8*e123-81*e1;}{%
\[
 - 48\,\mathit{e3} + 79\,\mathit{Id} + 13\,\mathit{e13} - 6\,
\mathit{e12} + 81\,\mathit{e2} - 8\,\mathit{e123} - 81\,\mathit{
e1}
\]
}
\end{maplelatex}

\end{maplegroup}

\subsection{Symbolic computations when $n=3$}
\label{sec:symb3}

A purely symbolic computation when $n=3$ will look as follows. Matrix $B$ is now defined as an arbitrary symbolic $3 \times 3$ matrix with a symmetric part $g$ and an antisymmetric part $A,$ and $F$ is again computed using the procedure {\tt makeF}. The exterior exponentials of $F$ and $-F$ are again denoted respectively by $F_1$ and $F_2.$ All symbolic parameters in $B$ are assumed to be real or complex.
\begin{maplegroup}
\begin{mapleinput}
\mapleinline{active}{1d}{dim:=3:eval(makealiases(dim)):
B:=matrix(dim,dim,[g11,g12+A12,g13+A13,g12-A12,g22,g23+A23,g13-A13,g23-A23,g33]);}{%
}
\end{mapleinput}

\mapleresult
\begin{maplelatex}
\mapleinline{inert}{2d}{B := matrix([[g11, g12+A12, g13+A13], [g12-A12, g22, g23+A23],
[g13-A13, g23-A23, g33]]);}{%
\[
B :=  \left[ 
{\begin{array}{ccc}
\mathit{g11} & \mathit{g12} + \mathit{A12} & \mathit{g13} + 
\mathit{A13} \\
\mathit{g12} - \mathit{A12} & \mathit{g22} & \mathit{g23} + 
\mathit{A23} \\
\mathit{g13} - \mathit{A13} & \mathit{g23} - \mathit{A23} & 
\mathit{g33}
\end{array}}
 \right] 
\]
}
\end{maplelatex}

\end{maplegroup}
\begin{maplegroup}
\begin{mapleinput}
\mapleinline{active}{1d}{g,A:=splitB(B);}{%
}
\end{mapleinput}

\mapleresult
\begin{maplelatex}
\mapleinline{inert}{2d}{g, A := matrix([[g11, g12, g13], [g12, g22, g23], [g13, g23, g33]]),
matrix([[0, A12, A13], [-A12, 0, A23], [-A13, -A23, 0]]);}{%
\[
g, \,A :=  \left[ 
{\begin{array}{ccc}
\mathit{g11} & \mathit{g12} & \mathit{g13} \\
\mathit{g12} & \mathit{g22} & \mathit{g23} \\
\mathit{g13} & \mathit{g23} & \mathit{g33}
\end{array}}
 \right] , \, \left[ 
{\begin{array}{ccc}
0 & \mathit{A12} & \mathit{A13} \\
 - \mathit{A12} & 0 & \mathit{A23} \\
 - \mathit{A13} &  - \mathit{A23} & 0
\end{array}}
 \right] 
\]
}
\end{maplelatex}

\end{maplegroup}
\begin{maplegroup}
\begin{mapleinput}
\mapleinline{active}{1d}{F:=map(normal,clicollect(makeF(dim)));}{%
}
\end{mapleinput}

\mapleresult
\begin{maplelatex}
\mapleinline{inert}{2d}{F :=
-(A23*g11+A12*g13-A13*g12)*e23/(-g33*g22*g11+g22*g13^2-2*g13*g23*g12+g
33*g12^2+g23^2*g11)+(-A23*g13+A13*g23-A12*g33)*e12/(-g33*g22*g11+g22*g
13^2-2*g13*g23*g12+g33*g12^2+g23^2*g11)+(A23*g12-A13*g22+A12*g23)*e13/
(-g33*g22*g11+g22*g13^2-2*g13*g23*g12+g33*g12^2+g23^2*g11);}{%
\maplemultiline{
F :=  - {\displaystyle \frac {(\mathit{A23}\,\mathit{g11} + 
\mathit{A12}\,\mathit{g13} - \mathit{A13}\,\mathit{g12})\,
\mathit{e23}}{\mathrm{\%1}}}  + {\displaystyle \frac {( - 
\mathit{A23}\,\mathit{g13} + \mathit{A13}\,\mathit{g23} - 
\mathit{A12}\,\mathit{g33})\,\mathit{e12}}{\mathrm{\%1}}}  \\
\mbox{} + {\displaystyle \frac {(\mathit{A23}\,\mathit{g12} - 
\mathit{A13}\,\mathit{g22} + \mathit{A12}\,\mathit{g23})\,
\mathit{e13}}{\mathrm{\%1}}}  \\
\mathrm{\%1} :=  - \mathit{g33}\,\mathit{g22}\,\mathit{g11} + 
\mathit{g22}\,\mathit{g13}^{2} - 2\,\mathit{g13}\,\mathit{g23}\,
\mathit{g12} + \mathit{g33}\,\mathit{g12}^{2} + \mathit{g23}^{2}
\,\mathit{g11} }
}
\end{maplelatex}

\end{maplegroup}
\begin{maplegroup}
\begin{mapleinput}
\mapleinline{active}{1d}{N:=floor(dim/2):F1:=wexp(F,N):F2:=wexp(-F,N):}{%
}
\end{mapleinput}
\end{maplegroup}
\noindent
We will now define two general elements $u$ and $v$ in $\cl(B)$ by decomposing them over a Grassmann basis (provided by a procedure {\tt cbasis}). Coefficients in these two expansions are assumed to be real or complex.\footnote{For technical reasons, in Maple these coefficients cannot be defined as $u_i$ and $v_i;$ that is why we use here $uu_i$ and $vv_i.$}
\begin{maplegroup}
\begin{mapleinput}
\mapleinline{active}{1d}{cbasis(dim);}{%
}
\end{mapleinput}

\mapleresult
\begin{maplelatex}
\mapleinline{inert}{2d}{[Id, e1, e2, e3, e12, e13, e23, e123];}{%
\[
[\mathit{Id}, \,\mathit{e1}, \,\mathit{e2}, \,\mathit{e3}, \,
\mathit{e12}, \,\mathit{e13}, \,\mathit{e23}, \,\mathit{e123}]
\]
}
\end{maplelatex}

\end{maplegroup}
\begin{maplegroup}
\begin{mapleinput}
\mapleinline{active}{1d}{u:=add(uu[k]*cbasis(dim)[k],k=1..2^dim);}{%
}
\end{mapleinput}

\mapleresult
\begin{maplelatex}
\mapleinline{inert}{2d}{u :=
uu[1]*Id+uu[2]*e1+uu[3]*e2+uu[4]*e3+uu[5]*e12+uu[6]*e13+uu[7]*e23+uu[8
]*e123;}{%
\[
u := {\mathit{uu}_{1}}\,\mathit{Id} + {\mathit{uu}_{2}}\,\mathit{
e1} + {\mathit{uu}_{3}}\,\mathit{e2} + {\mathit{uu}_{4}}\,
\mathit{e3} + {\mathit{uu}_{5}}\,\mathit{e12} + {\mathit{uu}_{6}}
\,\mathit{e13} + {\mathit{uu}_{7}}\,\mathit{e23} + {\mathit{uu}_{
8}}\,\mathit{e123}
\]
}
\end{maplelatex}

\end{maplegroup}
\begin{maplegroup}
\begin{mapleinput}
\mapleinline{active}{1d}{v:=add(vv[k]*cbasis(dim)[k],k=1..2^dim);}{%
}
\end{mapleinput}

\mapleresult
\begin{maplelatex}
\mapleinline{inert}{2d}{v :=
vv[1]*Id+vv[2]*e1+vv[3]*e2+vv[4]*e3+vv[5]*e12+vv[6]*e13+vv[7]*e23+vv[8
]*e123;}{%
\[
v := {\mathit{vv}_{1}}\,\mathit{Id} + {\mathit{vv}_{2}}\,\mathit{
e1} + {\mathit{vv}_{3}}\,\mathit{e2} + {\mathit{vv}_{4}}\,
\mathit{e3} + {\mathit{vv}_{5}}\,\mathit{e12} + {\mathit{vv}_{6}}
\,\mathit{e13} + {\mathit{vv}_{7}}\,\mathit{e23} + {\mathit{vv}_{
8}}\,\mathit{e123}
\]
}
\end{maplelatex}
\end{maplegroup}
\noindent
The Clifford product of $u$ and $v$ in $\cl(B)$ is then collected and assigned to a constant $res_1$ which we won't display due to its length.
\begin{maplegroup}
\begin{mapleinput}
\mapleinline{active}{1d}{res1:=clicollect(cmul(u,v)):}{%
}
\end{mapleinput}
\end{maplegroup}
\noindent
As before, we finish by computing the right hand side of (\ref{eq:Bprod}). By assigning it to $res_2,$ we can then easily find that $res_1 - res_2 = 0$ as expected.
\begin{maplegroup}
\begin{mapleinput}
\mapleinline{active}{1d}{
res2:=clicollect(RCg(cmulg(RCg(u,F1),RCg(v,F1)),F2)):res1-res2;}{%
}
\end{mapleinput}

\mapleresult
\begin{maplelatex}
\mapleinline{inert}{2d}{0;}{%
\[
0
\]
}
\end{maplelatex}

\end{maplegroup}

\section{Hecke algebra computations}
\label{sec:hecke}


In \cite{AblFauser99a} it was shown that the symmetric group $S_n$ and its group deformation, the Hecke algebra $H_{\BF}(n,q)$ could be constructed as a subalgebra of a Clifford algebra $\cl(B)$ for a suitably chosen $q$-dependent non-symmetric bilinear form $B.$ $q$-Young operators were constructed as Clifford idempotents and the Hecke algebra representations in ideals generated by these idempotents were computed. Appropriate $q$-Young diagrams and tableaux representing symmetrizers, antisymmetrizers, and operators of mixed symmetries were realized inside the Hecke algebra, while the ordinary case of the symmetric group was obtained in the limit $q \rightarrow 1.$

The Hecke algebra is the generalization of the group algebra of the symmetric group $S_n$ by adding the requirement that transpositions of adjacent elements $i,i+1$ are no longer involutions. Following \cite{AblFauser99a} we set $t_i^2 = (1-q) t_i + q$ which reduces to $s_i^2= 1$ in the limit $q\rightarrow 1.$ The defining relations of the Hecke algebra will be given according to Bourbaki \cite{Bourbaki}. Let $\{{\bf 1}, t_1, \ldots, t_n\}$ be a set of generators which fulfill these relations:
\begin{eqnarray}
t_i^2 &=& (1-q) t_i + q,                            \label{eq: t1} \\
t_i t_j &=& t_j t_i,   \quad \vert i-j \vert \ge 2, \label{eq: t2} \\
t_i t_{i+1} t_i &=& t_{i+1} t_i t_{i+1},            \label{eq: t3} 
\end{eqnarray}
then their algebraic span is the Hecke algebra $H_{\BF}(n,q).$ The algebra morphism $\rho$ which maps the Hecke algebra into the even part of an appropriate Clifford algebra was found in \cite{Fauser-hecke}. In particular, $\rho(t_i) = b_i := \e_i \w \e_{i+n},\,i=1,\ldots,n,$ where $\e_1,\ldots,\e_{2n}$ are the generators of the Clifford algebra 
$\cl(B,V),\, V=\span\{\e_i\},$ with the following non-symmetric bilinear form~$B:$
\begin{maplegroup}
\begin{mapleinput}
\mapleinline{active}{1d}{dim:=8:n:=dim/2:eval(makealiases(dim,'ordered')):B:=defB(dim);}{%
}
\end{mapleinput}
\mapleresult
\begin{maplelatex}
\mapleinline{inert}{2d}{B := matrix([[0, 0, 0, 0, q, -1-q, 1, 1], [0, 0, 0, 0, -1-q, q, -1-q,
1], [0, 0, 0, 0, 1, -1-q, q, -1-q], [0, 0, 0, 0, 1, 1, -1-q, q], [1,
1, -1, -1, 0, 0, 0, 0], [q, 1, 1, -1, 0, 0, 0, 0], [-1, q, 1, 1, 0, 0,
0, 0], [-1, -1, q, 1, 0, 0, 0, 0]]);}{%
\[
B :=  \left[ 
{\begin{array}{cccrcccc}
0 & 0 & 0 & 0 & q &  - 1 - q & 1 & 1 \\
0 & 0 & 0 & 0 &  - 1 - q & q &  - 1 - q & 1 \\
0 & 0 & 0 & 0 & 1 &  - 1 - q & q &  - 1 - q \\
0 & 0 & 0 & 0 & 1 & 1 &  - 1 - q & q \\
1 & 1 & -1 & -1 & 0 & 0 & 0 & 0 \\
q & 1 & 1 & -1 & 0 & 0 & 0 & 0 \\
-1 & q & 1 & 1 & 0 & 0 & 0 & 0 \\
-1 & -1 & q & 1 & 0 & 0 & 0 & 0
\end{array}}
 \right] 
\]
}
\end{maplelatex}
\end{maplegroup}
\noindent
Then, the form of $B$ guarantees that the following relations hold:
\begin{eqnarray}
b_i^2 &=& (1-q) b_i + q,                            \label{eq: b1} \\
b_i b_j &=& b_j b_i,   \quad \vert i-j \vert \ge 2, \label{eq: b2} \\
b_i b_{i+1} b_i &=& b_{i+1} b_i b_{i+1}.            \label{eq: b3} 
\end{eqnarray}
Following \cite{Fauser-hecke}, we define the Hecke generators $b_i$ as {\it balanced\/} basis Grassmann monomials of order~$2:$
\begin{maplegroup}
\begin{mapleinput}
\mapleinline{active}{1d}{for i from 1 to n do b.i:=(e.i) &w (e.(n+i)) od:}{%
}
\end{mapleinput}
\end{maplegroup}
\noindent
Using procedure {\tt cliexpand} we can expand these generators in terms of the unevaluated Clifford product which is denoted in Maple as $\&C:$
\begin{maplegroup}
\begin{mapleinput}
\mapleinline{active}{1d}{seq(cat('b',i)=cliexpand(b.i),i=1..4);}{%
}
\end{mapleinput}
\mapleresult
\begin{maplelatex}
\mapleinline{inert}{2d}{
b1 = `&C`(e1,e5)-q*Id, b2 = `&C`(e2,e6)-q*Id, b3 = `&C`(e3,e7)-q*Id, b4 = `&C`(e4,e8)-q*Id;}{%
\maplemultiline{
\mathit{b1}=(\mathit{e1}\,\mathrm{\&C}\,\mathit{e5}) - q\,\mathit{Id}, \,
\mathit{b2}=(\mathit{e2}\,\mathrm{\&C}\,\mathit{e6}) - q\,\mathit{Id}, \,
\mathit{b3}=(\mathit{e3}\,\mathrm{\&C}\,\mathit{e7}) - q\,\mathit{Id}, \,
\mathit{b4}=(\mathit{e4}\,\mathrm{\&C}\,\mathit{e8}) - q\,\mathit{Id}. }
}
\end{maplelatex}
\end{maplegroup}
\noindent
Checking if the Hecke generators satisfy the defining relations can be done as follows:
\begin{maplegroup}
\begin{mapleinput}
\mapleinline{active}{1d}{map(evalb,[seq(simplify(CS(cmul(b.i,b.i)-(1-q)*b.i-q))=0,i=1..n)]);}{
}
\end{mapleinput}

\mapleresult
\begin{maplelatex}
\mapleinline{inert}{2d}{[true, true, true, true];}{%
\[
[\mathit{true}, \,\mathit{true}, \,\mathit{true}, \,\mathit{true}
]
\]
}
\end{maplelatex}

\end{maplegroup}
\begin{maplegroup}
\begin{mapleinput}
\mapleinline{active}{1d}{map(evalb,[seq(seq(CS(cmul(b.i,b.j)-cmul(b.j,b.i))=0,i=j+2..n),j=1..n)]);}{
}
\end{mapleinput}

\mapleresult
\begin{maplelatex}
\mapleinline{inert}{2d}{[true, true, true];}{%
\[
[\mathit{true}, \,\mathit{true}, \,\mathit{true}]
\]
}
\end{maplelatex}

\end{maplegroup}
\begin{maplegroup}
\begin{mapleinput}
\mapleinline{active}{1d}{map(evalb@factor@normal,[seq(CS(cmul(cmul(b.i,cat(b,i+1)),b.i)-
cmul(cmul(cat(b,i+1),b.i),cat(b,i+1)))=0,i=1..n-1)]);}{%
}
\end{mapleinput}

\mapleresult
\begin{maplelatex}
\mapleinline{inert}{2d}{[true, true, true];}{%
\[
[\mathit{true}, \,\mathit{true}, \,\mathit{true}]
\]
}
\end{maplelatex}
\end{maplegroup}
\noindent
Let's define the remaining basis elements of the Hecke image $\rho(H_{\BF}(4,q))$ in $\cl(B):$
\begin{maplegroup}
\begin{mapleinput}
\mapleinline{active}{1d}{b12:=CS(cmul(b1,b2)):b21:=CS(cmul(b2,b1)):b23:=CS(cmul(b2,b3)):
b32:=CS(cmul(b3,b2)):b34:=CS(cmul(b3,b4)):b43:=CS(cmul(b4,b3)):
b121:=CS(cmul(b1,b21)):b123:=CS(cmul(b12,b3)):b234:=CS(cmul(b23,b4)):
b321:=CS(cmul(b32,b1)):b432:=CS(cmul(b43,b2)):b4321:=CS(cmul(b432,b1)):
b1234:=CS(cmul(b123,b4)):}{%
}
\end{mapleinput}
\end{maplegroup}
\noindent
Thus, $\rho$ defines a homomorphism from the Hecke algebra $H_{\BF}(4,q)$ into the Clifford algebra $\cl(B).$ It was shown in \cite{Fauser-hecke} that $\rho$ is not injective for 
$n \geq 4,$ and that its kernel contains all Young diagrams which are not $L$-shaped.

\subsection{Hecke Algebra $H_{\BF}(2,q)$}
\label{sec:hecke2}

We begin with the Hecke algebra $H_{\BF}(2,q)$ generated by $\{Id,b_1\}$ which reduces to $S_2$ in the limit $q\rightarrow 1.$ We have thus only one $q$-transposition $b_1$ from which we can calculate a $q$-symmetrizer $R(12)$ and a $q$-antisymmetrizer $C(12).$ One of the features of the construction presented in \cite{AblFauser99a} is that $R(12)$ and $C(12)$ are related by the reversion in $\cl(B):$
\begin{maplegroup}
\begin{mapleinput}
\mapleinline{active}{1d}{R.12:=CS(q*Id+b1);   #S_2 symmetrizer
C.12:=CS(Id-b1);     #S_2 antisymmetrizer
evalb(reversion(C.12)=R.12);}{%
}
\end{mapleinput}

\mapleresult
\begin{maplelatex}
\mapleinline{inert}{2d}{R12 := e15+q*Id;c12 := Id-e15;}{%
\[
\mathit{R12} := \mathit{e15} + q\,\mathit{Id},\quad \mathit{c12} := \mathit{Id} - \mathit{e15}
\]
}
\end{maplelatex}
\begin{maplelatex}
\mapleinline{inert}{2d}{true;}{%
\[
\mathit{true}
\]
}
\end{maplelatex}
\end{maplegroup}
\noindent
Observe, that  $R(12)$ is the reversion of $C(12),$ that is, $R(12) = \revside{C(12)}$ where
$\tilde{~}$ is the reversion in the Clifford algebra $\cl_{1,1}.$  Operators $C(12)$ and
$R(12)$ are almost idempotent and they annihilate each other:
\begin{maplegroup}
\begin{mapleinput}
\mapleinline{active}{1d}{factor(CS(R.12 &c R.12));}{%
}
\end{mapleinput}

\mapleresult
\begin{maplelatex}
\mapleinline{inert}{2d}{(q+1)*(q*Id+e15);}{%
\[
(q + 1)\,(q\,\mathit{Id} + \mathit{e15})
\]
}
\end{maplelatex}

\end{maplegroup}
\begin{maplegroup}
\begin{mapleinput}
\mapleinline{active}{1d}{factor(CS(C.12 &c C.12));}{%
}
\end{mapleinput}

\mapleresult
\begin{maplelatex}
\mapleinline{inert}{2d}{(q+1)*(Id-e15);}{%
\[
(q + 1)\,(\mathit{Id} - \mathit{e15})
\]
}
\end{maplelatex}

\end{maplegroup}
\begin{maplegroup}
\begin{mapleinput}
\mapleinline{active}{1d}{R.12 &c C.12, C.12 &c R.12;}{%
}
\end{mapleinput}

\mapleresult
\begin{maplelatex}
\mapleinline{inert}{2d}{0, 0;}{%
\[
0, \,0
\]
}
\end{maplelatex}
\end{maplegroup}
\noindent
Upon normalization, the symmetrizer $R(12)$ becomes the $q$-Young operator $Y^{(2)}_{1,2}$ while the antisymmetrizer $C(12)$ becomes $Y^{(11)}_{1,2}.$ The following computation verifies that  $Y^{(2)}_{1,2}$ and $Y^{(11)}_{1,2}$ are mutually annihilating idempotents adding up to the identity element.
\begin{maplegroup}
\begin{mapleinput}
\mapleinline{active}{1d}{Y2.12 := R.12/(1+q);Y11.12:= C.12/(1+q);}{%
}
\end{mapleinput}

\mapleresult
\begin{maplelatex}
\mapleinline{inert}{2d}{Y212 := (q*Id+e15)/(q+1);Y1112 := (Id-e15)/(q+1);}{%
\[
\mathit{Y212} := {\displaystyle \frac {q\,\mathit{Id} + \mathit{e15}}{q + 1}},\quad 
\mathit{Y1112} := {\displaystyle \frac {\mathit{Id} - \mathit{e15}}{q + 1}} 
\]
}
\end{maplelatex}
\end{maplegroup}
\begin{maplegroup}
\begin{mapleinput}
\mapleinline{active}{1d}{ybas:=[Y2.12,Y11.12]:f:=(i,j)->CS(cmul(ybas[i],ybas[j]));
YM:=evalm(linalg[matrix](2,2,(i,j)->f(i,j)));CS(Y11.12+Y2.12);}{%
}
\end{mapleinput}

\mapleresult
\begin{maplelatex}
\mapleinline{inert}{2d}{f := proc (i, j) options operator, arrow; CS(climul(ybas[i],ybas[j]))
end;}{%
\[
f := (i, \,j)\rightarrow \mathrm{CS}(\mathrm{climul}({\mathit{
ybas}_{i}}, \,{\mathit{ybas}_{j}}))
\]
}
\end{maplelatex}

\begin{maplelatex}
\mapleinline{inert}{2d}{YM := matrix([[q*Id/(q+1)+e15/(q+1), 0], [0,
Id/(q+1)-e15/(q+1)]]);}{%
\[
\mathit{YM} :=  \left[ 
{\begin{array}{cc}
{\displaystyle \frac {q\,\mathit{Id}}{q + 1}}  + {\displaystyle 
\frac {\mathit{e15}}{q + 1}}  & 0 \\ [2ex]
0 & {\displaystyle \frac {\mathit{Id}}{q + 1}}  - {\displaystyle 
\frac {\mathit{e15}}{q + 1}} 
\end{array}}
 \right] 
\]
}
\end{maplelatex}

\begin{maplelatex}
\mapleinline{inert}{2d}{Id;}{%
\[
\mathit{Id}
\]
}
\end{maplelatex}
\end{maplegroup}
In our construction, the Hecke algebra $H_{\BF}(2,q)$ is a subalgebra of the even part  $\cl^+_{1,1}$ of $\cl_{1,1}$ and it is generated by $\{\e_1,\e_5\}.$ That
is, the bilinear form $G$ on the vector space spanned by $\{\e_1,\e_5\}$~is
\begin{maplegroup}
\begin{mapleinput}
\mapleinline{active}{1d}{G:=matrix(2,2,[B[1,1],B[1,5],B[5,1],B[5,5]]);}{%
}
\end{mapleinput}

\mapleresult
\begin{maplelatex}
\mapleinline{inert}{2d}{G := matrix([[0, q], [1, 0]]);}{%
\[
G :=  \left[ 
{\begin{array}{rc}
0 & q \\
1 & 0
\end{array}}
 \right]. 
\]
}
\end{maplelatex}
\end{maplegroup}
\noindent
After symmetrization, diagonalization, and the limit $q\rightarrow 1,$ $G$ becomes
$\diag(1,-1).$ Therefore, due to the isomorphism $\cl(B)\simeq \cl(g)$ (as associative algebras) we can view $H_{\BF}(2,q)$ as the subalgebra of $\cl^+_{1,1}.$ This embedding implies that any idempotent of the Hecke algebra $H_{\BF}(2,q)$ must be an even Clifford element.  The following computation shows that the Young operators found above are the only two nontrivial mutually annihilating idempotents in $H_{\BF}(2,q).$
\begin{maplegroup}
\begin{mapleinput}
\mapleinline{active}{1d}{xx:=CS(a*Id+b*b1);}{%
}
\end{mapleinput}

\mapleresult
\begin{maplelatex}
\mapleinline{inert}{2d}{xx := a*Id+b*e15;}{%
\[
\mathit{xx} := a\,\mathit{Id} + b\,\mathit{e15}
\]
}
\end{maplelatex}

\end{maplegroup}
\begin{maplegroup}
\begin{mapleinput}
\mapleinline{active}{1d}{sol:=clisolve2(CS(cmul(xx,xx)-xx),[a,b]);}{%
}
\end{mapleinput}

\mapleresult
\begin{maplelatex}
\mapleinline{inert}{2d}{sol := [{b = 0, a = 0}, {a = 1, b = 0}, {b = -1/(q+1), a =
1/(q+1)}, {b = 1/(q+1), a = q/(q+1)}];}{%
\[
\mathit{sol} := [{b=0, \,a=0}, \,{a=1, \,b=0}, \,{b= - 
{\displaystyle \frac {1}{q + 1}} , \,a={\displaystyle \frac {1}{q
 + 1}} }, \,{b={\displaystyle \frac {1}{q + 1}} , \,a=
{\displaystyle \frac {q}{q + 1}} }]
\]
}
\end{maplelatex}

\end{maplegroup}
\begin{maplegroup}
\begin{mapleinput}
\mapleinline{active}{1d}{f.1:=normal(subs(sol[4],xx));f.2:=normal(subs(sol[3],xx));}{%
}
\end{mapleinput}

\mapleresult
\begin{maplelatex}
\mapleinline{inert}{2d}{f1 := (q*Id+e15)/(q+1);f2 := (Id-e15)/(q+1);}{%
\[
\mathit{f1} := {\displaystyle \frac {q\,\mathit{Id} + \mathit{e15}}{q + 1}},\quad
\mathit{f2} := {\displaystyle \frac {\mathit{Id} - \mathit{e15}}{q + 1}}. 
\]
}
\end{maplelatex}
\end{maplegroup}
\noindent
In this case, the Young operators happen to be the two even primitive idempotents $f_1,f_2$ in  $\cl^{+}_{1,1}.$ Notice, that in the case when $q=1,$ the idempotents $f_1$ and $f_2$ reduce to the well-known primitive idempotents
\begin{maplegroup}
\begin{mapleinput}
\mapleinline{active}{1d}{f11:=subs(q=1,f1);}{%
}
\end{mapleinput}

\mapleresult
\begin{maplelatex}
\mapleinline{inert}{2d}{f11 := 1/2*Id+1/2*e15;}{%
\[
\mathit{f11} := {\displaystyle \frac {1}{2}} \,\mathit{Id} + 
{\displaystyle \frac {1}{2}} \,\mathit{e15}
\]
}
\end{maplelatex}

\end{maplegroup}
\begin{maplegroup}
\begin{mapleinput}
\mapleinline{active}{1d}{f22:=subs(q=1,f2);}{%
}
\end{mapleinput}

\mapleresult
\begin{maplelatex}
\mapleinline{inert}{2d}{f22 := 1/2*Id-1/2*e15;}{%
\[
\mathit{f22} := {\displaystyle \frac {1}{2}} \,\mathit{Id} - 
{\displaystyle \frac {1}{2}} \,\mathit{e15}
\]
}
\end{maplelatex}
\end{maplegroup}
\noindent
in the Clifford algebra $\cl_{1,1}.$

\subsection{Hecke algebra $H_{\BF}(3,q)$}
\label{sec:hecke3}

In \cite{AblFauser99a} Young operators related to various symmetries were constructed as idempotent elements in the Hecke algebra $H_{\BF}(3,q)$ embedded into the even subalgebra $\cl^+_{2,2}$ of $\cl_{2,2}.$  As one of the main features of this construction, the Young operators corresponding to conjugate Young tableaux in the sense of McDonald \cite{Macdonald} were related through the reversion in the Clifford algebra. The goal was to find four Young operators known to exist from the general theory of the Hecke algebras for $n=3.$ The four $q$-Young operators found had one parameter and generalized the four Young operators of $S_3$ described in Hamermesh \cite{Hamermesh} on p. 245. One of them was a full symmetrizer, one was a full antisymmetrizer and two were of mixed symmetry. 

The construction began with finding the most general element $X$ in $H_{\BF}(3,q)$ such that\footnote{See formulas (24), (25) and (26) in \cite{AblFauser99a}.}:
\begin{eqnarray}
X+ \revtop{X}  &=& Id \label{eq:X1}\\
X^2            &=& X  \label{eq:X2}\\
X \revtop{X}   &=& 0  \label{eq:X3}
\end{eqnarray}
In the first step, the most general element in the Hecke algebra $H_{\BF}(3,q)$ was found that  satisfied~(\ref{eq:X1})\footnote{Procedure {\tt clisolve2} is capable of solving the equation $X+\revtop{X}-\Id = 0$ for the free parameters $K_1,\ldots,K_6$ appearing in $X.$ However, when used repeatedly, it randomly selects free parameters in the solution.}:
\begin{maplegroup}
\begin{mapleinput}
\mapleinline{active}{1d}{bset:=[Id,b1,b2,b12,b21,b121]:X:=bexpand(add(K[i]*bset[i],i=1..6)):
sollist[1]:=clisolve2(X+reversion(X)-Id,[seq(K[i],i=1..6)]):
X:=bexpand(subs(sollist[1],X));}{%
}
\end{mapleinput}
\mapleresult
\begin{maplelatex}
\mapleinline{inert}{2d}{X :=
(-1/2*K[6]*q+1/2*K[3]*q+1/2*K[6]*q^2+1/2*K[2]*q+1/2-1/2*K[3]-1/2*K[2])
*Id+K[2]*b1+K[3]*b2+K[4]*b12+(-K[6]+K[6]*q-K[4])*b21+K[6]*b121;}{%
\maplemultiline{
\mathit{X} := ( - {\displaystyle \frac {1}{2}} \,{K_{6}}\,q + 
{\displaystyle \frac {1}{2}} \,{K_{3}}\,q + {\displaystyle 
\frac {1}{2}} \,{K_{6}}\,q^{2} + {\displaystyle \frac {1}{2}} \,{
K_{2}}\,q + {\displaystyle \frac {1}{2}}  - {\displaystyle 
\frac {1}{2}} \,{K_{3}} - {\displaystyle \frac {1}{2}} \,{K_{2}})
\,\mathit{Id} + {K_{2}}\,\mathit{b1} + {K_{3}}\,\mathit{b2} + {K
_{4}}\,\mathit{b12} \\
\mbox{} + ( - {K_{6}} + {K_{6}}\,q - {K_{4}})\,\mathit{b21} + {K
_{6}}\,\mathit{b121} }
}
\end{maplelatex}
\end{maplegroup}
\noindent
In the above, procedure {\tt bexpand} expands elements in the Hecke algebra, which are normally expressed in the Grassmann basis of $\cl_{2,2},$ in terms of the Hecke basis $\{\Id,b_1,b_2,b_{12},b_{21},b_{121}\}.$

Thus, the solution to (\ref{eq:X1}) gives an element that belongs to a family parameterized by four real or complex parameters. In order to simplify Maple output, we define two aliases and then we substitute $X$ found above into the second equation (\ref{eq:X2}) .
\begin{maplegroup}
\begin{mapleinput}
\mapleinline{active}{1d}{alias(alpha=RootOf((1+q)*_Z^2+(-q^2*K[4]+K[4]+q*K[2]-1+K[2])*_Z+K[4]*K[2]+
K[2]^2-K[4]-K[2]+q*K[4]-q^2*K[4]^2-q*K[4]^2-q^2*K[2]*K[4]+q*K[2]^2)):}{%
}
\mapleinline{active}{1d}{alias(kappa=RootOf((1+q)*_Z^2+(-q^2*K[4]+K[4]+q*K[2]+1+K[2])*_Z+K[4]*
K[2]+K[2]^2+K[4]+K[2]-q*K[4]-q^2*K[4]^2-q*K[4]^2-q^2*K[2]*K[4]+q*K[2]^2));}{%
}
\end{mapleinput}
\mapleresult
\begin{maplelatex}
\mapleinline{inert}{2d}{I, alpha, kappa;}{%
\[
I, \,\alpha , \,\kappa 
\]
}
\end{maplelatex}
\end{maplegroup}
\begin{maplegroup}
\begin{mapleinput}
\mapleinline{active}{1d}{sollistxx1:=clisolve2(cmul(X,X)-X,[seq(K[i],i=1..6)]):
xxlist:=[seq(bexpand(subs(sollistxx1[i],X)),i=1..nops(sollistxx1))]:
for i from 1 to nops(xxlist) do r.i:=bexpand(xxlist[i]) od:}{%
}
\end{mapleinput}
\end{maplegroup}
\noindent
The six representatives $r_1,\ldots,r_6$ defined in this last Maple command are the same as the elements displayed in \cite{AblFauser99a} after formula (30). Using the command {\tt findbasis} we can determine that the set $\{r_i\},\,i=1,\ldots,6,$ is of rank $4:$
\begin{maplegroup}
\begin{mapleinput}
\mapleinline{active}{1d}{nops(findbasis([r.(1..6)]));}{%
}
\end{mapleinput}
\mapleresult
\begin{maplelatex}
\mapleinline{inert}{2d}{4;}{%
\[
4
\]
}
\end{maplelatex}
\end{maplegroup}
\noindent
Since the representatives are related by the reversion, we seek four linearly independent elements. For example, $r_1,\,r_2,\,r_3$ and $r_5$ are linearly independent. From now on we assign them to $f_1,f_2,f_3,f_4.$
\begin{maplegroup}
\begin{mapleinput}
\mapleinline{active}{1d}{nops(findbasis([xxlist[1],xxlist[2],xxlist[3],xxlist[5]]));}{%
}
\end{mapleinput}
\mapleresult
\begin{maplelatex}
\mapleinline{inert}{2d}{4;}{%
\[
4
\]
}
\end{maplelatex}

\end{maplegroup}
\begin{maplegroup}
\begin{mapleinput}
\mapleinline{active}{1d}{f1:=r.1:f2:=r.2:f3:=r.3:f4:=r.5:}{%
}
\end{mapleinput}
\end{maplegroup}
\noindent
We can easily verify that the four elements $\{f_1, f_2, f_3, f_4\}$ satisfy equation (\ref{eq:X1}):
\begin{maplegroup}
\begin{mapleinput}
\mapleinline{active}{1d}{CS(f.1+reversion(f.1)-Id),CS(f.2+reversion(f.2)-Id),
CS(f.3+reversion(f.3)-Id),CS(f.4+reversion(f.4)-Id); }{%
}
\end{mapleinput}
\mapleresult
\begin{maplelatex}
\mapleinline{inert}{2d}{0, 0, 0, 0;}{%
\[
0, \,0, \,0, \,0
\]
}
\end{maplelatex}
\end{maplegroup}
\noindent
and equation (\ref{eq:X2}):
\begin{maplegroup}
\begin{mapleinput}
\mapleinline{active}{1d}{simplify(CS(f.1 &c f.1 - f.1)),simplify(CS(f.2 &c f.2 - f.2)),
simplify(CS(f.3 &c f.3 - f.3)),simplify(CS(f.4 &c f.4 - f.4));}{%
}
\end{mapleinput}
\mapleresult
\begin{maplelatex}
\mapleinline{inert}{2d}{0, 0, 0, 0;}{%
\[
0, \,0, \,0, \,0
\]
}
\end{maplelatex}

\end{maplegroup}
Each of the four non-primitive idempotents $f_i$ generates a three-dimensional one-sided ideal in  $H_{\BF}(3,q).$ In order to split such ideal into a one-dimensional space and a two-dimensional space, one had to find a way of splitting at least one of these idempotents into a sum of two mutually annihilating idempotents. It was observed in \cite{AblFauser99a} that $f_1$ contained a full symmetrizer. Since the full symmetrizer $\Y{(3)}{1,2,3}$ was defined there as
\begin{eqnarray}
\Y{(3)}{1,2,3}   &:=& \frac{q^3\Id + q^2 b_1 + q^2 b_2 + q b_{12} + q b_{21} + b_{121}}
                           {(1+q+q^2)(1+q)}, 
\label{eq:Y3}
\end{eqnarray}
the full antisymmetrizer $\Y{(111)}{1,2,3}$ was, by construction, the reversion of the symmetrizer, that is, $\Y{(111)}{1,2,3} := \revside{\Y{(3)}{1,2,3}},$
\begin{eqnarray}
\Y{(111)}{1,2,3} &:=& \frac{\Id - b_1 - b_2 + b_{12} + b_{21} - b_{121}}
                           {(1+q+q^2)(1+q)}. 
\label{eq:Y111}
\end{eqnarray}
Therefore, by subtracting the full antisymmetrizer $\Y{(111)}{1,2,3}$ from $f_1,$ the first Young operator $\Y{(21)}{1,3,2}$ of a mixed type was found. Then, the second Young operator $\Y{(21)}{1,2,3}$ of a mixed type was computed by applying the reversion (conjugate) to $\Y{(21)}{1,3,2},$ that is, $\Y{(21)}{1,2,3} := \reversion{\Y{(21)}{1,3,2}}.$ Namely\footnote{Due to lengthy displays, we refer Reader to formulas (33) and (34) in \cite{AblFauser99a}.},
\begin{maplegroup}
\begin{mapleinput}
\mapleinline{active}{1d}{Y21.132:=bexpand(f.1-Y111.123):Y21.123:=bexpand(reversion(Y21.132)):}{%
}
\end{mapleinput}
\end{maplegroup}
\noindent
Furthermore, $\Y{(111)}{1,2,3}$ annihilates $\Y{(21)}{1,3,2}$ when multiplied from both sides. This is a reflection of the fact that the left ideal (the representation space) in the Hecke algebra $H_{\BF}(3,q)$ generated by $f_1$ decomposes into a direct sum of one-dimensional left ideal and a two-dimensional left ideal.
\begin{maplegroup}
\begin{mapleinput}
\mapleinline{active}{1d}{simplify(Y.111.123+Y.21.132-f1);}{%
}
\end{mapleinput}

\mapleresult
\begin{maplelatex}
\mapleinline{inert}{2d}{0;}{%
\[
0
\]
}
\end{maplelatex}

\end{maplegroup}
\begin{maplegroup}
\begin{mapleinput}
\mapleinline{active}{1d}{simplify(Y.21.123+Y.3.123-reversion(f1));}{%
}
\end{mapleinput}

\mapleresult
\begin{maplelatex}
\mapleinline{inert}{2d}{0;}{%
\[
0
\]
}
\end{maplelatex}

\end{maplegroup}
\begin{maplegroup}
\begin{mapleinput}
\mapleinline{active}{1d}{CS(Y.111.123 &c Y.21.132),CS(Y.21.132 &c Y.111.123);}{%
}
\end{mapleinput}

\mapleresult
\begin{maplelatex}
\mapleinline{inert}{2d}{0, 0;}{%
\[
0, \,0
\]
}
\end{maplelatex}

\end{maplegroup}
\begin{maplegroup}
\begin{mapleinput}
\mapleinline{active}{1d}{CS(Y.21.132 &c Y.3.123),CS(Y.3.123 &c Y.21.123);}{%
}
\end{mapleinput}

\mapleresult
\begin{maplelatex}
\mapleinline{inert}{2d}{0, 0;}{%
\[
0, \,0
\]
}
\end{maplelatex}

\end{maplegroup}
\begin{maplegroup}
\begin{mapleinput}
\mapleinline{active}{1d}{
CS(Y.111.123 &c Y.111.123 - Y.111.123),CS(Y.21.132 &c Y.21.132 - Y.21.132);
}{%
}
\end{mapleinput}

\mapleresult
\begin{maplelatex}
\mapleinline{inert}{2d}{0, 0;}{%
\[
0, \,0
\]
}
\end{maplelatex}
\end{maplegroup}
\noindent
In addition to the symmetrizer $\Y{(3)}{1,2,3}$ and the antisymmetrizer $\Y{(111)}{1,2,3}$ displayed in (\ref{eq:Y3}) and (\ref{eq:Y111}), we have also two Young operators $\Y{(21)}{1,2,3}$ and $\Y{(21)}{1,3,2}$ of mixed symmetry:
\begin{maplegroup}
\begin{mapleinput}
\mapleinline{active}{1d}{'Y21.123'=bexpand(Y21.123);}{%
}
\end{mapleinput}

\mapleresult
\begin{maplelatex}
\mapleinline{inert}{2d}{Y21.123 =
q*Id/(q+1+q^2)+(K[4]*q^3+2*K[4]*q^2-q^2+2*K[4]*q+K[4])*b1/((q+1+q^2)*(
1+q))-(K[4]*q^3+2*K[4]*q^2+q+2*K[4]*q+K[4])*q*b2/((q+1+q^2)*(1+q))-(K[
4]*q^3+2*K[4]*q^2+q+2*K[4]*q+K[4])*b12/(q^3+2*q^2+2*q+1)+(q^5*K[4]+q^4
*K[4]+q^3+K[4]*q^3-q^2+K[4]*q^2+K[4]*q+K[4]-1)*b21/(q*(q^3+2*q^2+2*q+1
))+(q^4*K[4]+K[4]*q^3+q^2-K[4]*q-K[4]+1)*b121/((q+1+q^2)*q*(1+q));}{%
\maplemultiline{
\mathit{Y21}.123={\displaystyle \frac {q\,\mathit{Id}}{q + 1 + q
^{2}}}  + {\displaystyle \frac {({K_{4}}\,q^{3} + 2\,{K_{4}}\,q^{
2} - q^{2} + 2\,{K_{4}}\,q + {K_{4}})\,\mathit{b1}}{(q + 1 + q^{2
})\,(1 + q)}}  \\
\mbox{} - {\displaystyle \frac {({K_{4}}\,q^{3} + 2\,{K_{4}}\,q^{
2} + q + 2\,{K_{4}}\,q + {K_{4}})\,q\,\mathit{b2}}{(q + 1 + q^{2}
)\,(1 + q)}}  \\
\mbox{} - {\displaystyle \frac {({K_{4}}\,q^{3} + 2\,{K_{4}}\,q^{
2} + q + 2\,{K_{4}}\,q + {K_{4}})\,\mathit{b12}}{q^{3} + 2\,q^{2}
 + 2\,q + 1}}  \\
\mbox{} + {\displaystyle \frac {(q^{5}\,{K_{4}} + q^{4}\,{K_{4}}
 + q^{3} + {K_{4}}\,q^{3} - q^{2} + {K_{4}}\,q^{2} + {K_{4}}\,q
 + {K_{4}} - 1)\,\mathit{b21}}{q\,(q^{3} + 2\,q^{2} + 2\,q + 1)}
}  \\
\mbox{} + {\displaystyle \frac {(q^{4}\,{K_{4}} + {K_{4}}\,q^{3}
 + q^{2} - {K_{4}}\,q - {K_{4}} + 1)\,\mathit{b121}}{(q + 1 + q^{
2})\,q\,(1 + q)}}  }
}
\end{maplelatex}
\end{maplegroup}

\begin{maplegroup}
\begin{mapleinput}
\mapleinline{active}{1d}{'Y21.132'=bexpand(Y21.132);}{%
}
\end{mapleinput}

\mapleresult
\begin{maplelatex}
\mapleinline{inert}{2d}{Y21.132 =
q*Id/(q+1+q^2)-(K[4]*q^3+2*K[4]*q^2+2*K[4]*q-1+K[4])*b1/((q+1+q^2)*(1+
q))+(q^4*K[4]+2*K[4]*q^3+2*K[4]*q^2+K[4]*q+1)*b2/((q+1+q^2)*(1+q))+(K[
4]*q^3+2*K[4]*q^2+2*K[4]*q-1+K[4])*b12/(q^3+2*q^2+2*q+1)-(q^5*K[4]+q^4
*K[4]+K[4]*q^3+q^3+K[4]*q^2+K[4]*q+q+K[4]-1)*b21/(q*(q^3+2*q^2+2*q+1))
-(q^4*K[4]+K[4]*q^3+q^2-K[4]*q-K[4]+1)*b121/((q+1+q^2)*q*(1+q));}{%
\maplemultiline{
\mathit{Y21}.132={\displaystyle \frac {q\,\mathit{Id}}{q + 1 + q
^{2}}}  - {\displaystyle \frac {({K_{4}}\,q^{3} + 2\,{K_{4}}\,q^{
2} + 2\,{K_{4}}\,q - 1 + {K_{4}})\,\mathit{b1}}{(q + 1 + q^{2})\,
(1 + q)}}  \\
\mbox{} + {\displaystyle \frac {(q^{4}\,{K_{4}} + 2\,{K_{4}}\,q^{
3} + 2\,{K_{4}}\,q^{2} + {K_{4}}\,q + 1)\,\mathit{b2}}{(q + 1 + q
^{2})\,(1 + q)}}  \\
\mbox{} + {\displaystyle \frac {({K_{4}}\,q^{3} + 2\,{K_{4}}\,q^{
2} + 2\,{K_{4}}\,q - 1 + {K_{4}})\,\mathit{b12}}{q^{3} + 2\,q^{2}
 + 2\,q + 1}}  \\
\mbox{} - {\displaystyle \frac {(q^{5}\,{K_{4}} + q^{4}\,{K_{4}}
 + {K_{4}}\,q^{3} + q^{3} + {K_{4}}\,q^{2} + {K_{4}}\,q + q + {K
_{4}} - 1)\,\mathit{b21}}{q\,(q^{3} + 2\,q^{2} + 2\,q + 1)}}  \\
\mbox{} - {\displaystyle \frac {(q^{4}\,{K_{4}} + {K_{4}}\,q^{3}
 + q^{2} - {K_{4}}\,q - {K_{4}} + 1)\,\mathit{b121}}{(q + 1 + q^{
2})\,q\,(1 + q)}}  }
}
\end{maplelatex}
\end{maplegroup}
\noindent
In order to represent our Young operator $\Y{(21)}{1,3,2}$ as a product of a row-symmetrizer
$R(13)$ and a column-antisymmetrizer $C(12)=f_1,$ we use $f_1$ defined above and compute $R(13)$ from the equation
\beq
\Y{(21)}{132}=R(13)f_1.
\label{eq:R13}
\eeq
In order to solve the above equation in Maple for $R(13),$ we need to find an element $Y$ in the Hecke algebra which would not only satisfy equation (\ref{eq:R13}) but also such that $Y+\revside{Y}=\Id.$ This is because we want the column antisymmetrizer $C(12)$ to remain related to $R(12)$ through the reversion. Notice also that we are justified in defining $C(12)$ as equal to $f_1$ (modulo a
normalizing factor) because $f_1$ generalizes the antisymmetrizer $C(12)$ from $S_2$ to $S_3:$
\begin{maplegroup}
\begin{mapleinput}
\mapleinline{active}{1d}{bexpand(c12/(1+q));}{%
}
\end{mapleinput}

\mapleresult
\begin{maplelatex}
\mapleinline{inert}{2d}{Id/(1+q)-b1/(1+q);}{%
\[
{\displaystyle \frac {\mathit{Id}}{1 + q}}  - {\displaystyle 
\frac {\mathit{b1}}{1 + q}} 
\]
}
\end{maplelatex}

\end{maplegroup}
\begin{maplegroup}
\begin{mapleinput}
\mapleinline{active}{1d}{bexpand(f1);}{%
}
\end{mapleinput}

\mapleresult
\begin{maplelatex}
\mapleinline{inert}{2d}{Id/(1+q)-K[4]*b1+K[4]*q*b2+K[4]*b12-(K[4]*q^3+q+K[4]-1)*b21/(q*(1+q))
-(-K[4]+K[4]*q^2+1)*b121/(q*(1+q));}{%
\[
{\displaystyle \frac {\mathit{Id}}{1 + q}}  - {K_{4}}\,\mathit{b1
} + {K_{4}}\,q\,\mathit{b2} + {K_{4}}\,\mathit{b12} - 
{\displaystyle \frac {({K_{4}}\,q^{3} + q + {K_{4}} - 1)\,
\mathit{b21}}{q\,(1 + q)}}  - {\displaystyle \frac {( - {K_{4}}
 + {K_{4}}\,q^{2} + 1)\,\mathit{b121}}{q\,(1 + q)}} 
\]
}
\end{maplelatex}
\end{maplegroup}
\noindent
That is, $f_1$ is seen to contain $C(12)/(1+q)$ if we replace $K_4$ with $1/(1+q).$ Thus, we first express $Y$ in the Hecke basis contained in the list {\tt bset} below and then we make sure that $Y+\revside{Y}=\Id:$
\begin{maplegroup}
\begin{mapleinput}
\mapleinline{active}{1d}{bset:=[Id,b1,b2,b12,b21,b121]:
Y:=bexpand(add(P[i]*bset[i],i=1..6)):
Y:=bexpand(op(clisolve2(Y+reversion(Y)-Id,Y)));}{%
}
\end{mapleinput}

\mapleresult
\begin{maplelatex}
\mapleinline{inert}{2d}{Y :=
(1/2*P[3]*q-1/2*P[6]*q+1/2*P[6]*q^2+1/2*P[2]*q-1/2*P[2]+1/2-1/2*P[3])*
Id+P[2]*b1+P[3]*b2+(-P[5]-P[6]+P[6]*q)*b12+P[5]*b21+P[6]*b121;}{%
\maplemultiline{
\mathit{Y} := ({\displaystyle \frac {1}{2}} \,{P_{3}}\,q - 
{\displaystyle \frac {1}{2}} \,{P_{6}}\,q + {\displaystyle 
\frac {1}{2}} \,{P_{6}}\,q^{2} + {\displaystyle \frac {1}{2}} \,{
P_{2}}\,q - {\displaystyle \frac {1}{2}} \,{P_{2}} + 
{\displaystyle \frac {1}{2}}  - {\displaystyle \frac {1}{2}} \,{P
_{3}})\,\mathit{Id} + {P_{2}}\,\mathit{b1} + {P_{3}}\,\mathit{b2}
 \\
\mbox{} + ( - {P_{5}} - {P_{6}} + {P_{6}}\,q)\,\mathit{b12} + {P
_{5}}\,\mathit{b21} + {P_{6}}\,\mathit{b121} }
}
\end{maplelatex}
\end{maplegroup}
\noindent
Next we require that the above found element $Y$ satisfies equation (\ref{eq:R13}):
\begin{maplegroup}
\begin{mapleinput}
\mapleinline{active}{1d}{R.13:=bexpand(op(clisolve2(Y21.132-cmul(Y,f_1),Y)));}{%
}
\end{mapleinput}

\mapleresult
\begin{maplelatex}
\mapleinline{inert}{2d}{R13 :=
q*Id/(1+q)-(-q^2+q^2*P[3]+P[3]*q-1+P[3])*b1/((q+1+q^2)*q)+P[3]*b2+(q^2
*P[3]+P[3]*q+P[3]-1)*b12/(q*(q+1+q^2))-(q^5*P[3]+q^4*P[3]+q^3*P[3]-q^2
+q^2*P[3]+P[3]*q-1+P[3])*b21/((1+q)*q^2*(q+1+q^2))-(q^4*P[3]+q^3*P[3]-
P[3]*q+1-P[3])*b121/(q^2*(q^3+2*q^2+2*q+1));}{%
\maplemultiline{
\mathit{R13} := {\displaystyle \frac {q\,\mathit{Id}}{1 + q}}  - 
{\displaystyle \frac {( - q^{2} + q^{2}\,{P_{3}} + {P_{3}}\,q - 1
 + {P_{3}})\,\mathit{b1}}{(q + 1 + q^{2})\,q}}  + {P_{3}}\,
\mathit{b2} + {\displaystyle \frac {(q^{2}\,{P_{3}} + {P_{3}}\,q
 + {P_{3}} - 1)\,\mathit{b12}}{q\,(q + 1 + q^{2})}}  \\
\mbox{} - {\displaystyle \frac {(q^{5}\,{P_{3}} + q^{4}\,{P_{3}}
 + q^{3}\,{P_{3}} - q^{2} + q^{2}\,{P_{3}} + {P_{3}}\,q - 1 + {P
_{3}})\,\mathit{b21}}{(1 + q)\,q^{2}\,(q + 1 + q^{2})}}  \\
\mbox{} - {\displaystyle \frac {(q^{4}\,{P_{3}} + q^{3}\,{P_{3}}
 - {P_{3}}\,q + 1 - {P_{3}})\,\mathit{b121}}{q^{2}\,(q^{3} + 2\,q
^{2} + 2\,q + 1)}}  }
}
\end{maplelatex}

\end{maplegroup}
\noindent
We verify that $R(13)$ is an idempotent:
\begin{maplegroup}
\begin{mapleinput}
\mapleinline{active}{1d}{CS(R13 &c R13-R13);}{%
}
\end{mapleinput}

\mapleresult
\begin{maplelatex}
\mapleinline{inert}{2d}{0;}{%
\[
0
\]
}
\end{maplelatex}
\end{maplegroup}
\noindent
It was pointed out in \cite{AblFauser99a} that when the Clifford product $R(13)f_1$ is computed, the free parameter $P_3$ disappears:
\begin{maplegroup}
\begin{mapleinput}
\mapleinline{active}{1d}{bexpand(cmul(R13,f1));}{%
}
\end{mapleinput}

\mapleresult
\begin{maplelatex}
\mapleinline{inert}{2d}{q*Id/(q+1+q^2)-(K[4]*q^3+2*K[4]*q^2+2*K[4]*q-1+K[4])*b1/(q^3+2*q^2+2*
q+1)+(q^4*K[4]+2*K[4]*q^3+2*K[4]*q^2+K[4]*q+1)*b2/(q^3+2*q^2+2*q+1)+(K
[4]*q^3+2*K[4]*q^2+2*K[4]*q-1+K[4])*b12/(q^3+2*q^2+2*q+1)-(q^5*K[4]+q^
4*K[4]+K[4]*q^3+q^3+K[4]*q^2+K[4]*q+q+K[4]-1)*b21/(q*(q^3+2*q^2+2*q+1)
)-(q^4*K[4]+K[4]*q^3+q^2-K[4]*q-K[4]+1)*b121/((q+1+q^2)*q*(1+q));}{%
\maplemultiline{
{\displaystyle \frac {q\,\mathit{Id}}{q + 1 + q^{2}}}  - 
{\displaystyle \frac {({K_{4}}\,q^{3} + 2\,{K_{4}}\,q^{2} + 2\,{K
_{4}}\,q - 1 + {K_{4}})\,\mathit{b1}}{\mathrm{\%1}}}  \\
\mbox{} + {\displaystyle \frac {(q^{4}\,{K_{4}} + 2\,{K_{4}}\,q^{
3} + 2\,{K_{4}}\,q^{2} + {K_{4}}\,q + 1)\,\mathit{b2}}{\mathrm{
\%1}}}  \\
\mbox{} + {\displaystyle \frac {({K_{4}}\,q^{3} + 2\,{K_{4}}\,q^{
2} + 2\,{K_{4}}\,q - 1 + {K_{4}})\,\mathit{b12}}{\mathrm{\%1}}} 
 \\
\mbox{} - {\displaystyle \frac {(q^{5}\,{K_{4}} + q^{4}\,{K_{4}}
 + {K_{4}}\,q^{3} + q^{3} + {K_{4}}\,q^{2} + {K_{4}}\,q + q + {K
_{4}} - 1)\,\mathit{b21}}{q\,\mathrm{\%1}}}  \\
\mbox{} - {\displaystyle \frac {(q^{4}\,{K_{4}} + {K_{4}}\,q^{3}
 + q^{2} - {K_{4}}\,q - {K_{4}} + 1)\,\mathit{b121}}{(q + 1 + q^{
2})\,q\,(1 + q)}}  \\
\mathrm{\%1} := q^{3} + 2\,q^{2} + 2\,q + 1 }
}
\end{maplelatex}

\end{maplegroup}
\noindent
which gives $\Y{(21)}{1,3,2}$ computed above:
\begin{maplegroup}
\begin{mapleinput}
\mapleinline{active}{1d}{CS(\%-Y21.132);}{%
}
\end{mapleinput}

\mapleresult
\begin{maplelatex}
\mapleinline{inert}{2d}{0;}{%
\[
0
\]
}
\end{maplelatex}
\end{maplegroup}
In order to construct the representation spaces from the four Young operators, one needs to find at least one {\it Garnir\/} $\G{(\lambda_i)}{i,j}$ element in the Hecke algebra \cite{AblFauser99a,KingWybourne}. All Garnir elements can be seen to act as row or column cycles in Young tableaux thereby generating non-standard tableaux which correspond to the basis vectors of the representation space. A Garnir element has the following defining properties:
\begin{eqnarray}
\Y{(21)}{1,2,3}\,  \G{(21)}{1,1}   &   = & 0, \label{eq:G1} \\[1ex]
\G{(21)}{1,1}\,  \Y{(21)}{1,2,3}   & \neq& 0. \label{eq:G2} 
\end{eqnarray}
The reason for requiring (\ref{eq:G2}) is that we want $\G{(21)}{1,1} \Y{(21)}{1,2,3}$ to be a second basis element in the left Hecke ideal generated by $\Y{(21)}{1,2,3}.$ In order to solve (\ref{eq:G1}) and (\ref{eq:G2}), we begin by assigning to $X$ a general element of the Hecke algebra, that is, $X$ is a linear combination of the Hecke basis elements with some undefined coefficients $K_i,i=1,\ldots,6:$
\begin{maplegroup}
\begin{mapleinput}
\mapleinline{active}{1d}{X:=bexpand(K[1]*Id+K[2]*b1+K[3]*b2+K[4]*b12+K[5]*b21+K[6]*b121);}{%
}
\end{mapleinput}

\mapleresult
\begin{maplelatex}
\mapleinline{inert}{2d}{X := K[1]*Id+K[2]*b1+K[3]*b2+K[4]*b12+K[5]*b21+K[6]*b121;}{%
\[
X={K_{1}}\,\mathit{Id} + {K_{2}}\,\mathit{b1} + {K_{3}}\,\mathit{
b2} + {K_{4}}\,\mathit{b12} + {K_{5}}\,\mathit{b21} + {K_{6}}\,
\mathit{b121}
\]
}
\end{maplelatex}

\end{maplegroup}
\noindent
We use {\tt clisolve2} to solve equation (\ref{eq:G1}). All three solutions returned by Maple are assigned to a list {\tt sol}.
\begin{maplegroup}
\begin{mapleinput}
\mapleinline{active}{1d}{sol:=clisolve2(cmul(Y21.123,X),[seq(K[i],i=1..6)]):nops(sol);}{%
}
\end{mapleinput}

\mapleresult
\begin{maplelatex}
\mapleinline{inert}{2d}{3;}{%
\[
3
\]
}
\end{maplelatex}

\end{maplegroup}
\begin{maplegroup}
\begin{mapleinput}
\mapleinline{active}{1d}{for i from 1 to nops(sol) do X.i:=bexpand(subs(sol[i],X)) od;}{%
}
\end{mapleinput}

\mapleresult
\begin{maplelatex}
\mapleinline{inert}{2d}{X1 :=
(K[2]*q-K[6]*q+K[4])*Id+K[2]*b1+K[3]*b2+K[4]*b12+(K[6]*q^2+K[4]*q^2-K[
3]*q-K[6]*q-K[4]*q+K[2]+K[4])*b21/q+K[6]*b121;}{%
\maplemultiline{
\mathit{X1} := ({K_{2}}\,q - {K_{6}}\,q + {K_{4}})\,\mathit{Id}
 + {K_{2}}\,\mathit{b1} + {K_{3}}\,\mathit{b2} + {K_{4}}\,
\mathit{b12} \\
\mbox{} + {\displaystyle \frac {({K_{6}}\,q^{2} + {K_{4}}\,q^{2}
 - {K_{3}}\,q - {K_{6}}\,q - {K_{4}}\,q + {K_{2}} + {K_{4}})\,
\mathit{b21}}{q}}  + {K_{6}}\,\mathit{b121} }
}
\end{maplelatex}

\begin{maplelatex}
\mapleinline{inert}{2d}{X2 :=
K[1]*Id+(K[6]*q^3+q^2*K[1]+K[6]*q^2+q*K[1]+1)*b1/(q^2*(1+q))+K[3]*b2-b
12/(q*(1+q))+(q^5*K[6]-q^4*K[3]-q^3*K[3]-q^3+K[6]*q^2+q^2*K[1]+q^2+q*K
[1]-q+1)*b21/(q^3*(1+q))+K[6]*b121;}{%
\maplemultiline{
\mathit{X2} := {K_{1}}\,\mathit{Id} + {\displaystyle \frac {({K
_{6}}\,q^{3} + q^{2}\,{K_{1}} + {K_{6}}\,q^{2} + q\,{K_{1}} + 1)
\,\mathit{b1}}{q^{2}\,(1 + q)}}  + {K_{3}}\,\mathit{b2} - 
{\displaystyle \frac {\mathit{b12}}{q\,(1 + q)}}  \\
\mbox{} + {\displaystyle \frac {(q^{5}\,{K_{6}} - q^{4}\,{K_{3}}
 - q^{3}\,{K_{3}} - q^{3} + {K_{6}}\,q^{2} + q^{2}\,{K_{1}} + q^{
2} + q\,{K_{1}} - q + 1)\,\mathit{b21}}{q^{3}\,(1 + q)}}  + {K_{6
}}\,\mathit{b121} }
}
\end{maplelatex}

\begin{maplelatex}
\mapleinline{inert}{2d}{X3 :=
K[1]*Id+(q-1+K[6]*q+2*K[6]*q^2+2*K[6]*q^3+K[1]+2*q*K[1]+q^4*K[6]+2*q^2
*K[1]+q^3*K[1])*b1/(q*(q^3+2*q^2+2*q+1))+K[3]*b2-(q-1)*b12/(q^3+2*q^2+
2*q+1)+(q^6*K[6]+q^5*K[6]-q^5*K[3]-2*q^4*K[3]-q^4+q^4*K[6]+2*q^3+q^3*K
[1]+K[6]*q^3-2*q^3*K[3]+K[6]*q^2+2*q^2*K[1]-q^2*K[3]-2*q^2+K[6]*q+2*q*
K[1]+2*q+K[1]-1)*b21/(q^2*(q^3+2*q^2+2*q+1))+K[6]*b121;}{%
\maplemultiline{
\mathit{X3} := {K_{1}}\,\mathit{Id} \\
\mbox{} + {\displaystyle \frac {(q - 1 + {K_{6}}\,q + 2\,{K_{6}}
\,q^{2} + 2\,{K_{6}}\,q^{3} + {K_{1}} + 2\,q\,{K_{1}} + q^{4}\,{K
_{6}} + 2\,q^{2}\,{K_{1}} + q^{3}\,{K_{1}})\,\mathit{b1}}{q\,(q^{
3} + 2\,q^{2} + 2\,q + 1)}}  \\
\mbox{} + {K_{3}}\,\mathit{b2} - {\displaystyle \frac {(q - 1)\,
\mathit{b12}}{q^{3} + 2\,q^{2} + 2\,q + 1}}  + ((q^{6}\,{K_{6}}
 + q^{5}\,{K_{6}} - q^{5}\,{K_{3}} - 2\,q^{4}\,{K_{3}} - q^{4} + 
q^{4}\,{K_{6}} \\
\mbox{} + 2\,q^{3} + q^{3}\,{K_{1}} + {K_{6}}\,q^{3} - 2\,q^{3}\,
{K_{3}} + {K_{6}}\,q^{2} + 2\,q^{2}\,{K_{1}} - q^{2}\,{K_{3}} - 2
\,q^{2} + {K_{6}}\,q + 2\,q\,{K_{1}} \\
\mbox{} + 2\,q + {K_{1}} - 1)\mathit{b21})/(q^{2}\,(q^{3} + 2\,q
^{2} + 2\,q + 1))\mbox{} + {K_{6}}\,\mathit{b121} }
}
\end{maplelatex}

\end{maplegroup}
\noindent
One way to find out if the three solutions returned by Maple are really linearly independent is to use a procedure {\tt findbasis} which from the given list of Clifford polynomials extracts all polynomials which are linearly independent.
\begin{maplegroup}
\begin{mapleinput}
\mapleinline{active}{1d}{nops(findbasis([X.1,X.2,X.3]));}{%
}
\end{mapleinput}

\mapleresult
\begin{maplelatex}
\mapleinline{inert}{2d}{3;}{%
\[
3
\]
}
\end{maplelatex}
\end{maplegroup}
\noindent
The three solutions $X_1,X_2,X_3$ are therefore linearly independent. Another way to verify that fact would be to try to find three coefficients $c_1,c_2,c_3,$ not all equal to zero, that would satisfy the following linear combination:
\beq
c_1 \, X_1 + c_2 \, X_2 + c_3 \, X_3=0.
\label{eq:Xcomb}
\eeq
Thus, we can use again the procedure {\tt clisolve2} and try to solve equation (\ref{eq:Xcomb}) as follows:
\begin{maplegroup}
\begin{mapleinput}
\mapleinline{active}{1d}{clisolve2(c[1]*X.1+c[2]*X.2+c[3]*X.3,[c[1],c[2],c[3]]);}{%
}
\end{mapleinput}

\mapleresult
\begin{maplelatex}
\mapleinline{inert}{2d}{[\{c[1] = 0, c[2] = 0, c[3] = 0\}];}{%
\[
[\{{\mathit{c}_{1}}=0, \,{\mathit{c}_{2}}=0, \,{\mathit{c}_{3}}=0\}]
\]
}
\end{maplelatex}

\end{maplegroup}
\noindent
As expected, all coefficients are zero. Thus, $X_1,X_2,X_3$ are three linearly independent solutions of (\ref{eq:G1}), that is, $\Y{(21)}{1,2,3}X_i = 0$ for $i=1,2,3.$

We proceed now to verify whether the solutions we have just obtained satisfy also equation (\ref{eq:G2}). In particular, we will try to see if there any non-zero values of the parameters $K_2,K_4,K_5,K_6$ in $X_1$ so that $X_1 \Y{(21)}{1,2,3}$ would be $0.$ We will again use the procedure {\tt clisolve2}.
\begin{maplegroup}
\begin{mapleinput}
\mapleinline{active}{1d}{var1:=select(type,indets(X1),indexed);}{%
}
\end{mapleinput}

\mapleresult
\begin{maplelatex}
\mapleinline{inert}{2d}{var1 := \{K[2], K[4], K[5], K[6]\};}{%
\[
\mathit{var1} := \{{K_{2}}, \,{K_{4}}, \,{K_{5}}, \,{K_{6}}\}
\]
}
\end{maplelatex}

\end{maplegroup}
\begin{maplegroup}
\begin{mapleinput}
\mapleinline{active}{1d}{clisolve2(X.1 &c Y.21.123,vars1);}{%
}
\end{mapleinput}

\mapleresult
\begin{maplelatex}
\mapleinline{inert}{2d}{[ ];}{%
\[
[\,]
\]
}
\end{maplelatex}

\end{maplegroup}
\noindent
As expected, Maple returns an empty solution set.

\subsection{Automorphism $\alphaq$ and the Garnir elements in the Hecke algebra}
\label{sec:Garnir}

In \cite{AblFauser99a}, an automorphism $\alphaq$ was introduced in the Hecke algebra via the  formula $(48).$ $\alphaq$ replaces the reversion $\tilde{~}$ and gives the inverse of the basis element $b_K$ for any multi-index $K.$ The automorphism $\alphaq$ is then extended to the whole Hecke algebra by the following definition:
\begin{equation}
\begin{array}{rcl}
\alphaq(b_{i_1} \ldots b_{i_s}) 
                      & = & 
        (\frac{-1}{q})^s (b_{i_1} \cdots b_{i_s})\tilde{~}       \\[1.5ex]
                      & = & 
        (\frac{-1}{q})^s (\tilde{b_{i_s}} \cdots \tilde{b_{i_1}})\\[1.5ex]
                      & = &  
\alphaq(b_{i_s}) \cdots \alphaq(b_{i_1}).
\end{array}
\label{eq:alphaq}
\end{equation}
\noindent
In {\sc CLIFFORD}, the $\alphaq$ automorphism has been programmed as a Maple procedure {\tt alpha2}.  For example, when $s=2,$ verification of (\ref{eq:alphaq}) can be done as follows:\footnote{Procedure {\tt alpha2} is part of a package {\sc suppl}. Rather than displaying and comparing long expressions, it is often convenient to use Maple's built-in Boolean procedure {\tt evalb} which returns {\tt true} or {\tt false} when the two expressions are equal or not.}
\begin{maplegroup}
\begin{mapleinput}
\mapleinline{active}{1d}{evalb(bexpand(alpha2(b1 &c b2))=bexpand(alpha2(b2) &c alpha2(b1)));
evalb(bexpand(alpha2(b2 &c b1))=bexpand(alpha2(b1) &c alpha2(b2)));}{%
}
\end{mapleinput}

\mapleresult
\begin{maplelatex}
\mapleinline{inert}{2d}{true;true;}{%
\[
\mathit{true},\;\mathit{true},
\]
}
\end{maplelatex}
\end{maplegroup}
\noindent
The following is a verification that indeed 
$\alphaq(b_i)=b_i^{-1}=\frac{(q-1)}{q}\Id + \frac{b_i}{q},\,i=1,2.$\footnote{In {\sc CLIFFORD}, the symbolic inverse of any element can be found using a procedure {\tt cinv}.}

\begin{maplegroup}
\begin{mapleinput}
\mapleinline{active}{1d}{bexpand(alpha2(b1)),evalb(bexpand(cinv(b1)) = bexpand(alpha2(b1)));
bexpand(alpha2(b2)),evalb(bexpand(cinv(b2)) = bexpand(alpha2(b2)));}{%
}
\end{mapleinput}

\mapleresult
\begin{maplelatex}
\mapleinline{inert}{2d}{(q-1)*Id/q+b1/q, true;}{%
\[
{\displaystyle \frac {(q - 1)\,\mathit{Id}}{q}}  + 
{\displaystyle \frac {\mathit{b1}}{q}} , \,\mathit{true}
\]
}
\end{maplelatex}

\begin{maplelatex}
\mapleinline{inert}{2d}{(q-1)*Id/q+b2/q, true;}{%
\[
{\displaystyle \frac {(q - 1)\,\mathit{Id}}{q}}  + 
{\displaystyle \frac {\mathit{b2}}{q}} , \,\mathit{true}
\]
}
\end{maplelatex}
\end{maplegroup}
\noindent
On the other hand, (\ref{eq:alphaq}) implies that $\alphaq(b_{12}) = (\frac{-1}{q})^2 \, \tilde{b_{12}}=b_{12}^{-1}:$
\begin{maplegroup}
\begin{mapleinput}
\mapleinline{active}{1d}{bexpand(alpha2(b12));
map(normal,bexpand((-1/q)^2*(reversion(b12))));
evalb(bexpand(cinv(b12))=bexpand(alpha2(b12)));
}{%
}

\end{mapleinput}

\mapleresult
\begin{maplelatex}
\mapleinline{inert}{2d}{(1-2*q+q^2)*Id/(q^2)+(q-1)*b1/(q^2)+(q-1)*b2/(q^2)+b21/(q^2);}{%
\[
{\displaystyle \frac {(1 - 2\,q + q^{2})\,\mathit{Id}}{q^{2}}} 
 + {\displaystyle \frac {(q - 1)\,\mathit{b1}}{q^{2}}}  + 
{\displaystyle \frac {(q - 1)\,\mathit{b2}}{q^{2}}}  + 
{\displaystyle \frac {\mathit{b21}}{q^{2}}} 
\]
}
\end{maplelatex}

\begin{maplelatex}
\mapleinline{inert}{2d}{(1-2*q+q^2)*Id/(q^2)+(q-1)*b1/(q^2)+(q-1)*b2/(q^2)+b21/(q^2);}{%
\[
{\displaystyle \frac {(1 - 2\,q + q^{2})\,\mathit{Id}}{q^{2}}} 
 + {\displaystyle \frac {(q - 1)\,\mathit{b1}}{q^{2}}}  + 
{\displaystyle \frac {(q - 1)\,\mathit{b2}}{q^{2}}}  + 
{\displaystyle \frac {\mathit{b21}}{q^{2}}} 
\]
}
\end{maplelatex}

\begin{maplelatex}
\mapleinline{inert}{2d}{true;}{%
\[
\mathit{true}
\]
}
\end{maplelatex}

\end{maplegroup}
\noindent
and similarly for $b_{21}.$ For completeness we only show that
$\alphaq(b_{121}) = (\frac{-1}{q})^3 \, \tilde{b_{121}}=b_{121}^{-1}:$ 
\begin{maplegroup}
\begin{mapleinput}
\mapleinline{active}{1d}{bexpand(alpha2(b121));
map(normal,bexpand((-1/q)^3*(reversion(b121))));
evalb(bexpand(cinv(b121))=bexpand(alpha2(b121)));}{%
}
\end{mapleinput}

\mapleresult
\begin{maplelatex}
\mapleinline{inert}{2d}{(-2*q^2+q^3+2*q-1)*Id/(q^3)+(1-2*q+q^2)*b1/(q^3)+(1-2*q+q^2)*b2/(q^3)
+(q-1)*b12/(q^3)+(q-1)*b21/(q^3)+b121/(q^3);}{%
\maplemultiline{
{\displaystyle \frac {( - 2\,q^{2} + q^{3} + 2\,q - 1)\,\mathit{
Id}}{q^{3}}}  + {\displaystyle \frac {(1 - 2\,q + q^{2})\,
\mathit{b1}}{q^{3}}}  + {\displaystyle \frac {(1 - 2\,q + q^{2})
\,\mathit{b2}}{q^{3}}}  + {\displaystyle \frac {(q - 1)\,\mathit{
b12}}{q^{3}}}  \\
\mbox{} + {\displaystyle \frac {(q - 1)\,\mathit{b21}}{q^{3}}} 
 + {\displaystyle \frac {\mathit{b121}}{q^{3}}}  }
}
\end{maplelatex}

\begin{maplelatex}
\mapleinline{inert}{2d}{(-2*q^2+q^3+2*q-1)*Id/(q^3)+(1-2*q+q^2)*b1/(q^3)+(1-2*q+q^2)*b2/(q^3)
+(q-1)*b12/(q^3)+(q-1)*b21/(q^3)+b121/(q^3);}{%
\maplemultiline{
{\displaystyle \frac {( - 2\,q^{2} + q^{3} + 2\,q - 1)\,\mathit{
Id}}{q^{3}}}  + {\displaystyle \frac {(1 - 2\,q + q^{2})\,
\mathit{b1}}{q^{3}}}  + {\displaystyle \frac {(1 - 2\,q + q^{2})
\,\mathit{b2}}{q^{3}}}  + {\displaystyle \frac {(q - 1)\,\mathit{
b12}}{q^{3}}}  \\
\mbox{} + {\displaystyle \frac {(q - 1)\,\mathit{b21}}{q^{3}}} 
 + {\displaystyle \frac {\mathit{b121}}{q^{3}}}  }
}
\end{maplelatex}

\begin{maplelatex}
\mapleinline{inert}{2d}{true;}{%
\[
\mathit{true}
\]
}
\end{maplelatex}

\end{maplegroup}
\noindent
Thus, we have verified that $\alphaq(b_{K}) = (\frac{-1}{q})^{|K|} \, \tilde{b_{K}}=b_{K}^{-1}$ for a multi-index $K$ of length $|K|.$ However, this property does not extend to non-homogeneous (non-versor like) elements of the Hecke algebra. Let {\tt hecke} be a Maple variable representing an arbitrary element in $H_\BF(3,q)$ expanded in the Hecke basis $\{\Id,b_1,b_2,b_{12},b_{21},b_{121}\}:$
\begin{maplegroup}
\begin{mapleinput}
\mapleinline{active}{1d}{
hecke:=h[1]*Id + h[2]*'b1' + h[3]*'b2' + h[4]*'b12' + h[5]*'b21' + h[6]*'b121';}{%
}
\end{mapleinput}

\mapleresult
\begin{maplelatex}
\mapleinline{inert}{2d}{hecke := h[1]*Id+h[2]*b1+h[3]*b2+h[4]*b12+h[5]*b21+h[6]*b121;}{%
\[
\mathit{hecke} := {h_{1}}\,\mathit{Id} + {h_{2}}\,\mathit{b1} + {
h_{3}}\,\mathit{b2} + {h_{4}}\,\mathit{b12} + {h_{5}}\,\mathit{
b21} + {h_{6}}\,\mathit{b121}
\]
}
\end{maplelatex}
\end{maplegroup}
\noindent
Then the action of $\alphaq$ on {\tt hecke} can be written as follows:
\begin{maplegroup}
\begin{mapleinput}
\mapleinline{active}{1d}{alpha2hecke:=bexpand(alpha2(hecke));}{%
}
\end{mapleinput}

\mapleresult
\begin{maplelatex}
\mapleinline{inert}{2d}{alpha2hecke :=
H1*Id/(q^3)+H2*b1/(q^3)+H3*b2/(q^3)+H4*b12/(q^3)+H5*b21/(q^3)+H6*b121/
(q^3);}{%
\[
\mathit{alpha2hecke} := {\displaystyle \frac {\mathit{H1}\,
\mathit{Id}}{q^{3}}}  + {\displaystyle \frac {\mathit{H2}\,
\mathit{b1}}{q^{3}}}  + {\displaystyle \frac {\mathit{H3}\,
\mathit{b2}}{q^{3}}}  + {\displaystyle \frac {\mathit{H4}\,
\mathit{b12}}{q^{3}}}  + {\displaystyle \frac {\mathit{H5}\,
\mathit{b21}}{q^{3}}}  + {\displaystyle \frac {\mathit{H6}\,
\mathit{b121}}{q^{3}}} 
\]
}
\end{maplelatex}
\end{maplegroup}
\noindent
where the expressions $H_i,i=1,\ldots,6,$ are $q$-polynomials with coefficients expressed in terms of $h_i,i=1,\ldots,6.$\footnote{In order to get the display for $\alphaq(hecke)$ shown above expressions $H_i,i=1,\ldots,6,$ have been defined in Maple as {\tt aliases}. See Appendix 2 for their definitions.} With a little experimentation it can be easily verified that, for example, $\alphaq$ does not give the inverse of the element $\Id+b_1:$
\begin{maplegroup}
\begin{mapleinput}
\mapleinline{active}{1d}{h[1],h[2],h[3],h[4],h[5],h[6]:=1,1,0,0,0,0:bexpand(hecke);}{%
}
\end{mapleinput}

\mapleresult
\begin{maplelatex}
\mapleinline{inert}{2d}{Id+b1;}{%
\[
\mathit{Id} + \mathit{b1}
\]
}
\end{maplelatex}

\end{maplegroup}
\begin{maplegroup}
\begin{mapleinput}
\mapleinline{active}{1d}{bexpand(CS(cinv(hecke))),bexpand(alpha2(hecke));}{%
}
\end{mapleinput}

\mapleresult
\begin{maplelatex}
\mapleinline{inert}{2d}{1/2*(q-2)*Id/(q-1)+1/2*b1/(q-1), (2*q-1)*Id/q+b1/q;}{%
\[
{\displaystyle \frac {1}{2}} \,{\displaystyle \frac {(q - 2)\,\mathit{Id}}{q - 1}}  + {\displaystyle \frac {1}{2}} \,{\displaystyle \frac {\mathit{b1}}{q - 1}}, \quad
{\displaystyle \frac {(2\,q - 1)\,\mathit{Id}}{q}}  + {\displaystyle \frac {\mathit{b1}}{q}}. 
\]
}
\end{maplelatex}
\end{maplegroup}

We return now to the problem of finding the Garnir element $\G{(21)}{1,1}$ that would satisfy equations (\ref{eq:G1}) and (\ref{eq:G2}). Recall that equation (\ref{eq:G1}) had three linearly independent solutions $X_1,X_2,X_3.$ Notice, that $X_1$ does not annihilate the Young operator $\Y{(21)}{1,2,3}$ and $\alphaq(X_1)$ does not annihilate the Young operator $\Y{(21)}{1,3,2}$ no matter what values are assigned to the parameters $K_2,K_4,K_5,K_6:$
\begin{maplegroup}
\begin{mapleinput}
\mapleinline{active}{1d}{
clisolve2(X.1 &c Y.21.123,vars1),clisolve2(alpha2(X.1) &c Y.21.132,vars1);}{%
}
\end{mapleinput}

\mapleresult
\begin{maplelatex}
\mapleinline{inert}{2d}{[], [];}{%
\[
[\,], \,[\,]
\]
}
\end{maplelatex}
\end{maplegroup}
\noindent
Furthermore, the following shows that the elements $X_1 \Y{(21)}{1,2,3}$ and
$\alphaq(X_1)\Y{(21)}{1,3,2}$ are linearly independent:
\begin{maplegroup}
\begin{mapleinput}
\mapleinline{active}{1d}{nops(findbasis([X.1 &c Y.21.123,alpha2(X.1) &c Y.21.132]));}{%
}
\end{mapleinput}

\mapleresult
\begin{maplelatex}
\mapleinline{inert}{2d}{2;}{%
\[
2
\]
}
\end{maplelatex}
\end{maplegroup}
\noindent
Since the six elements 
$$
\{\Y{(3)}{1,2,3},\Y{(21)}{1,2,3},X_1 \Y{(21)}{1,2,3},
\alphaq(X_1)\Y{(21)}{1,3,2},\Y{(21)}{1,3,2},\Y{(111)}{1,2,3}\}
$$
are linearly independent, 
\begin{maplegroup}
\begin{mapleinput}
\mapleinline{active}{1d}{
nops(findbasis([Y.3.123,Y.21.123,X.1 &c Y.21.123,alpha2(X.1) &c Y.21.132,
Y.21.132,Y.111.123]));}{%
}
\end{mapleinput}

\mapleresult
\begin{maplelatex}
\mapleinline{inert}{2d}{6;}{%
\[
6
\]
}
\end{maplelatex}
\end{maplegroup}
\noindent
they may form a basis $S$ for the Hecke algebra $H_{\BF}(3,q).$ We define the Garnir element $\G{(21)}{1,1}$ to be equal to $X_1.$\footnote{The above computations could have been performed with $X_2$ or $X_3$ instead of $X_1,$ and the results would have been similar.}
\begin{maplegroup}
\begin{mapleinput}
\mapleinline{active}{1d}{
alias(t1=K[6]*q^2+K[4]*q^2-K[5]*q-K[6]*q-K[4]*q+K[2]+K[4]):G21.1.1:=bexpand(X1);}{%
}
\end{mapleinput}

\mapleresult
\begin{maplelatex}
\mapleinline{inert}{2d}{G2111 :=
(K[2]*q-K[6]*q+K[4])*Id+K[2]*b1+t1*b2/q+K[4]*b12+K[5]*b21+K[6]*b121;}{
\[
\mathit{G2111} := ({K_{2}}\,q - {K_{6}}\,q + {K_{4}})\,\mathit{Id
} + {K_{2}}\,\mathit{b1} + {\displaystyle \frac {\mathit{t1}\,
\mathit{b2}}{q}}  + {K_{4}}\,\mathit{b12} + {K_{5}}\,\mathit{b21}
 + {K_{6}}\,\mathit{b121}
\]
}
\end{maplelatex}

\end{maplegroup}
\begin{maplegroup}
\begin{mapleinput}
\mapleinline{active}{1d}{S:=[Y.3.123,Y.21.123,G21.1.1 &c Y.21.123,
    alpha2(G21.1.1) &c Y.21.132,Y.21.132,Y.111.123]:}{%
}
\end{mapleinput}

\end{maplegroup}
\noindent
Procedure {\tt yexpand} (See Appendix 2) is used to expand elements in the Hecke algebra $H_{\BF}(3,q)$ in terms of the Young basis $S.$ For example, 

\begin{maplegroup}
\begin{mapleinput}
\mapleinline{active}{1d}{yexpand(Y3.123);}{%
}
\end{mapleinput}

\mapleresult
\begin{maplelatex}
\mapleinline{inert}{2d}{Y3.123;}{%
\[
\mathit{Y3}.123
\]
}
\end{maplelatex}

\end{maplegroup}
\begin{maplegroup}
\begin{mapleinput}
\mapleinline{active}{1d}{yexpand(&c(G21.1.1,Y21.123));}{%
}
\end{mapleinput}

\mapleresult
\begin{maplelatex}
\mapleinline{inert}{2d}{`&c`(G21.1.1,Y21.123);}{%
\[
\mathit{G21}.1.1\,\mathrm{\&c}\,\mathit{Y21}.123
\]
}
\end{maplelatex}

\end{maplegroup}
\begin{maplegroup}
\begin{mapleinput}
\mapleinline{active}{1d}{yexpand(&c(alpha2(G21.1.1),Y21.132));}{%
}
\end{mapleinput}

\mapleresult
\begin{maplelatex}
\mapleinline{inert}{2d}{`&c`(alpha2(G21.1.1),Y21.132);}{%
\[
\alpha 2(\mathit{G21}.1.1)\,\mathrm{\&c}\,\mathit{Y21}.132
\]
}
\end{maplelatex}

\end{maplegroup}
\noindent
Recall that the original basis in the Hecke algebra was: $\{\Id,b_1,b_2,b_{12},b_{21},b_{121}\}.$ 
Each original basis element should be representable in terms of the Young basis $S:$
\begin{maplegroup}
\begin{mapleinput}
\mapleinline{active}{1d}{'Id'=yexpand(Id);}{%
}
\end{mapleinput}

\mapleresult
\begin{maplelatex}
\mapleinline{inert}{2d}{Id = Y3.123+Y21.123+Y21.132+Y111.123;}{%
\[
\mathit{Id}=\mathit{Y3}.123 + \mathit{Y21}.123 + \mathit{Y21}.132
 + \mathit{Y111}.123
\]
}
\end{maplelatex}
\end{maplegroup}

\begin{maplegroup}
\begin{mapleinput}
\mapleinline{active}{1d}{'b1'=yexpand(b1);}{%
}
\end{mapleinput}

\mapleresult
\begin{maplelatex}
\mapleinline{inert}{2d}{b1 =
Y3.123+w2*Y21.123/(1+q)+q*w1*`&c`(G21.1.1,Y21.123)/(w3*(1+q))-q^2*`&c`
(alpha2(G21.1.1),Y21.132)/w4+q*w5*Y21.132/w6-q*Y111.123;}{%
\maplemultiline{
\mathit{b1}=\mathit{Y3}.123 + {\displaystyle \frac {\mathit{w2}\,
\mathit{Y21}.123}{1 + q}}  + {\displaystyle \frac {q\,\mathit{w1}
\,(\mathit{G21}.1.1\,\mathrm{\&c}\,\mathit{Y21}.123)}{\mathit{w3}
\,(1 + q)}}  \\
\mbox{} - {\displaystyle \frac {q^{2}\,(\alpha 2(\mathit{G21}.1.1
)\,\mathrm{\&c}\,\mathit{Y21}.132)}{\mathit{w4}}}  + 
{\displaystyle \frac {q\,\mathit{w5}\,\mathit{Y21}.132}{\mathit{
w6}}}  - q\,\mathit{Y111}.123 }
}
\end{maplelatex}

\end{maplegroup}
\begin{maplegroup}
\begin{mapleinput}
\mapleinline{active}{1d}{'b2'=yexpand(b2);}{%
}
\end{mapleinput}

\mapleresult
\begin{maplelatex}
\mapleinline{inert}{2d}{b2 =
Y3.123-w7*q*Y21.123/(1+q)-w8*q*`&c`(G21.1.1,Y21.123)/(w3*(1+q))+q^3*`&
c`(alpha2(G21.1.1),Y21.132)/w4-w9*Y21.132/w6-q*Y111.123;}{%
\maplemultiline{
\mathit{b2}=\mathit{Y3}.123 - {\displaystyle \frac {\mathit{w7}\,
q\,\mathit{Y21}.123}{1 + q}}  - {\displaystyle \frac {\mathit{w8}
\,q\,(\mathit{G21}.1.1\,\mathrm{\&c}\,\mathit{Y21}.123)}{\mathit{
w3}\,(1 + q)}}  \\
\mbox{} + {\displaystyle \frac {q^{3}\,(\alpha 2(\mathit{G21}.1.1
)\,\mathrm{\&c}\,\mathit{Y21}.132)}{\mathit{w4}}}  - 
{\displaystyle \frac {\mathit{w9}\,\mathit{Y21}.132}{\mathit{w6}}
}  - q\,\mathit{Y111}.123 }
}
\end{maplelatex}

\end{maplegroup}
\noindent
where $w_i,i=1,\ldots,9,$ are polynomials in $q$ parameterized in terms of $K_2,K_4,K_5,K_6.$ They have been defined as Maple aliases and are shown in Appendix 2.

\section{Singular Value Decomposition}
\label{sec:svd}

Our next application of Clifford algebras will be to the Singular Value Decomposition (SVD) of a matrix \cite{Strang}. There are many uses of SVD such as in image processing, description of the so called {\it principal gains\/} in a multivariable system \cite{Maciejowski}, or in an automated data indexing known as {\it Latent Semantic Indexing\/} (or LSI). LSI presents a very interesting and useful technique in information retrieval models and it is based on the SVD \cite{LSI}. While in these practical cases computations are done numerically, it may be of interest to ask whether such decomposition of a matrix can be performed in the framework of Clifford algebras. That is, if any new insights, theoretical or otherwise, into such decomposition could be gained when stated in the Clifford algebra language. In this section we will present examples of such computations. 

We will explore a well-known fact that when $p-q\not=1\bmod 4,$ Clifford algebra $\cl_{p,q}$ is a simple algebra of dimension $2^n,\,n=p+q,$ isomorphic to a full matrix algebra $\Mat(2^k,\BK)$ of $2^k \times 2^k$ matrices\footnote{The value of $k=q-r_{q-p},$ where $r_i$ is the Radon-Hurwitz number. The Radon-Hurwitz number is defined by a recursion as $r_{i+8}=r_{i}+4$ and these initial values: $r_{0}=0,\,r_{1}=1,\,r_{2}=r_{3}=2,\,r_{4}=r_{5}=r_{6}=r_{7}=3.$} with entries in $\BK$ which is $\BR,\,\BC,$ or $\BH$ (see \cite{Symbolic}). Thus, any operation performed on a matrix $A$ can be expressed as an operation on a corresponding to it element $p$ in $\cl_{p,q}.$ The choice of the signature $(p,q)$ depends on the size of $A$ and the division ring~$\BK.$ Of course, for computational reasons one should find the smallest Clifford algebra $\cl_{p,q}$ such that the given matrix $A$ can be embedded into 
$\Mat(2^k,\BK) \stackrel{\varphi}{\simeq} \cl_{p,q}.$ In the following we will use the same approach as in \cite{Matrix} where a technique for matrix exponentiation based on the isomorphism $\varphi$ was presented. In particular, we will use a faithful spinor representation of $\cl_{p,q}$ in a minimal left ideal $S=\cl_{p,q}f$ generated by a primitive idempotent $f.$ Symbolic computations of such representations with {\sc CLIFFORD} were shown in \cite{Symbolic}.
 
Following \cite{Strang}, let $A$ be an $m \times n$ real matrix of rank $r.$ Then the SVD of $A$ is defined a factorization of $A$ into a product of three matrices $U, \Sigma, V^{-1}$ where $U$ and $V$ are orthogonal matrices $m \times m$ and $n \times n$ respectively, and $\Sigma$ is a $m \times n$ matrix containing {\it singular values\/} of $A$ on its ``diagonal".
\begin{equation}
A=U\Sigma V^{-1}, \quad \UT U = I, \quad \VT V = I.
\label{eq:svd}
\end{equation}

The matrices $V=[v_1|v_2|\ldots|v_n]$ and $U=[u_1|u_2|\ldots|u_m]$ contain orthonormal bases for all four fundamental spaces of $A.$ Namely, the first $r$ columns $v_1,v_2,\ldots,v_r$ of $V$ provide a basis for the row space $\mathcal{R}(\AT)$ while the remaining $n-r$ columns of $V$ provide a basis for the null space $\mathcal{N}(A).$ Likewise, the first $r$ columns $u_1,u_2,\ldots,u_r$ of $U$ provide a basis for the column space $\mathcal{C}(A)$ while the remaining $m-r$ columns of $U$ provide a basis for the left-null space $\mathcal{N}(A^T).$ 
Vectors $v_i$ are the normalized eigenvectors of $\ATA$ while vectors $u_i$ are the normalized eigenvectors of $\AAT.$ For $i=1,\ldots,r,$ these vectors can be chosen to be related via the positive singular values $\sigma_i$ of $A$ which are just the square roots of the eigenvalues of $\ATA$ (or of 
$\AAT.)$ Namely,
\begin{equation}
A v_i = \sigma_i u_i, \quad i=1,\ldots,r.
\label{eq:avsu}
\end{equation}
It is a little tricky to make sure that the above relation is satisfied: this is because the choice of vectors $u_i$ is independent of the choice of vectors $v_i.$ However, it is always possible to do so as we will see below (see also \cite{Strang}). In order to complete the picture, the orthonormal set $\{v_1,\ldots,v_r\}$ needs to be completed to a full orthonormal basis for $\BR^n$ while $\{u_1,\ldots,u_r\}$ needs to be completed to a full orthonormal basis for $\BR^m.$ Since the additional vectors are being annihilated by $A$ and $A^T$ respectively, that is, they are eigenvectors of $A$ and $\AT$ (or of $\ATA$ and $\AAT)$ that correspond to the eigenvalue $0,$ care has to be exercised when finding them. For example, while the eigenvectors of the symmetric matrix $A A^T$ are automatically orthogonal provided they correspond to different eigenvalues, eigenvectors of $A A^T$ that correspond to the $0$ eigenvalue don't need to be orthogonal: in this case the Gram-Schmidt orthogonalization process is used to complete the two sets.    

\subsection{Singular Value Decomposition of a $2 \times 2$ matrix of rank $2$}
\label{sec:svd1}

In this section we present our first example of SVD applied to a $2 \times 2$ real matrix of rank $2.$ The purpose of this example is just to show step by step how finding the SVD of a matrix can be done in the Clifford algebra language. Reader is encouraged to perform these computations with {\sc CLIFFORD} and an additional package {\sc asvd} which is described in Appendix~3.
\begin{maplegroup}
\begin{mapleinput}
\mapleinline{active}{1d}{A:=matrix(2,2,[2,3,1,2]);#defining A
m:=rowdim(A): #number of rows of A is m
n:=coldim(A): #number of columns of A is n}{%
}
\end{mapleinput}

\mapleresult
\begin{maplelatex}
\mapleinline{inert}{2d}{A := matrix([[2, 3], [1, 2]]);}{%
\[
A :=  \left[ 
{\begin{array}{rr}
2 & 3 \\
1 & 2
\end{array}}
 \right] 
\]
}
\end{maplelatex}

\end{maplegroup}
\noindent
Since $A \in \Mat(2,\BR),$ we need to find $(p,q)$ such that $\cl_{p,q} \stackrel{\varphi}{\simeq}\Mat(2,\BR).$ As shown next, we have two choices for the signature:
\begin{maplegroup}
\begin{mapleinput}
\mapleinline{active}{1d}{all_sigs(2..2,real,simple);}{%
}
\end{mapleinput}

\mapleresult
\begin{maplelatex}
\mapleinline{inert}{2d}{[[1, 1], [2, 0]];}{%
\[
[[1, \,1], \,[2, \,0]]
\]
}
\end{maplelatex}
\end{maplegroup}
\noindent
Thus, we can pick either $\cl_{1,1}$ or $\cl_{2,0}.$ Our choice is $\cl_{2,0}.$ We define a bilinear form $B=\diag(1,1)$ and display information about $\cl_{2,0}.$ 
\begin{maplegroup}
\begin{mapleinput}
\mapleinline{active}{1d}{dim:=2:B:=diag(1,1):eval(makealiases(dim)):data:=clidata();}{%
}
\end{mapleinput}

\mapleresult
\begin{maplelatex}
\mapleinline{inert}{2d}{data := [real, 2, simple, 1/2*Id+1/2*e1, [Id, e2], [Id], [Id,e2]];}{%
\[
\mathit{data} := [\mathit{real}, \,2, \,\mathit{simple}, \,
{\displaystyle \frac {1}{2}} \,\mathit{Id} + {\displaystyle 
\frac {1}{2}} \,\mathit{e1}, \,[\mathit{Id}, \,\mathit{e2}], \,[
\mathit{Id}], \,[\mathit{Id}, \,\mathit{e2}]]
\]
}
\end{maplelatex}
\end{maplegroup}
\noindent
The above output means that $\cl_{2,0}$ is a simple algebra isomorphic with $\Mat(2,\BR);$ that the element $\frac12 + \frac12 \e_1$ is a primitive idempotent which we will call $f;$ that the list $[Id,e2]$ shown as the fourth entry displays generators of a minimal left-ideal $\cl_{2,0}f$ considered as vector space over $\BR;$ that the division ring $\BK=f\cl_{2,0}f = <Id>_\BR \simeq \BR;$ and that the last list $[Id,e2]$ gives generators of $\cl_{2,0}f$ over $\BK,$ and since $\BK \simeq \BR,$ it is the same as the fourth list.\footnote{For more information see \cite{Symbolic} and {\sc CLIFFORD}'s help pages.}. In the following, we define a Grassmann basis in $\cl_{2,0},$ assign the primitive idempotent to $f,$ and generate a spinor basis 
in~$\cl_{2,0}f.$
\begin{maplegroup}
\begin{mapleinput}
\mapleinline{active}{1d}{clibas:=cbasis(dim); #ordered basis in Cl(2,0)}{%
}
\end{mapleinput}

\mapleresult
\begin{maplelatex}
\mapleinline{inert}{2d}{clibas := [Id, e1, e2, e12];}{%
\[
\mathit{clibas} := [\mathit{Id}, \,\mathit{e1}, \,\mathit{e2}, \,\mathit{e12}]
\]
}
\end{maplelatex}

\end{maplegroup}
\begin{maplegroup}
\begin{mapleinput}
\mapleinline{active}{1d}{f:=data[4]:#a primitive idempotent in Cl(2,0)
SBgens:=data[5]:#generators for a real basis in S
FBgens:=data[6]:#generators for K}{%
}
\end{mapleinput}
\end{maplegroup}
\noindent
Here {\tt SBgens} is a $\BK$-basis for $S=\cl_{2,0}f.$ Since for the signature $(2,0)$ we have 
$\BK \simeq \BR,$ $S \simeq \BR^2,$ and $\cl_{2,0}\simeq\Mat(2,\BR),$ the output from 
the procedure {\tt spinorKbasis} shown below has two basis elements and their generators 
modulo~$f:$ 
\begin{maplegroup}
\begin{mapleinput}
\mapleinline{active}{1d}{Kbasis:=spinorKbasis(SBgens,f,FBgens,'left');}{%
}
\end{mapleinput}

\mapleresult
\begin{maplelatex}
\mapleinline{inert}{2d}{Kbasis := [[1/2*Id+1/2*e1, 1/2*e2-1/2*e12], [Id, e2], left];}{%
\[
\mathit{Kbasis} := [[{\displaystyle \frac {1}{2}} \,\mathit{Id}
 + {\displaystyle \frac {1}{2}} \,\mathit{e1}, \,{\displaystyle 
\frac {1}{2}} \,\mathit{e2} - {\displaystyle \frac {1}{2}} \,
\mathit{e12}], \,[\mathit{Id}, \,\mathit{e2}], \,\mathit{left}]
\]
}
\end{maplelatex}
\end{maplegroup}
\noindent
Thus, the real spinor basis in $S$ consists of the following two polynomials: 
\begin{maplegroup}
\begin{mapleinput}
\mapleinline{active}{1d}{for i from 1 to nops(Kbasis[1]) do f.i:=Kbasis[1][i] od;}{%
}
\end{mapleinput}

\mapleresult
\begin{maplelatex}
\mapleinline{inert}{2d}{f1 := 1/2*Id+1/2*e1;f2 := 1/2*e2-1/2*e12;}{%
\[
\mathit{f1} := {\displaystyle \frac {1}{2}} \,\mathit{Id} + 
{\displaystyle \frac {1}{2}} \,\mathit{e1}, \quad 
\mathit{f2} := {\displaystyle \frac {1}{2}} \,\mathit{e2} - 
{\displaystyle \frac {1}{2}} \,\mathit{e12}.
\]
}
\end{maplelatex}
\end{maplegroup}
\noindent
Now, we compute matrices $M_1,M_2,M_3,M_4$ representing each of the four basis elements $\{\Id,\e_1,\e_2,\e_{12}\}$ in $\cl_{2,0}.$\footnote{Since a similar computation was done in \cite{Matrix,Symbolic}, we won't display the matrices.}
\begin{maplegroup}
\begin{mapleinput}
\mapleinline{active}{1d}{for i from 1 to nops(clibas) do M[i]:=subs(Id=1,matKrepr(clibas[i])) od:
}{%
}
\end{mapleinput}
\end{maplegroup}
\noindent
We will use a procedure {\tt phi} which gives the isomorphism $\varphi$ from $\Mat(2,\BR)$ to $\cl_{2,0}.$  This way we can find the image $p$ in $\cl_{2,0}$ of any real $2 \times 2$ real matrix $A.$ Knowing the image $\varphi(M_i)$ of each matrix $M_i$ in terms of some Clifford polynomial  in $\cl_{2,0},$ we can easily find the image $p=\varphi(A)$ of $A$ as follows:\footnote{From now on, as use Maple {\tt alias(t=transpose)}, that is, {\tt t(A)} denotes the matrix transposition in Maple.}
\begin{maplegroup}
\begin{mapleinput}
\mapleinline{active}{1d}{p:=phi(A,M); #finding image of A in Cl(2,0)}{%
}
\end{mapleinput}

\mapleresult
\begin{maplelatex}
\mapleinline{inert}{2d}{p := 2*Id+2*e2+e12;}{%
\[
p := 2\,\mathit{Id} + 2\,\mathit{e2} + \mathit{e12}
\]
}
\end{maplelatex}

\end{maplegroup}
\begin{maplegroup}
\begin{mapleinput}
\mapleinline{active}{1d}{pT:=phi(t(A),M); #finding image of t(A) in Cl(2,0)}{%
}
\end{mapleinput}

\mapleresult
\begin{maplelatex}
\mapleinline{inert}{2d}{pT := 2*Id+2*e2-e12;}{%
\[
\mathit{pT} := 2\,\mathit{Id} + 2\,\mathit{e2} - \mathit{e12}
\]
}
\end{maplelatex}

\end{maplegroup}
Next, we compute a symmetric matrix $A^T A$ (denoted in Maple as {\tt ATA}), its characteristic polynomial, eigenvalues, and its orthonormal eigenvectors $v_1,v_2.$  Vectors $v_1$ and $v_2$ will become columns of an orthogonal matrix $V$ needed for SVD of $A$\footnote{We use procedure {\tt radsimplify} to simplify radicals in matrices and vectors.}:
\begin{maplegroup}
\begin{mapleinput}
\mapleinline{active}{1d}{ATA:=evalm(t(A) &* A); #finding matrix ATA}{%
}
\end{mapleinput}

\mapleresult
\begin{maplelatex}
\mapleinline{inert}{2d}{ATA := matrix([[5, 8], [8, 13]]);}{%
\[
\mathit{ATA} :=  \left[ 
{\begin{array}{rr}
5 & 8 \\
8 & 13
\end{array}}
 \right] 
\]
}
\end{maplelatex}

\end{maplegroup}
\begin{maplegroup}
\begin{mapleinput}
\mapleinline{active}{1d}{pTp:=phi(ATA,M); #finding image of ATA in Cl(2,0)}{%
}
\end{mapleinput}

\mapleresult
\begin{maplelatex}
\mapleinline{inert}{2d}{pTp := 9*Id-4*e1+8*e2;}{%
\[
\mathit{pTp} := 9\,\mathit{Id} - 4\,\mathit{e1} + 8\,\mathit{e2}
\]
}
\end{maplelatex}

\end{maplegroup}
\noindent
which should be the same as
\begin{maplegroup}
\begin{mapleinput}
\mapleinline{active}{1d}{'pTp'=cmul(pT,p);}{%
}
\end{mapleinput}

\mapleresult
\begin{maplelatex}
\mapleinline{inert}{2d}{pTp = 9*Id-4*e1+8*e2;}{%
\[
\mathit{pTp}=9\,\mathit{Id} - 4\,\mathit{e1} + 8\,\mathit{e2}
\]
}
\end{maplelatex}

\end{maplegroup}
\noindent
The minimum polynomial of $A^T A$ (or $p\,Tp = \varphi(\ATA))$ is:
\begin{maplegroup}
\begin{mapleinput}
\mapleinline{active}{1d}{climinpoly(pTp);}{%
}
\end{mapleinput}

\mapleresult
\begin{maplelatex}
\mapleinline{inert}{2d}{x^2-18*x+1;}{%
\[
x^{2} - 18\,x + 1
\]
}
\end{maplelatex}

\end{maplegroup}
\noindent
and, in this case, it is the same as the characteristic polynomial of $A^T A:$
\begin{maplegroup}
\begin{mapleinput}
\mapleinline{active}{1d}{pol:=charpoly(ATA,x);#characteristic polynomial of ATA}{%
}
\end{mapleinput}

\mapleresult
\begin{maplelatex}
\mapleinline{inert}{2d}{pol := x^2-18*x+1;}{%
\[
\mathit{pol} := x^{2} - 18\,x + 1
\]
}
\end{maplelatex}

\end{maplegroup}
In order to find eigenvalues and eigenvectors of $\ATA,$ we will use Maple's procedure {\tt eigenvects} modified by our own sorting via a new procedure {\tt assignL}. The latter displays a list containing two lists: one has the eigenvalues while the second has the eigenvectors.\footnote{The first entry $2$ in the output is just the number of eigenvectors.} In the following, we will assign the eigenvalues to $\lambda_1,\lambda_2$ and the (un-normalized, but orthogonal) eigenvectors we assign to $v_1,v_2.$
\begin{maplegroup}
\begin{mapleinput}
\mapleinline{active}{1d}{P:=assignL(sort([eigenvects(ATA)],byeigenvals));N:=P[1]:}{%
}
\end{mapleinput}

\mapleresult
\begin{maplelatex}
\mapleinline{inert}{2d}{P := [2, [9+4*sqrt(5), 9-4*sqrt(5)], [vector([1, 1/2+1/2*sqrt(5)]),
vector([1, 1/2-1/2*sqrt(5)])]];}{%
\[
P := [2, \,[9 + 4\,\sqrt{5}, \,9 - 4\,\sqrt{5}], \,[ \left[  \! 1
, \,{\displaystyle \frac {1}{2}}  + {\displaystyle \frac {1}{2}} 
\,\sqrt{5} \!  \right] , \, \left[  \! 1, \,{\displaystyle 
\frac {1}{2}}  - {\displaystyle \frac {1}{2}} \,\sqrt{5} \! 
 \right] ]]
\]
}
\end{maplelatex}
\end{maplegroup}
\begin{maplegroup}
\begin{mapleinput}
\mapleinline{active}{1d}{
for i from 1 to N do lambda.i:=P[2][i];v.i:=map(simplify,normalize(P[3][i])) od:}{%
}
\end{mapleinput}
\end{maplegroup}
\noindent
We can now verify that vectors $v_1,v_2$ are eigenvectors of $A^T A$ with the eigenvalues $\lambda_1,\lambda_2.$ 

\begin{maplegroup}
\begin{mapleinput}
\mapleinline{active}{1d}{for i from 1 to N do map(simplify,evalm(ATA &* v.i - lambda.i*v.i)) od;
}{%
}
\end{mapleinput}

\mapleresult
\begin{maplelatex}
\mapleinline{inert}{2d}{vector([0, 0]);vector([0, 0]);}{%
\[
[0, \,0], \; [0, \,0]
\]
}
\end{maplelatex}
\end{maplegroup}
\noindent
Similar verification can be done in $\cl_{2,0}$ since one can view the $1$-column eigenvectors $v_1,\,v_2$ as one-column spinors in $S.$ We simply convert the two vectors to spinors  $sv_1,sv_2$ which we express in the previously computed spinor basis $f_1,f_2.$
\begin{maplegroup}
\begin{mapleinput}
\mapleinline{active}{1d}{spinorbasis:=[''f1'',''f2'']:
for i from 1 to N do sv.i:=convert(v.i,spinor,spinorbasis) od;}{%
}
\end{mapleinput}

\mapleresult
\begin{maplelatex}
\mapleinline{inert}{2d}{sv1 :=2*f1/(sqrt(10+2*sqrt(5)))+(1+sqrt(5))*f2/(sqrt(10+2*sqrt(5)));
sv2 := 2*f1/(sqrt(10-2*sqrt(5)))-(-1+sqrt(5))*f2/(sqrt(10-2*sqrt(5)));}{%
\[
\mathit{sv1} := 2\,{\displaystyle \frac {\mathit{f1}}{\sqrt{10 + 
2\,\sqrt{5}}}}  + {\displaystyle \frac {(1 + \sqrt{5})\,\mathit{
f2}}{\sqrt{10 + 2\,\sqrt{5}}}},\quad
\mathit{sv2} := 2\,{\displaystyle \frac {\mathit{f1}}{\sqrt{10 - 
2\,\sqrt{5}}}}  - {\displaystyle \frac {( - 1 + \sqrt{5})\,
\mathit{f2}}{\sqrt{10 - 2\,\sqrt{5}}}}. 
\]
}
\end{maplelatex}
\end{maplegroup}
\noindent
Since $v_1,v_2$ are eigenvectors of $AT A,$ spinors $sv_1, sv_2$ must be eigenspinors of 
$pTp = \varphi(A^T A).$
\begin{maplegroup}
\begin{mapleinput}
\mapleinline{active}{1d}{for i from 1 to N do simplify((pTp - lambda.i) &c sv.i) od;}{%
}
\end{mapleinput}

\mapleresult
\begin{maplelatex}
\mapleinline{inert}{2d}{0,0;}{%
\[
0,0
\]
}
\end{maplelatex}
\end{maplegroup}
\noindent
We are now in position to define the orthogonal matrix $V=[v_1|v_2].$
\begin{maplegroup}
\begin{mapleinput}
\mapleinline{active}{1d}{V:=radsimplify(augment(v.(1..N))); #defining matrix V}{%
}
\end{mapleinput}

\mapleresult
\begin{maplelatex}
\mapleinline{inert}{2d}{V := matrix([[2*1/(sqrt(10+2*sqrt(5))), 2*1/(sqrt(10-2*sqrt(5)))],
[(1+sqrt(5))/(sqrt(10+2*sqrt(5))),
(-sqrt(5)+1)/(sqrt(10-2*sqrt(5)))]]);}{%
\[
V :=  \left[ 
{\begin{array}{cc}
2\,{\displaystyle \frac {1}{\sqrt{10 + 2\,\sqrt{5}}}}  & 2\,
{\displaystyle \frac {1}{\sqrt{10 - 2\,\sqrt{5}}}}  \\ [2ex]
{\displaystyle \frac {1 + \sqrt{5}}{\sqrt{10 + 2\,\sqrt{5}}}}  & 
{\displaystyle \frac { - \sqrt{5} + 1}{\sqrt{10 - 2\,\sqrt{5}}}} 
\end{array}}
 \right] 
\]
}
\end{maplelatex}
\end{maplegroup}
Since later we will need images of $V$ and $V^T$ under $\varphi,$ we compute them now and store under the variables $pV$ and $pVt.$ The fact that $V$ is orthogonal can be easily verified in the  matrix language; in $\cl_{2,0}$ it can be done as follows:
\begin{maplegroup}
\begin{mapleinput}
\mapleinline{active}{1d}{simplify(cmul(pVt,pV));}{%
}
\end{mapleinput}

\mapleresult
\begin{maplelatex}
\mapleinline{inert}{2d}{Id;}{%
\[
\mathit{Id}
\]
}
\end{maplelatex}
\end{maplegroup}
Now we repeat the above steps and apply them to $AA^T.$ In the process, we will find its eigenvectors $u_1,u_2.$  We must make sure that $A v_i = \sigma_i u_i$ where $\sigma_i=\sqrt{\lambda_i},\,i=1,2.$  This will require extra checking and possibly redefining of the $u$'s.
\begin{maplegroup}
\begin{mapleinput}
\mapleinline{active}{1d}{AAT:=evalm(A &* transpose(A)); #computing AAT}{%
}
\end{mapleinput}

\mapleresult
\begin{maplelatex}
\mapleinline{inert}{2d}{AAT := matrix([[13, 8], [8, 5]]);}{%
\[
\mathit{AAT} :=  \left[ 
{\begin{array}{rr}
13 & 8 \\
8 & 5
\end{array}}
 \right] 
\]
}
\end{maplelatex}
\end{maplegroup}
\noindent
The image of $A A^T$ under $\varphi$ in $\cl_{2,0}$ we denote as $ppT.$
\begin{maplegroup}
\begin{mapleinput}
\mapleinline{active}{1d}{ppT:=phi(AAT,M); #finding image of AAT in Cl(2,0)}{%
}
\end{mapleinput}

\mapleresult
\begin{maplelatex}
\mapleinline{inert}{2d}{ppT := 9*Id+4*e1+8*e2;}{%
\[
\mathit{ppT} := 9\,\mathit{Id} + 4\,\mathit{e1} + 8\,\mathit{e2}
\]
}
\end{maplelatex}

\end{maplegroup}
\noindent
In this case, the minimal polynomial of $ppT$ and the characteristic polynomial of $A A^T$ are the same.
\begin{maplegroup}
\begin{mapleinput}
\mapleinline{active}{1d}{
pol2:=charpoly(AAT,lambda); #finding the characteristic polynomial of AAT
}{%
}
\end{mapleinput}

\mapleresult
\begin{maplelatex}
\mapleinline{inert}{2d}{pol2 := lambda^2-18*lambda+1;}{%
\[
\mathit{pol2} := \lambda ^{2} - 18\,\lambda  + 1
\]
}
\end{maplelatex}

\end{maplegroup}
\begin{maplegroup}
\begin{mapleinput}
\mapleinline{active}{1d}{'ppT'=climinpoly(ppT);}{%
}
\end{mapleinput}

\mapleresult
\begin{maplelatex}
\mapleinline{inert}{2d}{ppT = x^2-18*x+1;}{%
\[
\mathit{ppT}=x^{2} - 18\,x + 1
\]
}
\end{maplelatex}
\end{maplegroup}
\noindent
Since matrices $\ATA$ and $\AAT$ have the same characteristic polynomials, their eigenvalues will be the same.  We define therefore the {\it singular values\/} $\sigma_1$ and $\sigma_2$ of $A:$
\begin{maplegroup}
\begin{mapleinput}
\mapleinline{active}{1d}{for i from 1 to N do sigma.i:=sqrt(lambda.i) od;}{%
}
\end{mapleinput}

\mapleresult
\begin{maplelatex}
\mapleinline{inert}{2d}{sigma1 := sqrt(5)+2;sigma2 := sqrt(5)-2;}{%
\[
\sigma 1 := \sqrt{5} + 2, \quad \sigma 2 := \sqrt{5} - 2
\]
}
\end{maplelatex}

\end{maplegroup}
\noindent
When we compute the eigenvectors $u_1, u_2$ of $\AAT,$ we will not necessarily have $A v_i = \sigma_i u_i, i=1,2.$ This is because the choice of $u_1,u_2$ is not consistent with the choice of $v_1,v_2.$ 
\begin{maplegroup}
\begin{mapleinput}
\mapleinline{active}{1d}{P:=assignL(sort([eigenvects(AAT)],byeigenvals)):
for i from 1 to N do lambda.i:=P[2][i]; u.i:=map(simplify,normalize(P[3][i])) od:}{%
}
\end{mapleinput}
\end{maplegroup}
\noindent
However, $A v_1 = \sigma_1 u_1$ while $A v_2 = -\sigma_2 u_2:$ 
\begin{maplegroup}
\begin{mapleinput}
\mapleinline{active}{1d}{radsimplify(evalm(A &* v1-sigma1*u1)); #this one checks out but}{%
}
\end{mapleinput}

\mapleresult
\begin{maplelatex}
\mapleinline{inert}{2d}{vector([0, 0]);}{%
\[
[0, \,0]
\]
}
\end{maplelatex}

\end{maplegroup}
\begin{maplegroup}
\begin{mapleinput}
\mapleinline{active}{1d}{radsimplify(evalm(A &* v2-sigma2*u2)); #this one does not check out
radsimplify(evalm(A &* v2+sigma2*u2)); #this one does check out}{%
}
\end{mapleinput}

\mapleresult
\begin{maplelatex}
\mapleinline{inert}{2d}{vector([(14-6*sqrt(5))/(sqrt(10-2*sqrt(5))),
(8-4*sqrt(5))/(sqrt(10-2*sqrt(5)))]);}{%
\[
 \left[  \! {\displaystyle \frac {14 - 6\,\sqrt{5}}{\sqrt{10 - 2
\,\sqrt{5}}}} , \,{\displaystyle \frac {8 - 4\,\sqrt{5}}{\sqrt{10
 - 2\,\sqrt{5}}}}  \!  \right] 
\]
}
\end{maplelatex}

\begin{maplelatex}
\mapleinline{inert}{2d}{vector([0, 0]);}{%
\[
[0, \,0]
\]
}
\end{maplelatex}

\end{maplegroup}
\noindent
Notice that the set $\{u_1, u_2\}$ is orthonormal, but so is $\{u_1,-u_2\}.$  Let's re-define $u_2$ as $-u_2$ and call it $u_{22}.$  For completeness we rename $u_1$ as $u_{11}:$
\begin{maplegroup}
\begin{mapleinput}
\mapleinline{active}{1d}{u11:=evalm(u1):u22:=evalm(-u2):}{%
}
\end{mapleinput}
\end{maplegroup}

In the Clifford algebra $\cl_{2,0},$ we need to perform similar computations with $v_1,v_2.$ The images $\varphi(u_{11}),\varphi(u_{22})$ contained in the spinor ideal need to be found first. We call them $su_1$ and $su_2.$
\begin{maplegroup}
\begin{mapleinput}
\mapleinline{active}{1d}{for i from 1 to N do su.i:=convert(u.i.i,spinor,spinorbasis) od;}{%
}
\end{mapleinput}

\mapleresult
\begin{maplelatex}
\mapleinline{inert}{2d}{su1 := (1+sqrt(5))*f1/(sqrt(10+2*sqrt(5)))+2*f2/(sqrt(10+2*sqrt(5)));
su2 := (-1+sqrt(5))*f1/(sqrt(10-2*sqrt(5)))-2*f2/(sqrt(10-2*sqrt(5)));}{%
\[
\mathit{su1} := {\displaystyle \frac {(1 + \sqrt{5})\,\mathit{f1}
}{\sqrt{10 + 2\,\sqrt{5}}}}  + 2\,{\displaystyle \frac {\mathit{
f2}}{\sqrt{10 + 2\,\sqrt{5}}}}, \quad 
\mathit{su2} := {\displaystyle \frac {( - 1 + \sqrt{5})\,\mathit{
f1}}{\sqrt{10 - 2\,\sqrt{5}}}}  - 2\,{\displaystyle \frac {
\mathit{f2}}{\sqrt{10 - 2\,\sqrt{5}}}}. 
\]
}
\end{maplelatex}
\end{maplegroup}
\noindent
The verification of the condition (\ref{eq:avsu}) in $\cl_{2,0}$ looks as follows:
\begin{maplegroup}
\begin{mapleinput}
\mapleinline{active}{1d}{for i from 1 to N do simplify(p &c sv.i-sigma.i*su.i) od; }{%
}
\end{mapleinput}

\mapleresult
\begin{maplelatex}
\mapleinline{inert}{2d}{0,0;}{%
\[
0, \; 0
\]
}
\end{maplelatex}
\end{maplegroup}
\noindent
Now we may define the orthogonal matrix $U=[u_{11}|u_{22}]$ and its image $\varphi(U)$ in $\cl_{2,0}$ which we call $pU:$\footnote{For a later verification we will also need $pUt = \varphi(U^T).$ Expressions $\%1$ and $\%2$ showing up in the Maple output for $pU$ are just place holders for $\sqrt{10 + 2\,\sqrt{5}}$ and $\sqrt{10 - 2\,\sqrt{5}}$ respectively as shown at the end of the display.}
\begin{maplegroup}
\begin{mapleinput}
\mapleinline{active}{1d}{U:=radsimplify(augment(u11,u22)); #defining matrix U}{%
}
\end{mapleinput}

\mapleresult
\begin{maplelatex}
\mapleinline{inert}{2d}{U := matrix([[(1+sqrt(5))/(sqrt(10+2*sqrt(5))),
(-1+sqrt(5))/(sqrt(10-2*sqrt(5)))], [2*1/(sqrt(10+2*sqrt(5))),
-2*1/(sqrt(10-2*sqrt(5)))]]);}{%
\[
U :=  \left[ 
{\begin{array}{cc}
{\displaystyle \frac {1 + \sqrt{5}}{\sqrt{10 + 2\,\sqrt{5}}}}  & 
{\displaystyle \frac { - 1 + \sqrt{5}}{\sqrt{10 - 2\,\sqrt{5}}}} 
 \\ [2ex]
2\,{\displaystyle \frac {1}{\sqrt{10 + 2\,\sqrt{5}}}}  &  - 2\,
{\displaystyle \frac {1}{\sqrt{10 - 2\,\sqrt{5}}}} 
\end{array}}
 \right] 
\]
}
\end{maplelatex}

\end{maplegroup}
\begin{maplegroup}
\begin{mapleinput}
\mapleinline{active}{1d}{pU:=phi(U,M);pUt:=phi(t(U),M):}{%
}
\end{mapleinput}

\mapleresult
\begin{maplelatex}
\mapleinline{inert}{2d}{pU :=
Id*(1/20*sqrt(10+2*sqrt(5))*sqrt(5)-1/8*sqrt(10-2*sqrt(5))-1/40*sqrt(1
0-2*sqrt(5))*sqrt(5))+(1/20*sqrt(10+2*sqrt(5))*sqrt(5)+1/8*sqrt(10-2*s
qrt(5))+1/40*sqrt(10-2*sqrt(5))*sqrt(5))*e1+(1/20*sqrt(10-2*sqrt(5))*s
qrt(5)+1/8*sqrt(10+2*sqrt(5))-1/40*sqrt(10+2*sqrt(5))*sqrt(5))*e2+(1/2
0*sqrt(10-2*sqrt(5))*sqrt(5)-1/8*sqrt(10+2*sqrt(5))+1/40*sqrt(10+2*sqr
t(5))*sqrt(5))*e12;}{%
\maplemultiline{
\mathit{pU} := \mathit{Id}\,({\displaystyle \frac {1}{20}} \,
\mathrm{\%1}\,\sqrt{5} - {\displaystyle \frac {1}{8}} \,\mathrm{
\%2} - {\displaystyle \frac {1}{40}} \,\mathrm{\%2}\,\sqrt{5}) + 
({\displaystyle \frac {1}{20}} \,\mathrm{\%1}\,\sqrt{5} + 
{\displaystyle \frac {1}{8}} \,\mathrm{\%2} + {\displaystyle 
\frac {1}{40}} \,\mathrm{\%2}\,\sqrt{5})\,\mathit{e1} \\
\mbox{} + ({\displaystyle \frac {1}{20}} \,\mathrm{\%2}\,\sqrt{5}
 + {\displaystyle \frac {1}{8}} \,\mathrm{\%1} - {\displaystyle 
\frac {1}{40}} \,\mathrm{\%1}\,\sqrt{5})\,\mathit{e2} + (
{\displaystyle \frac {1}{20}} \,\mathrm{\%2}\,\sqrt{5} - 
{\displaystyle \frac {1}{8}} \,\mathrm{\%1} + {\displaystyle 
\frac {1}{40}} \,\mathrm{\%1}\,\sqrt{5})\,\mathit{e12} \\
\mathrm{\%1} := \sqrt{10 + 2\,\sqrt{5}} \\
\mathrm{\%2} := \sqrt{10 - 2\,\sqrt{5}} }
}
\end{maplelatex}

\end{maplegroup}
\noindent
The fact that $U$ is an orthogonal matrix can be easily now checked both in the matrix language and in the Clifford language:
\begin{maplegroup}
\begin{mapleinput}
\mapleinline{active}{1d}{radsimplify(evalm(t(U) &* U));#U is an orthogonal matrix}{%
}
\end{mapleinput}

\mapleresult
\begin{maplelatex}
\mapleinline{inert}{2d}{matrix([[1, 0], [0, 1]]);}{%
\[
 \left[ 
{\begin{array}{rr}
1 & 0 \\
0 & 1
\end{array}}
 \right] 
\]
}
\end{maplelatex}
\end{maplegroup}
\begin{maplegroup}
\begin{mapleinput}
\mapleinline{active}{1d}{simplify(pUt &c pU);}{%
}
\end{mapleinput}

\mapleresult
\begin{maplelatex}
\mapleinline{inert}{2d}{Id;}{%
\[
\mathit{Id}
\]
}
\end{maplelatex}
\end{maplegroup}
\noindent
Finally, we define matrix $\Sigma$ using a procedure {\tt makediag}. Recall \cite{Strang} that $\Sigma$ has the same dimensions as the original matrix $A$ and that $\STS,\SST$ are the diagonal forms of $\ATA$ and $\AAT$ respectively. In this example matrices $\STS$ and $\SST$ are the same since $\Sigma$ is a square diagonal matrix. Normally these matrices are different although their nonzero ``diagonal" entries are the same. 
\begin{equation}
\ATA = V \STS V^T,\quad \AAT = U \SST U^T,\quad
\Sigma = \left[\matrix{\sigma_1 & 0 \cr
	                 0 & \sigma_2 }\right], \quad
\STS = \SST = \left[\matrix{\sigma_1^2 & 0 \cr
	                      0 & \sigma_2^2 }\right].
\label{eq:sigmas}
\end{equation}
\begin{maplegroup}
\begin{mapleinput}
\mapleinline{active}{1d}{Sigma:=makediag(m,n,[seq(sigma.i,i=1..N)]);
STS,SST:=evalm(t(Sigma) &* Sigma),evalm(Sigma &* t(Sigma));}{%
}
\end{mapleinput}

\mapleresult
\begin{maplelatex}
\mapleinline{inert}{2d}{Sigma := matrix([[sqrt(5)+2, 0], [0, sqrt(5)-2]]);}{%
\[
\Sigma  :=  \left[ 
{\begin{array}{cc}
\sqrt{5} + 2 & 0 \\
0 & \sqrt{5} - 2
\end{array}}
 \right] 
\]
}
\end{maplelatex}

\begin{maplelatex}
\mapleinline{inert}{2d}{STS, SST := matrix([[(sqrt(5)+2)^2, 0], [0, (sqrt(5)-2)^2]]),
matrix([[(sqrt(5)+2)^2, 0], [0, (sqrt(5)-2)^2]]);}{%
\[
\mathit{STS}, \,\mathit{SST} :=  \left[ 
{\begin{array}{cc}
(\sqrt{5} + 2)^{2} & 0 \\
0 & (\sqrt{5} - 2)^{2}
\end{array}}
 \right] , \, \left[ 
{\begin{array}{cc}
(\sqrt{5} + 2)^{2} & 0 \\
0 & (\sqrt{5} - 2)^{2}
\end{array}}
 \right] 
\]
}
\end{maplelatex}

\end{maplegroup}
\begin{maplegroup}
\begin{mapleinput}
\mapleinline{active}{1d}{pSigma,pSTS,pSST:=phi(Sigma,M),phi(STS,M,FBgens),phi(SST,M);}{%
}
\end{mapleinput}

\mapleresult
\begin{maplelatex}
\mapleinline{inert}{2d}{pSigma, pSTS, pSST := sqrt(5)*Id+2*e1, 9*Id+4*sqrt(5)*e1,
9*Id+4*sqrt(5)*e1;}{%
\[
\mathit{pSigma}, \,\mathit{pSTS}, \,\mathit{pSST} := \sqrt{5}\,
\mathit{Id} + 2\,\mathit{e1}, \,9\,\mathit{Id} + 4\,\sqrt{5}\,
\mathit{e1}, \,9\,\mathit{Id} + 4\,\sqrt{5}\,\mathit{e1}
\]
}
\end{maplelatex}
\end{maplegroup}
\noindent 
We should be able to verify in $\cl_{2,0}$ the following two factorizations of $\AAT$ and $\ATA:$
\begin{eqnarray}
\ATA &=& V \STS \VT  \label{eq:ATA} \\
\AAT &=& U \SST \UT  \label{eq:AAT}
\end{eqnarray}
\noindent
like this:
\begin{maplegroup}
\begin{mapleinput}
\mapleinline{active}{1d}{
evalb(pTp=simplify(pV &c pSTS &c pVt)), evalb(ppT=simplify(pU &c pSST &c pUt));
}{%
}
\end{mapleinput}

\mapleresult
\begin{maplelatex}
\mapleinline{inert}{2d}{true,true;}{%
\[
\mathit{true},\;\mathit{true}
\]
}
\end{maplelatex}
\end{maplegroup}
\noindent
We check the SVD of $A,$ which is $A=U \Sigma \VT,$\footnote{The SVD of $A$ is not unique: For example, $A=(-U) \Sigma (-\VT)$ is another such factorization.} in the Clifford algebra language:
\begin{maplegroup}
\begin{mapleinput}
\mapleinline{active}{1d}{evalb(p=simplify(pU &c pSigma &c pVt)),}{%
}
\end{mapleinput}

\mapleresult
\begin{maplelatex}
\mapleinline{inert}{2d}{true;}{%
\[
\mathit{true},
\]
}
\end{maplelatex}

\end{maplegroup}
\noindent
where
\begin{maplegroup}
\begin{mapleinput}
\mapleinline{active}{1d}{'pU'=pU;}{%
}
\end{mapleinput}

\mapleresult
\begin{maplelatex}
\mapleinline{inert}{2d}{pU =
Id*(1/20*sqrt(10+2*sqrt(5))*sqrt(5)-1/8*sqrt(10-2*sqrt(5))-1/40*sqrt(1
0-2*sqrt(5))*sqrt(5))+(1/20*sqrt(10+2*sqrt(5))*sqrt(5)+1/8*sqrt(10-2*s
qrt(5))+1/40*sqrt(10-2*sqrt(5))*sqrt(5))*e1+(1/20*sqrt(10-2*sqrt(5))*s
qrt(5)+1/8*sqrt(10+2*sqrt(5))-1/40*sqrt(10+2*sqrt(5))*sqrt(5))*e2+(1/2
0*sqrt(10-2*sqrt(5))*sqrt(5)-1/8*sqrt(10+2*sqrt(5))+1/40*sqrt(10+2*sqr
t(5))*sqrt(5))*e12;}{%
\maplemultiline{
\mathit{pU}=\mathit{Id}\,({\displaystyle \frac {1}{20}} \,
\mathrm{\%1}\,\sqrt{5} - {\displaystyle \frac {1}{8}} \,\mathrm{
\%2} - {\displaystyle \frac {1}{40}} \,\mathrm{\%2}\,\sqrt{5}) + 
({\displaystyle \frac {1}{20}} \,\mathrm{\%1}\,\sqrt{5} + 
{\displaystyle \frac {1}{8}} \,\mathrm{\%2} + {\displaystyle 
\frac {1}{40}} \,\mathrm{\%2}\,\sqrt{5})\,\mathit{e1} \\
\mbox{} + ({\displaystyle \frac {1}{20}} \,\mathrm{\%2}\,\sqrt{5}
 + {\displaystyle \frac {1}{8}} \,\mathrm{\%1} - {\displaystyle 
\frac {1}{40}} \,\mathrm{\%1}\,\sqrt{5})\,\mathit{e2} + (
{\displaystyle \frac {1}{20}} \,\mathrm{\%2}\,\sqrt{5} - 
{\displaystyle \frac {1}{8}} \,\mathrm{\%1} + {\displaystyle 
\frac {1}{40}} \,\mathrm{\%1}\,\sqrt{5})\,\mathit{e12} \\
\mathrm{\%1} := \sqrt{10 + 2\,\sqrt{5}} \\
\mathrm{\%2} := \sqrt{10 - 2\,\sqrt{5}} }
}
\end{maplelatex}

\end{maplegroup}
\begin{maplegroup}
\begin{mapleinput}
\mapleinline{active}{1d}{'pSigma'=pSigma;}{%
}
\end{mapleinput}

\mapleresult
\begin{maplelatex}
\mapleinline{inert}{2d}{pSigma = sqrt(5)*Id+2*e1;}{%
\[
\mathit{pSigma}=\sqrt{5}\,\mathit{Id} + 2\,\mathit{e1}
\]
}
\end{maplelatex}

\end{maplegroup}
\begin{maplegroup}
\begin{mapleinput}
\mapleinline{active}{1d}{'pVt'=pVt;}{%
}
\end{mapleinput}

\mapleresult
\begin{maplelatex}
\mapleinline{inert}{2d}{pVt =
Id*(-1/20*sqrt(10-2*sqrt(5))*sqrt(5)+1/8*sqrt(10+2*sqrt(5))-1/40*sqrt(
10+2*sqrt(5))*sqrt(5))+(1/20*sqrt(10-2*sqrt(5))*sqrt(5)+1/8*sqrt(10+2*
sqrt(5))-1/40*sqrt(10+2*sqrt(5))*sqrt(5))*e1+(1/20*sqrt(10+2*sqrt(5))*
sqrt(5)+1/8*sqrt(10-2*sqrt(5))+1/40*sqrt(10-2*sqrt(5))*sqrt(5))*e2+(1/
20*sqrt(10+2*sqrt(5))*sqrt(5)-1/8*sqrt(10-2*sqrt(5))-1/40*sqrt(10-2*sq
rt(5))*sqrt(5))*e12;}{%
\maplemultiline{
\mathit{pVt}=\mathit{Id}\,( - {\displaystyle \frac {1}{20}} \,
\mathrm{\%1}\,\sqrt{5} + {\displaystyle \frac {1}{8}} \,\mathrm{
\%2} - {\displaystyle \frac {1}{40}} \,\mathrm{\%2}\,\sqrt{5}) + 
({\displaystyle \frac {1}{20}} \,\mathrm{\%1}\,\sqrt{5} + 
{\displaystyle \frac {1}{8}} \,\mathrm{\%2} - {\displaystyle 
\frac {1}{40}} \,\mathrm{\%2}\,\sqrt{5})\,\mathit{e1} \\
\mbox{} + ({\displaystyle \frac {1}{20}} \,\mathrm{\%2}\,\sqrt{5}
 + {\displaystyle \frac {1}{8}} \,\mathrm{\%1} + {\displaystyle 
\frac {1}{40}} \,\mathrm{\%1}\,\sqrt{5})\,\mathit{e2} + (
{\displaystyle \frac {1}{20}} \,\mathrm{\%2}\,\sqrt{5} - 
{\displaystyle \frac {1}{8}} \,\mathrm{\%1} - {\displaystyle 
\frac {1}{40}} \,\mathrm{\%1}\,\sqrt{5})\,\mathit{e12} \\
\mathrm{\%1} := \sqrt{10 - 2\,\sqrt{5}} \\
\mathrm{\%2} := \sqrt{10 + 2\,\sqrt{5}} }
}
\end{maplelatex}
\end{maplegroup}

\subsection{Singular Value Decomposition of a $3 \times 2$ matrix of rank $2$}
\label{sec:svd2}

In this section we will show our second example of SVD applied to a non-square matrix. The matrix  will need to be embedded first into an appropriate matrix algebra before its image can be found in a suitable Clifford algebra. 
\begin{maplegroup}
\begin{mapleinput}
\mapleinline{active}{1d}{C:=matrix(3,2,[3,0,0,-1,0,1]);
m:=rowdim(C): #number of rows of C is m
n:=coldim(C): #number of columns of C is n}{%
}
\end{mapleinput}

\mapleresult
\begin{maplelatex}
\mapleinline{inert}{2d}{C := matrix([[3, 0], [0, -1], [0, 1]]);}{%
\[
C :=  \left[ 
{\begin{array}{rr}
3 & 0 \\
0 & -1 \\
0 & 1
\end{array}}
 \right] 
\]
}
\end{maplelatex}
\end{maplegroup}
\noindent
Since our matrix $C$ is $3 \times 2,$ we will embed it into $\Mat(4,\BR) \simeq \cl_{p,q}$ where the signature $(p,q)$ could be either $(2,2)$ or $(3,1).$ Since the symbolic spinor representation of $\cl_{3,1}$ was already computed in \cite{Symbolic}, we will work with the signature $(3,1).$ 
\begin{maplegroup}
\begin{mapleinput}
\mapleinline{active}{1d}{dim:=4:B:=diag(1$3,-1$1):eval(makealiases(dim)):}{%
}
\end{mapleinput}
\end{maplegroup}
\noindent
Let's recall information about $\cl_{3,1}$ stored in {\sc CLIFFORD}\footnote{In the following display, $\hbox{'}cmulQ\hbox{'}(\frac12*Id+\frac12*e1,\frac12*Id+\frac12*e34)$ denotes an unevaluated product of two non-primitive idempotents $\frac12*Id+\frac12*e1$ and $\frac12*Id+\frac12*e34$ in $\cl_{3,1},$ and {\tt cmulQ} is a name of a simplified version of the procedure {\tt cmul}. While {\tt cmul} gives the Clifford product in $\cl(B),$ {\tt cmulQ} gives the Clifford product in $\cl(Q).$}:

\begin{maplegroup}
\begin{mapleinput}
\mapleinline{active}{1d}{data:=clidata();}{%
}
\end{mapleinput}
\mapleresult
\begin{maplelatex}
\mapleinline{inert}{2d}{data := [real, 4, simple, 'cmulQ'(1/2*Id+1/2*e1,1/2*Id+1/2*e34), [Id,
e2, e3, e23], [Id], [Id, e2, e3, e23]];}{%
\maplemultiline{
\mathit{data} := [\mathit{real}, \,4, \,\mathit{simple}, \,
\hbox{'}\mathit{cmulQ}\hbox{'}({\displaystyle \frac {1}{2}} \,
\mathit{Id} + {\displaystyle \frac {1}{2}} \,\mathit{e1}, \,
{\displaystyle \frac {1}{2}} \,\mathit{Id} + {\displaystyle 
\frac {1}{2}} \,\mathit{e34}), \,[\mathit{Id}, \,\mathit{e2}, \,
\mathit{e3}, \,\mathit{e23}], \,[\mathit{Id}],  \\
[\mathit{Id}, \,\mathit{e2}, \,\mathit{e3}, \,\mathit{e23}]] }
}
\end{maplelatex}

\end{maplegroup}
\noindent
We begin by defining a Grassmann basis in $\cl_{3,1},$ by assigning the fourth entry in the list to a primitive idempotent $f,$ and by generating a spinor basis in $S=\cl_{3,1}f.$
\begin{maplegroup}
\begin{mapleinput}
\mapleinline{active}{1d}{clibas:=cbasis(dim): #ordered basis in Cl(3,1)
f:=data[4]; #a primitive idempotent in Cl(3,1)
SBgens:=data[5]; #generators for a real basis in S
N:=data[2]: #dimension of the spinor representation }{%
}
\end{mapleinput}

\mapleresult
\begin{maplelatex}
\mapleinline{inert}{2d}{f := cmulQ(1/2*Id+1/2*e1,1/2*Id+1/2*e34);}{%
\[
f := \mathrm{cmulQ}({\displaystyle \frac {1}{2}} \,\mathit{Id} + 
{\displaystyle \frac {1}{2}} \,\mathit{e1}, \,{\displaystyle 
\frac {1}{2}} \,\mathit{Id} + {\displaystyle \frac {1}{2}} \,
\mathit{e34})
\]
}
\end{maplelatex}

\begin{maplelatex}
\mapleinline{inert}{2d}{SBgens := [Id, e2, e3, e23];}{%
\[
\mathit{SBgens} := [\mathit{Id}, \,\mathit{e2}, \,\mathit{e3}, \,
\mathit{e23}]
\]
}
\end{maplelatex}
\end{maplegroup}
\noindent
Thus, the real spinor basis in $S=\cl_{3,1}f$ consists of the following four polynomials: 

\begin{maplegroup}
\begin{mapleinput}
\mapleinline{active}{1d}{for i from 1 to N do f.i:=SBgens[i] &c f od;}{%
}
\end{mapleinput}

\mapleresult
\begin{maplelatex}
\mapleinline{inert}{2d}{f1 := 1/4*Id+1/4*e34+1/4*e1+1/4*e134;f2 := 1/4*e2+1/4*e234-1/4*e12-1/4*e1234;}{%
\[
\mathit{f1} := {\displaystyle \frac {1}{4}} \,\mathit{Id} + 
{\displaystyle \frac {1}{4}} \,\mathit{e34} + {\displaystyle 
\frac {1}{4}} \,\mathit{e1} + {\displaystyle \frac {1}{4}} \,
\mathit{e134},\quad \mathit{f2} := {\displaystyle \frac {1}{4}} \,\mathit{e2} + 
{\displaystyle \frac {1}{4}} \,\mathit{e234} - {\displaystyle 
\frac {1}{4}} \,\mathit{e12} - {\displaystyle \frac {1}{4}} \,
\mathit{e1234}
\]
}
\end{maplelatex}

\begin{maplelatex}
\mapleinline{inert}{2d}{
f3 := 1/4*e3+1/4*e4-1/4*e13-1/4*e14;f4 := 1/4*e23+1/4*e24+1/4*e123+1/4*e124;}{%
\[
\mathit{f3} := {\displaystyle \frac {1}{4}} \,\mathit{e3} + 
{\displaystyle \frac {1}{4}} \,\mathit{e4} - {\displaystyle 
\frac {1}{4}} \,\mathit{e13} - {\displaystyle \frac {1}{4}} \,
\mathit{e14}, \quad \mathit{f4} := {\displaystyle \frac {1}{4}} \,\mathit{e23} + 
{\displaystyle \frac {1}{4}} \,\mathit{e24} + {\displaystyle 
\frac {1}{4}} \,\mathit{e123} + {\displaystyle \frac {1}{4}} \,
\mathit{e124}
\]
}
\end{maplelatex}
\end{maplegroup}
\noindent
We compute matrices $M_1,\ldots,M_{16}$ representing each basis element in $\cl_{3,1}.$ These matrices can be computed also with the help of a procedure {\tt matKrepr} which, being less general than {\tt spinorKrepr} used earlier, is also simpler to use. We won't display these matrices since they can be found in \cite{Symbolic}.
\begin{maplegroup}
\begin{mapleinput}
\mapleinline{active}{1d}{
for i from 1 to nops(clibas) do M[i]:=subs(Id=1,matKrepr(clibas[i]))od:}{%
}
\end{mapleinput}
\end{maplegroup}
\noindent
Before we can use the procedure {\tt phi} that realizes the isomorphism 
$\Mat(4,\BR) \stackrel{\varphi}{\simeq} \cl_{3,1},$ we need to embed matrix $C$ into $\Mat(4,\BR).$ This is accomplished with a procedure {\tt embed}. 
\begin{maplegroup}
\begin{mapleinput}
\mapleinline{active}{1d}{A:=embed(C);m1:=rowdim(A):n1:=coldim(A):}{%
}
\end{mapleinput}

\mapleresult
\begin{maplelatex}
\mapleinline{inert}{2d}{A := matrix([[3, 0, 0, 0], [0, -1, 0, 0], [0, 1, 0, 0], [0, 0, 0,
0]]);}{%
\[
A :=  \left[ 
{\begin{array}{rrrr}
3 & 0 & 0 & 0 \\
0 & -1 & 0 & 0 \\
0 & 1 & 0 & 0 \\
0 & 0 & 0 & 0
\end{array}}
 \right] 
\]
}
\end{maplelatex}
\end{maplegroup}
\noindent

\begin{maplegroup}
\begin{mapleinput}
\mapleinline{active}{1d}{p:=phi(A,M); #finding image of A (and C) in Cl(3,1)
pT:=phi(t(A),M); #finding image of AT (and CT) in Cl(3,1)}{%
}
\end{mapleinput}

\mapleresult
\begin{maplelatex}
\mapleinline{inert}{2d}{p := 1/2*Id+e1-1/4*e23-1/4*e24+1/2*e34+1/4*e123+1/4*e124+e134;}{%
\[
p := {\displaystyle \frac {1}{2}} \,\mathit{Id} + \mathit{e1} - 
{\displaystyle \frac {1}{4}} \,\mathit{e23} - {\displaystyle 
\frac {1}{4}} \,\mathit{e24} + {\displaystyle \frac {1}{2}} \,
\mathit{e34} + {\displaystyle \frac {1}{4}} \,\mathit{e123} + 
{\displaystyle \frac {1}{4}} \,\mathit{e124} + \mathit{e134}
\]
}
\end{maplelatex}

\begin{maplelatex}
\mapleinline{inert}{2d}{pT := 1/2*Id+e1+1/4*e23-1/4*e24+1/2*e34-1/4*e123+1/4*e124+e134;}{%
\[
\mathit{pT} := {\displaystyle \frac {1}{2}} \,\mathit{Id} + 
\mathit{e1} + {\displaystyle \frac {1}{4}} \,\mathit{e23} - 
{\displaystyle \frac {1}{4}} \,\mathit{e24} + {\displaystyle 
\frac {1}{2}} \,\mathit{e34} - {\displaystyle \frac {1}{4}} \,
\mathit{e123} + {\displaystyle \frac {1}{4}} \,\mathit{e124} + 
\mathit{e134}
\]
}
\end{maplelatex}
\end{maplegroup}
Next we compute $\ATA,\AAT,$ their images in the Clifford algebra $\cl_{3,1}$ under $\varphi,$ their eigenvalues, and orthonormal eigenvectors $v_i$ and $u_i,$ $i=1,\ldots,4,$ which will become columns of two orthogonal matrices $V$ and $U.$
\begin{maplegroup}
\begin{mapleinput}
\mapleinline{active}{1d}{ATA,AAT:=evalm(t(A) &* A),evalm(A &* t(A));}{%
}
\end{mapleinput}

\mapleresult
\begin{maplelatex}
\mapleinline{inert}{2d}{ATA, AAT := matrix([[9, 0, 0, 0], [0, 2, 0, 0], [0, 0, 0, 0], [0, 0,
0, 0]]), matrix([[9, 0, 0, 0], [0, 1, -1, 0], [0, -1, 1, 0], [0, 0, 0,
0]]);}{%
\[
\mathit{ATA}, \,\mathit{AAT} :=  \left[ 
{\begin{array}{rrrr}
9 & 0 & 0 & 0 \\
0 & 2 & 0 & 0 \\
0 & 0 & 0 & 0 \\
0 & 0 & 0 & 0
\end{array}}
 \right] , \, \left[ 
{\begin{array}{rrrr}
9 & 0 & 0 & 0 \\
0 & 1 & -1 & 0 \\
0 & -1 & 1 & 0 \\
0 & 0 & 0 & 0
\end{array}}
 \right] 
\]
}
\end{maplelatex}

\end{maplegroup}
\begin{maplegroup}
\begin{mapleinput}
\mapleinline{active}{1d}{
pTp,ppT:=phi(ATA,M),phi(AAT,M); #finding images of ATA and AAT in Cl(3,1)
}{%
}
\end{mapleinput}

\mapleresult
\begin{maplelatex}
\mapleinline{inert}{2d}{pTp, ppT := 11/4*Id+7/4*e1+11/4*e34+7/4*e134,
11/4*Id+7/4*e1+1/2*e24+9/4*e34-1/2*e124+9/4*e134;}{%
\maplemultiline{
\mathit{pTp}, \,\mathit{ppT} :=  \\
{\displaystyle \frac {11}{4}} \,\mathit{Id} + {\displaystyle 
\frac {7}{4}} \,\mathit{e1} + {\displaystyle \frac {11}{4}} \,
\mathit{e34} + {\displaystyle \frac {7}{4}} \,\mathit{e134}, \,
{\displaystyle \frac {11}{4}} \,\mathit{Id} + {\displaystyle 
\frac {7}{4}} \,\mathit{e1} + {\displaystyle \frac {1}{2}} \,
\mathit{e24} + {\displaystyle \frac {9}{4}} \,\mathit{e34} - 
{\displaystyle \frac {1}{2}} \,\mathit{e124} + {\displaystyle 
\frac {9}{4}} \,\mathit{e134} }
}
\end{maplelatex}

\end{maplegroup}
\noindent
The characteristic polynomial of $\ATA$ and $\AAT$ is $(x-9)(x-2)x^2$ while the minimal polynomial of $p\,Tp$ and $ppT$ is $x(x-2)(x-9).$ 
\begin{maplegroup}
\begin{mapleinput}
\mapleinline{active}{1d}{P1:=assignL(sort([eigenvects(ATA)],byeigenvals));
P2:=assignL(sort([eigenvects(AAT)],byeigenvals));
for i from 1 to N do lambda.i:=P1[2][i];
\phtab{xxx} v.i:=map(simplify,normalize(P1[3][i]));
\phtab{xxx} u.i:=map(simplify,normalize(P2[3][i]));
od:}{%
}
\end{mapleinput}

\mapleresult
\begin{maplelatex}
\mapleinline{inert}{2d}{P1 := [4, [9, 2, 0, 0], [vector([1, 0, 0, 0]), vector([0, 1, 0, 0]),
vector([0, 0, 1, 0]), vector([0, 0, 0, 1])]];}{%
\[
\mathit{P1} := [4, \,[9, \,2, \,0, \,0], \,[[1, \,0, \,0, \,0], 
\,[0, \,1, \,0, \,0], \,[0, \,0, \,1, \,0], \,[0, \,0, \,0, \,1]]
]
\]
}
\end{maplelatex}

\begin{maplelatex}
\mapleinline{inert}{2d}{P2 := [4, [9, 2, 0, 0], [vector([1, 0, 0, 0]), vector([0, -1, 1, 0]),
vector([0, 1, 1, 0]), vector([0, 0, 0, 1])]];}{%
\[
\mathit{P2} := [4, \,[9, \,2, \,0, \,0], \,[[1, \,0, \,0, \,0], 
\,[0, \,-1, \,1, \,0], \,[0, \,1, \,1, \,0], \,[0, \,0, \,0, \,1]
]]
\]
}
\end{maplelatex}

\end{maplegroup}
\noindent
The third list in $P1$ and the third list in $P2$ contain vectors $v_i$ and $u_i$ that will make up matrices $V$ and $U$ respectively. Clearly, since $V$ is the $4 \times 4$ identity matrix, $\varphi(V)=\varphi(\VT)=\Id.$ In order to translate this matrix picture $\cl_{3,1},$ we need to find spinors $sv_i=\varphi(v_i)$ and $su_i=\varphi(u_i)$ in $S=\cl_{3,1}f.$
\begin{maplegroup}
\begin{mapleinput}
\mapleinline{active}{1d}{spinorbasis:=[''f1'',''f2'',''f3'',''f4'']:}{%
}
\mapleinline{active}{1d}{for i from 1 to N do sv.i:=convert(v.i,spinor,spinorbasis);
\phtab{xxxxxxxxxxxxxxxxxx} su.i:=convert(u.i,spinor,spinorbasis) od;
}{%
}
\end{mapleinput}

\mapleresult
\begin{maplelatex}
\mapleinline{inert}{2d}{sv1 := f1;sv2 := f2;sv3 := f3;sv4 := f4;}{%
\[
\mathit{sv1} := \mathit{f1}, \; \mathit{sv2} := \mathit{f2}, \; 
\mathit{sv3} := \mathit{f3}, \;\mathit{sv4} := \mathit{f4}
\]
}
\end{maplelatex}

\begin{maplelatex}
\mapleinline{inert}{2d}{su1 := f1;su2 := -1/2*sqrt(2)*f2+1/2*sqrt(2)*f3;
su3 := 1/2*sqrt(2)*f2+1/2*sqrt(2)*f3;su4 := f4;}{%
\[
\mathit{su1} := \mathit{f1}, \;\mathit{su2} :=  - {\displaystyle \frac {1}{2}} \,\sqrt{2}\,
\mathit{f2} + {\displaystyle \frac {1}{2}} \,\sqrt{2}\,\mathit{f3}, \; \mathit{su3} := {\displaystyle \frac {1}{2}} \,\sqrt{2}\,\mathit{f2} + {\displaystyle \frac {1}{2}} \,\sqrt{2}\,\mathit{f3}, \;\mathit{sv4} := \mathit{f4}
\]
}
\end{maplelatex}
\end{maplegroup}
\noindent
Spinors $sv_i$ are eigenspinors of $p\,Tp,$ while $su_i$ are eigenspinors of $ppT$ with the eigenvalues~$\lambda_1=9,$ $\lambda_2=2,\lambda_3=0,\lambda_4=0,$ namely:\begin{maplegroup}
\begin{mapleinput}
\mapleinline{active}{1d}{for i from 1 to N do simplify((pTp - lambda.i) &c sv.i);
\phtab{xxxxxxxxxxxxxxxxxx} simplify((ppT - lambda.i) &c su.i) od;}{%
}
\end{mapleinput}

\mapleresult
\begin{maplelatex}
\mapleinline{inert}{2d}{0;0;0;0;0;0;0;0;}{%
\[
0,\,0,\,0,\,0,\,0,\,0,\,0,\,0
\]
}
\end{maplelatex}
\end{maplegroup}
\begin{maplegroup}
\begin{mapleinput}
\mapleinline{active}{1d}{
V:=radsimplify(augment(v.(1..N))): #defining matrix V (identity matrix)
pV,pVt:=phi(V,M),phi(t(V),M): #finding images of V and t(V) in Cl(3,1)
}{%
}
\end{mapleinput}
\end{maplegroup}
\noindent
We check that $A v_i = \sigma_i u_i$ where $\sigma_i=\sqrt{\lambda_i}$ by verifying this fact in $\cl_{3,1}.$\footnote{As before, we are following engineering practice \cite{Maciejowski} of considering square roots of zero eigenvalues of $\AAT$ and $\ATA$ also as the singular values of a matrix. While it is  convenient to do so for computational reasons, it is not conceptually correct since by the definition singular values of a matrix are always positive.}
\begin{maplegroup}
\begin{mapleinput}
\mapleinline{active}{1d}{for i from 1 to N do sigma.i:=sqrt(lambda.i) od;}{%
}
\end{mapleinput}

\mapleresult
\begin{maplelatex}
\mapleinline{inert}{2d}{sigma1 := 3;sigma2 := sqrt(2);sigma3 := 0;sigma4 := 0;}{%
\[
\sigma 1 := 3, \; \sigma 2 := \sqrt{2}, \; \sigma 3 := 0, \; \sigma 4 := 0
\]
}
\end{maplelatex}

\end{maplegroup}
\begin{maplegroup}
\begin{mapleinput}
\mapleinline{active}{1d}{for i from 1 to N do simplify(p &c sv.i-sigma.i*su.i) od; }{%
}
\end{mapleinput}

\mapleresult
\begin{maplelatex}
\mapleinline{inert}{2d}{0;0;0;0;}{%
\[
0,\;0,\;0,\;0
\]
}
\end{maplelatex}
\end{maplegroup}
\noindent
Thus, we may now define an orthogonal matrix $U=[u_1|u_2|u_3|u_4]$ and find its image $\varphi(U)$ and the image of its transpose in $\cl_{3,1}:$
\begin{maplegroup}
\begin{mapleinput}
\mapleinline{active}{1d}{U:=radsimplify(augment(seq(u.i,i=1..N))); #defining matrix U}{%
}
\end{mapleinput}

\mapleresult
\begin{maplelatex}
\mapleinline{inert}{2d}{U := matrix([[1, 0, 0, 0], [0, -1/2*sqrt(2), 1/2*sqrt(2), 0], [0,
1/2*sqrt(2), 1/2*sqrt(2), 0], [0, 0, 0, 1]]);}{%
\[
U :=  \left[ 
{\begin{array}{rccr}
1 & 0 & 0 & 0 \\
0 &  - {\displaystyle \frac {1}{2}} \,\sqrt{2} & {\displaystyle 
\frac {1}{2}} \,\sqrt{2} & 0 \\ [2ex]
0 & {\displaystyle \frac {1}{2}} \,\sqrt{2} & {\displaystyle 
\frac {1}{2}} \,\sqrt{2} & 0 \\ [2ex]
0 & 0 & 0 & 1
\end{array}}
 \right] 
\]
}
\end{maplelatex}

\end{maplegroup}
\begin{maplegroup}
\begin{mapleinput}
\mapleinline{active}{1d}{pU:=phi(U,M);pUt:=phi(t(U),M);}{%
}
\end{mapleinput}

\mapleresult
\begin{maplelatex}
\mapleinline{inert}{2d}{pU :=
1/2*Id+1/2*e1-1/4*sqrt(2)*e24-1/4*sqrt(2)*e34+1/4*sqrt(2)*e124+1/4*sqr
t(2)*e134;}{%
\[
\mathit{pU} := {\displaystyle \frac {1}{2}} \,\mathit{Id} + 
{\displaystyle \frac {1}{2}} \,\mathit{e1} - {\displaystyle 
\frac {1}{4}} \,\sqrt{2}\,\mathit{e24} - {\displaystyle \frac {1
}{4}} \,\sqrt{2}\,\mathit{e34} + {\displaystyle \frac {1}{4}} \,
\sqrt{2}\,\mathit{e124} + {\displaystyle \frac {1}{4}} \,\sqrt{2}
\,\mathit{e134}
\]
}
\end{maplelatex}

\begin{maplelatex}
\mapleinline{inert}{2d}{pUt :=
1/2*Id+1/2*e1-1/4*sqrt(2)*e24-1/4*sqrt(2)*e34+1/4*sqrt(2)*e124+1/4*sqr
t(2)*e134;}{%
\[
\mathit{pUt} := {\displaystyle \frac {1}{2}} \,\mathit{Id} + 
{\displaystyle \frac {1}{2}} \,\mathit{e1} - {\displaystyle 
\frac {1}{4}} \,\sqrt{2}\,\mathit{e24} - {\displaystyle \frac {1
}{4}} \,\sqrt{2}\,\mathit{e34} + {\displaystyle \frac {1}{4}} \,
\sqrt{2}\,\mathit{e124} + {\displaystyle \frac {1}{4}} \,\sqrt{2}
\,\mathit{e134}
\]
}
\end{maplelatex}

\end{maplegroup}
\noindent
The fact that $U$ is orthogonal is reflected in the following:
\begin{maplegroup}
\begin{mapleinput}
\mapleinline{active}{1d}{pUt &c pU;}{%
}
\end{mapleinput}

\mapleresult
\begin{maplelatex}
\mapleinline{inert}{2d}{Id;}{%
\[
\mathit{Id}
\]
}
\end{maplelatex}
\end{maplegroup}
\noindent
Finally, it is just enough to find matrix $\Sigma$ and verify SVD for $\varphi(A)$ in the Clifford algebra $\cl_{3,1}.$
\begin{maplegroup}
\begin{mapleinput}
\mapleinline{active}{1d}{Sigma:=makediag(m1,n1,[seq(sigma.i,i=1..N)]);pSigma:=phi(Sigma,M);}{%
}
\end{mapleinput}

\mapleresult
\begin{maplelatex}
\mapleinline{inert}{2d}{Sigma := matrix([[3, 0, 0, 0], [0, sqrt(2), 0, 0], [0, 0, 0, 0], [0,
0, 0, 0]]);}{%
\[
\Sigma  :=  \left[ 
{\begin{array}{rcrr}
3 & 0 & 0 & 0 \\
0 & \sqrt{2} & 0 & 0 \\
0 & 0 & 0 & 0 \\
0 & 0 & 0 & 0
\end{array}}
 \right] 
\]
}
\end{maplelatex}

\begin{maplelatex}
\mapleinline{inert}{2d}{pSigma :=
Id*(3/4+1/4*sqrt(2))+(-1/4*sqrt(2)+3/4)*e1+(3/4+1/4*sqrt(2))*e34+(-1/4
*sqrt(2)+3/4)*e134;}{%
\[
\mathit{pSigma} := \mathit{Id}\,({\displaystyle \frac {3}{4}}  + 
{\displaystyle \frac {1}{4}} \,\sqrt{2}) + ( - {\displaystyle 
\frac {1}{4}} \,\sqrt{2} + {\displaystyle \frac {3}{4}} )\,
\mathit{e1} + ({\displaystyle \frac {3}{4}}  + {\displaystyle 
\frac {1}{4}} \,\sqrt{2})\,\mathit{e34} + ( - {\displaystyle 
\frac {1}{4}} \,\sqrt{2} + {\displaystyle \frac {3}{4}} )\,
\mathit{e134}
\]
}
\end{maplelatex}
\end{maplegroup}
\begin{maplegroup}
\begin{mapleinput}
\mapleinline{active}{1d}{evalb(p=pU &c pSigma &c pVt); #SVD of phi(A)}{%
}
\end{mapleinput}

\mapleresult
\begin{maplelatex}
\mapleinline{inert}{2d}{true;}{%
\[
\mathit{true}
\]
}
\end{maplelatex}

\end{maplegroup}
\noindent
Since the original matrix $C$ was $3$ by $2$ and not $4$ by $4,$ in order to find SVD of $C$ we need to project out certain columns and rows out of the matrices $U,\Sigma,\mbox{ and }V.$ This can also be done internally in the Clifford algebra. The original matrix $C$ has therefore this factorization:
\begin{maplegroup}
\begin{mapleinput}
\mapleinline{active}{1d}{U2,Sigma2,V2t:=submatrix(U,1..m,1..m),submatrix(Sigma,1..m,1..n),
submatrix(t(V),1..n,1..n);}{%
}
\end{mapleinput}

\mapleresult
\begin{maplelatex}
\mapleinline{inert}{2d}{U2, Sigma2, V2t := matrix([[1, 0, 0], [0, -1/2*sqrt(2), 1/2*sqrt(2)],
[0, 1/2*sqrt(2), 1/2*sqrt(2)]]), matrix([[3, 0], [0, sqrt(2)], [0,
0]]), matrix([[1, 0], [0, 1]]);}{%
\[
\mathit{U2}, \,\Sigma 2, \,\mathit{V2t} :=  \left[ 
{\begin{array}{rcc}
1 & 0 & 0 \\
0 &  - {\displaystyle \frac {1}{2}} \,\sqrt{2} & {\displaystyle 
\frac {1}{2}} \,\sqrt{2} \\ [2ex]
0 & {\displaystyle \frac {1}{2}} \,\sqrt{2} & {\displaystyle 
\frac {1}{2}} \,\sqrt{2}
\end{array}}
 \right] , \, \left[ 
{\begin{array}{rc}
3 & 0 \\
0 & \sqrt{2} \\
0 & 0
\end{array}}
 \right] , \, \left[ 
{\begin{array}{rr}
1 & 0 \\
0 & 1
\end{array}}
 \right] 
\]
}
\end{maplelatex}

\end{maplegroup}
\begin{maplegroup}
\begin{mapleinput}
\mapleinline{active}{1d}{evalm(C)=radsimplify(evalm(U2 &* Sigma2 &* V2t));}{%
}
\end{mapleinput}

\mapleresult
\begin{maplelatex}
\mapleinline{inert}{2d}{matrix([[3, 0], [0, -1], [0, 1]]) = matrix([[3, 0], [0, -1], [0,
1]]);}{%
\[
 \left[ 
{\begin{array}{rr}
3 & 0 \\
0 & -1 \\
0 & 1
\end{array}}
 \right] = \left[ 
{\begin{array}{rr}
3 & 0 \\
0 & -1 \\
0 & 1
\end{array}}
 \right] 
\]
}
\end{maplelatex}

\end{maplegroup}

\subsection{Additional comments}
\label{sec:comments}

In this section we have shown that it is possible to translate the matrix algebra picture of the Singular Value Decomposition of a matrix $A$ into the Clifford algebra language. Although we have not abandoned entirely the linear algebra formalism in our examples, e.g., we have computed the eigenvalues and the eigenvectors of $\ATA$ and $\AAT,$ and only then we have found images of the eigenvectors in the spinor space $\cl_{p,q}f,$ these computations including solving the eigenvalue problems can be done entirely in $\cl_{p,q}$ without using matrices. For example, if we consider element $p\,Tp = \varphi(\ATA)$ in $\cl_{3,1}$ from the last example, we can find its eigenvalues from its minimal polynomial while its eigenvectors can be found by solving the eigenvalue equation. For example, we know that the minimal polynomial of $p\,Tp$ is
\begin{maplegroup}
\begin{mapleinput}
\mapleinline{active}{1d}{pol:=factor(climinpoly(pTp));}{%
}
\end{mapleinput}

\mapleresult
\begin{maplelatex}
\mapleinline{inert}{2d}{pol := x*(x-2)*(x-9);}{%
\[
\mathit{pol} := x\,(x - 2)\,(x - 9)
\]
}
\end{maplelatex}

\end{maplegroup}
\noindent
hence the eigenvalues of $p\,Tp$ are $9,2,0.$  Of course, the minimal polynomial doesn't give us their geometric multiplicities. However, we can find them by solving the eigenvalue equation directly in $\cl_{3,1}$ for each of the eigenvalues. Let's assign these three known eigenvalues to $\lambda_1,\lambda_2,\lambda_3:$
\begin{maplegroup}
\begin{mapleinput}
\mapleinline{active}{1d}{lambda1,lambda2,lambda3:=9,2,0;}{%
}
\end{mapleinput}

\mapleresult
\begin{maplelatex}
\mapleinline{inert}{2d}{lambda1, lambda2, lambda3 := 9, 2, 0;}{%
\[
\lambda 1, \,\lambda 2, \,\lambda 3 := 9, \,2, \,0
\]
}
\end{maplelatex}
\end{maplegroup}
\noindent
Let $\psi \in \cl_{3,1}f$ be an arbitrary spinor expressed in the spinor basis $\{f_1,f_2,f_3,f_4\}$ computed earlier, and let $c_1,\ldots,c_4$ be its undefined (real) coefficients.
\begin{maplegroup}
\begin{mapleinput}
\mapleinline{active}{1d}{psi:=c[1]*'f1'+c[2]*'f2'+c[3]*'f3'+c[4]*'f4';}{%
}
\end{mapleinput}

\mapleresult
\begin{maplelatex}
\mapleinline{inert}{2d}{psi := c[1]*f1+c[2]*f2+c[3]*f3+c[4]*f4;}{%
\[
\psi  := {c_{1}}\,\mathit{f1} + {c_{2}}\,\mathit{f2} + {c_{3}}\,
\mathit{f3} + {c_{4}}\,\mathit{f4}
\]
}
\end{maplelatex}
\end{maplegroup}
\noindent
We will now use the procedure {\tt clisolve2} from the {\sc suppl} package (see Section 2) to solve the eigenvalue equation 
\begin{equation}
p\,Tp \, \psi = \lambda \psi.
\label{eq:eigen}
\end{equation}
\noindent
First we solve it when $\lambda=\lambda_1=9.$  
\begin{maplegroup}
\begin{mapleinput}
\mapleinline{active}{1d}{sol1:=clisolve2(pTp &c psi - lambda1*psi,[c[1],c[2],c[3],c[4]]);}{%
}
\end{mapleinput}

\mapleresult
\begin{maplelatex}
\mapleinline{inert}{2d}{sol1 := [\{c[2] = 0, c[3] = 0, c[4] = 0, c[1] = c[1]\}];}{%
\[
\mathit{sol1} := [\{{c_{2}}=0, \,{c_{3}}=0, \,{c_{4}}=0, \,{c_{1}
}={c_{1}}\}]
\]
}
\end{maplelatex}
\end{maplegroup}
\begin{maplegroup}
\begin{mapleinput}
\mapleinline{active}{1d}{psi1:=c[1]*'f1'; #eigenspinor of pTp with eigenvalue lambda1=9}{%
}
\end{mapleinput}

\mapleresult
\begin{maplelatex}
\mapleinline{inert}{2d}{psi1 := c[1]*f1;}{%
\[
\psi 1 := {c_{1}}\,\mathit{f1}
\]
}
\end{maplelatex}
\end{maplegroup}
\noindent
Thus, the first solution to (\ref{eq:eigen}) for $\lambda_1=9$ is a one-parameter solution that belongs to a one-dimensional subspace spanned by $f_1.$ That is, the geometric multiplicity of the eigenvalue $\lambda_1=9$ is $1.$ Similarly for $\lambda=\lambda_2=2:$
\begin{maplegroup}
\begin{mapleinput}
\mapleinline{active}{1d}{sol2:=clisolve2(pTp &c psi - lambda2*psi,[c[1],c[2],c[3],c[4]]);}{%
}
\end{mapleinput}

\mapleresult
\begin{maplelatex}
\mapleinline{inert}{2d}{sol2 := [\{c[3] = 0, c[4] = 0, c[1] = 0, c[2] = c[2]\}];}{%
\[
\mathit{sol2} := [\{{c_{3}}=0, \,{c_{4}}=0, \,{c_{1}}=0, \,{c_{2}
}={c_{2}}\}]
\]
}
\end{maplelatex}

\end{maplegroup}
\begin{maplegroup}
\begin{mapleinput}

\mapleinline{active}{1d}{psi2:=c[2]*'f2'; #eigenspinor of pTp with eigenvalue lambda2=2}{%
}
\end{mapleinput}

\mapleresult
\begin{maplelatex}
\mapleinline{inert}{2d}{psi2 := c[2]*f2;}{%
\[
\psi 2 := {c_{2}}\,\mathit{f2}
\]
}
\end{maplelatex}
\end{maplegroup}
\noindent
The second solution to (\ref{eq:eigen}) for $\lambda_2=2$ is also a one-parameter solution that belongs to a one-dimensional subspace spanned by $f_2.$ The geometric multiplicity of the eigenvalue $\lambda_2=2$ is also~$1.$
\begin{maplegroup}
\begin{mapleinput}
\mapleinline{active}{1d}{sol34:=clisolve2(pTp &c psi - lambda3*psi,[c[1],c[2],c[3],c[4]]);}{%
}
\end{mapleinput}

\mapleresult
\begin{maplelatex}
\mapleinline{inert}{2d}{sol34 := [\{c[3] = c[3], c[4] = c[4], c[2] = 0, c[1] = 0\}];}{%
\[
\mathit{sol34} := [\{{c_{3}}={c_{3}}, \,{c_{4}}={c_{4}}, \,{c_{2}
}=0, \,{c_{1}}=0\}]
\]
}
\end{maplelatex}

\end{maplegroup}
\noindent
The third solution to (\ref{eq:eigen}) for $\lambda=\lambda_3=0$ is parameterized by two parameters $c_3$ and $c_4$ and belongs therefore to a two-dimensional subspace of $S$ spanned by $\{f_3,f_4\}.$ This implies that the geometric multiplicity of $\lambda=\lambda3=0$ is two.
\begin{maplegroup}
\begin{mapleinput}
\mapleinline{active}{1d}{psi3:=c[3]*'f3'; #eigenspinor of pTp with eigenvalue lambda3=0}{%
}
\end{mapleinput}

\mapleresult
\begin{maplelatex}
\mapleinline{inert}{2d}{psi3 := c[3]*f3;}{%
\[
\psi 3 := {c_{3}}\,\mathit{f3}
\]
}
\end{maplelatex}

\end{maplegroup}
\begin{maplegroup}
\begin{mapleinput}
\mapleinline{active}{1d}{psi4:=c[4]*'f4'; #eigenspinor of pTp with eigenvalue lambda3=0}{%
}
\end{mapleinput}

\mapleresult
\begin{maplelatex}
\mapleinline{inert}{2d}{psi4 := c[4]*f4;}{%
\[
\psi 4 := {c_{4}}\,\mathit{f4}
\]
}
\end{maplelatex}
\end{maplegroup}
\noindent
So, up to a normalizing scalar, the eigenspinors $\psi_1,\psi_2,\psi_3,\psi_4$ give previously computed normalized spinors $sv_1,sv_2,sv_3,sv_4.$ Likewise for the image $ppT$ of $\AAT:$
\begin{maplegroup}
\begin{mapleinput}
\mapleinline{active}{1d}{sol1:=clisolve2(ppT &c psi - lambda1*psi,[c[1],c[2],c[3],c[4]]);}{%
}
\end{mapleinput}

\mapleresult
\begin{maplelatex}
\mapleinline{inert}{2d}{sol1 := [\{c[2] = 0, c[3] = 0, c[4] = 0, c[1] = c[1]\}];}{%
\[
\mathit{sol1} := [\{{c_{2}}=0, \,{c_{3}}=0, \,{c_{4}}=0, \,{c_{1}}={c_{1}}\}]
\]
}
\end{maplelatex}

\end{maplegroup}
\begin{maplegroup}
\begin{mapleinput}
\mapleinline{active}{1d}{phi1:=c[1]*'f1'; #eigenspinor of ppT with eigenvalue lambda1=9}{%
}
\end{mapleinput}

\mapleresult
\begin{maplelatex}
\mapleinline{inert}{2d}{phi1 := c[1]*f1;}{%
\[
\phi 1 := {c_{1}}\,\mathit{f1}
\]
}
\end{maplelatex}

\end{maplegroup}
\noindent
The first one-parameter eigenspace of $ppT$ corresponding to the eigenvalue $\lambda_1=9$ is  one-dimensional and it is spanned by~$f_1.$
\begin{maplegroup}
\begin{mapleinput}
\mapleinline{active}{1d}{sol2:=clisolve2(ppT &c psi - lambda2*psi,[c[1],c[2],c[3],c[4]]);}{%
}
\end{mapleinput}

\mapleresult
\begin{maplelatex}
\mapleinline{inert}{2d}{sol2 := [\{c[3] = c[3], c[4] = 0, c[1] = 0, c[2] = -c[3]\}];}{%
\[
\mathit{sol2} := [\{{c_{3}}={c_{3}}, \,{c_{4}}=0, \,{c_{1}}=0, \,
{c_{2}}= - {c_{3}}\}]
\]
}
\end{maplelatex}

\end{maplegroup}
\begin{maplegroup}
\begin{mapleinput}
\mapleinline{active}{1d}{phi2:=-c[3]*'f2'+c[3]*'f3'; #eigenspinor of ppT with eigenvalue lambda2=2
}{%
}
\end{mapleinput}

\mapleresult
\begin{maplelatex}
\mapleinline{inert}{2d}{phi2 := -c[3]*f2+c[3]*f3;}{%
\[
\phi 2 :=  - {c_{3}}\,\mathit{f2} + {c_{3}}\,\mathit{f3}
\]
}
\end{maplelatex}

\end{maplegroup}
\noindent
The second one-parameter eigenspace of $ppT$ corresponding to the eigenvalue $\lambda_2=0$ is also one-dimensional and it is spanned by $\{f_2,f_3\}.$
\begin{maplegroup}
\begin{mapleinput}
\mapleinline{active}{1d}{sol34:=clisolve2(ppT &c psi - lambda3*psi,[c[1],c[2],c[3],c[4]]);}{%
}
\end{mapleinput}

\mapleresult
\begin{maplelatex}
\mapleinline{inert}{2d}{sol34 := [\{c[3] = c[3], c[4] = c[4], c[1] = 0, c[2] = c[3]\}];}{%
\[
\mathit{sol34} := [\{{c_{3}}={c_{3}}, \,{c_{4}}={c_{4}}, \,{c_{1}
}=0, \,{c_{2}}={c_{3}}\}]
\]
}
\end{maplelatex}

\end{maplegroup}
\noindent
The third solution to (\ref{eq:eigen}) is parameterized by two parameters $c_3$ and $c_4.$ Thus, we get a two dimensional vectors space spanned by $\{f_3,f_4\}.$
\begin{maplegroup}
\begin{mapleinput}
\mapleinline{active}{1d}{phi3:=c[3]*'f2'+c[3]*'f3'; #eigenspinor of ppT with eigenvalue lambda3=0
}{%
}
\end{mapleinput}

\mapleresult
\begin{maplelatex}
\mapleinline{inert}{2d}{phi3 := c[3]*f2+c[3]*f3;}{%
\[
\phi 3 := {c_{3}}\,\mathit{f2} + {c_{3}}\,\mathit{f3}
\]
}
\end{maplelatex}

\end{maplegroup}
\begin{maplegroup}
\begin{mapleinput}
\mapleinline{active}{1d}{phi4:=c[4]*'f4'; #eigenspinor of ppT with eigenvalue lambda3=0}{%
}
\end{mapleinput}

\mapleresult
\begin{maplelatex}
\mapleinline{inert}{2d}{phi4 := c[4]*f4;}{%
\[
\phi 4 := {c_{4}}\,\mathit{f4}
\]
}
\end{maplelatex}
\end{maplegroup}
\noindent
So, again up to a normalizing scalar, $\phi_1,\phi_2,\phi_3,\phi_4$ give the eigenspinors $su_1,su_2,su_3,su_4$ computed earlier.
 
\section{Clifford algebras in robotics}
\label{sec:robotics}

Clifford algebras $\cl(V,Q)$ on a quadratic space $(V,Q)$ endowed with a {\it degenerate\/} quadratic form $Q$ and associated groups $\Spin,\,\Pin,$ Clifford, etc., were studied in \cite{Deformation,AblLoun85} and \cite{Structure,Brooke7880,Brooke85}. In contrast to the Clifford algebras of a non-degenerate quadratic form, these algebras possess a non-trivial two-sided nilpotent ideal called {\it Jacobson radical\/}. The Jacobson radical $J$ is generated by the null-vectors in $V$ which are orthogonal to the entire space $V$ (that is, $J$ is generated by the orthogonal complement $V^\bot$ of $V).$ It is known \cite{CurtisReiner,Herstein} that $J$ contains every nilpotent left and right ideal in $\cl(V,Q).$ From the point of view of the spinorial representation theory of Clifford algebras used in Section \ref{sec:svd}, an important difference is that $\cl(V,Q)$ does not possess faithful matrix representation when $Q$ is degenerate.

Let $V=V' \bot V^\bot$ where $V'$ is endowed with a non-degenerate part $Q'$ of $Q$ of signature $(p,q).$ Let $\dim(V^{\bot})=d,$ hence $p+q+d=\dim(V).$ Let's denote $\cl(Q)$ as $\cl_{d,p,q}.$ Then we have a direct sum decomposition $\cl_{d,p,q} = \cl_{p,q} \oplus J$ into $\cl_{p,q}$-modules. It was shown in \cite{Structure} that when $d=1$ this decomposition is responsible for a semi-direct product structure of the group of units $\cl_{1,p,q}^\ast$ of $\cl_{1,p,q}$ and of all of its subgroups such as the {\it Clifford group\/} $\Gamma(1,p,k)$ and the {\it special Clifford groups\/} $\Gamma^{\pm}(1,p,k) = \Gamma(1,p,k) \cap \cl_{1,p,q}^{\pm},$ where $\cl_{1,p,q}^{+}$ (resp. $\cl_{1,p,q}^{-})$ denotes the even (resp. odd) part of $\cl_{1,p,q}.$  The Clifford group was defined as $\Gamma(1,p,k)=\{g \in \cl_{1,p,q}^\ast \,|\, g v g^{-1} \in V, v \in V\},$ that is, without a twist.\footnote{See Crumeyrolle \cite{Crum90} for a definition of the Clifford group with the twist given by $\alpha:$ in Crumeyrolle's notation $\alpha$ denotes the principal automorphism or the grade involution in $\cl(Q).$ Then the {\it twisted Clifford group\/} is defined as $\Gamma_{\alpha}(1,p,k)=\{g \in \cl_{1,p,q}^\ast \,|\, \alpha(g) v g^{-1} \in V, v \in V\}.$}
Let $N: \Gamma(1,p,k) \rightarrow \cl_{1,p,k}^{\ast}$ be defined as $N(g)=\revtop{g}g.$ Then we define the {\it reduced Clifford groups\/} as $\Gamma^{\pm}_{0}(1,p,k) = \ker N \cap \Gamma^{\pm}(1,p,k).$  The $\Pin(1,p,q)$ and $\Spin(1,p,k)$ groups are then:
\begin{equation}
\Pin(1,p,q) = \{g \in \Gamma(1,p,q) \,|\, N(g)=\pm 1\},\quad 
\Spin(1,p,q) = \{g \in \Gamma^{+}(1,p,q) \,|\, N(g)=\pm 1\}.
\end{equation}
In preparation for our computations below, from now on we assume that $p+q$ is an odd positive integer. Let $G=\{1 + v \e_1 | v \in V', \e_1^2=0\}$ be a subgroup of $\cl_{1,p,q}^{\ast}.$ Then it was proven in \cite{Structure} that
\begin{equation}
\Pin(1,p,q) = G \rtimes \Gamma^{\pm}_0(p,k), \quad 
\Spin(1,p,q) = G \rtimes \Spin(p,k).
\label{eq:semiprod}
\end{equation}
In the above, symbol $\rtimes$ denotes a semi-direct product with the group on the right acting on the group on the left. For example, it will be of interest to us to note that the homogeneous Galilei group of rigid motions $G_6=\BR^3 \rtimes {\rm  SO}(3)$ in $\BR^3$ is isomorphic to ${\rm SO}^{+}(1,0,3)$ and it is doubly covered by $\Spin^{+}(1,0,3),$ the identity component of $\Spin(1,0,3).$ For a similar result to (\ref{eq:semiprod}) when one considers the twisted Clifford group $\Gamma_{\alpha}(1,p,k)$ and a twisted map 
$N_{\alpha}: \Gamma_{\alpha}(1,p,k) \rightarrow \cl_{1,p,k}^{\ast}$ defined as $N_{\alpha}(g)=\bar{g}g$ where $\bar{\phantom{x}}$ denotes the conjugation in $\cl_{1,p,q},$ see \cite{Brooke7880,Port95}. 

In the following two sections we will use approach and notation from Selig \cite{Selig} where the author denotes the degenerate Clifford algebra $\cl_{d,p,q}$ as $C(p,q,d).$ Furthermore Selig uses twisted groups and defines the $\Pin$ and $\Spin$ groups as follows:
\begin{eqnarray}
\Pin(n) &=&
\{\bg \in C(0,n,0): \bg\bg^{\ast}=\Id \mbox{ and } \alpha(\bg)\bx\bg^{\ast} \in V 
 \mbox{ for all }  \bx \in V\}, \label{eq:pindef} \\
\Spin(n) &=&
\{\bg \in C^{+}(0,n,0): \bg\bg^{\ast}=\Id \mbox{ and } \bg\bx\bg^{\ast} \in V 
 \mbox{ for all }  \bx \in V\}, \label{eq:spindef}
\end{eqnarray} 
with $\phantom{}^{\ast}$ denoting the conjugation $\bar{\phantom{x}}$ in the Clifford algebra $C(p,q,d).$ It is implicit in the definitions above that the actions of $\Pin$ and $\Spin$ on $V$ are $\bx \mapsto \alpha(\bg)\bx\bg^{\ast}$ and $\bx \mapsto \bg\bx\bg^{\ast}$ respectively.

\subsection{Group $\Pin(3)$}
\label{sec:pin3}

In this section we will perform some computations with $\pin3.$ In particular, we will find all possible forms of the elements in $\pin3$ and verify some facts about that group. We begin by assigning a diagonal matrix to the bilinear form $B.$ Grassmann basis for $\cl_{0,3}$ will be stored in the variable {\tt clibas}. Following Selig we re-name Clifford conjugation as a procedure {\tt star} and define a Euclidean norm on $V=\BR^3$ as a procedure {\tt Enorm}. We will also define some additional Maple procedures that will be useful below.
\begin{maplegroup}
\begin{mapleinput}
\mapleinline{active}{1d}{B:=diag(-1$3);eval(makealiases(3)):clibas:=cbasis(3);}{%
}
\end{mapleinput}
\mapleresult
\begin{maplelatex}
\mapleinline{inert}{2d}{B := matrix([[-1, 0, 0], [0, -1, 0], [0, 0, -1]]);}{%
\[
B :=  \left[ 
{\begin{array}{rrr}
-1 & 0 & 0 \\
0 & -1 & 0 \\
0 & 0 & -1
\end{array}}
 \right] 
\]
}
\end{maplelatex}

\begin{maplelatex}
\mapleinline{inert}{2d}{clibas := [Id, e1, e2, e3, e12, e13, e23, e123];}{%
\[
\mathit{clibas} := [\mathit{Id}, \,\mathit{e1}, \,\mathit{e2}, \,
\mathit{e3}, \,\mathit{e12}, \,\mathit{e13}, \,\mathit{e23}, \,
\mathit{e123}]
\]
}
\end{maplelatex}

\end{maplegroup}
\begin{maplegroup}
\begin{mapleinput}
\mapleinline{active}{1d}{
star:=proc(x) conjugation(x) end: #star (conjugation) operation in Cl(0,3)
Enorm:=v->simplify(scalarpart(v &c star(v))): #(pseudo)Euclidean norm in V
alpha:=proc(x) gradeinv(x) end: #alpha (grade involution) operation in Cl(0,3)
scalarprod:=(x,y)->scalarpart(1/2*(x &c star(y) + star(y) &c x)):
}{%
}
\end{mapleinput}
\end{maplegroup}
\begin{maplegroup}
\begin{mapleinput}
\mapleinline{active}{1d}{
Pin_action:=(x,g)->clicollect(simplify(alpha(g) &c x &c star(g)));#action of Pin(3)}{%
}
\end{mapleinput}

\mapleresult
\begin{maplelatex}
\mapleinline{inert}{2d}{Pin_action := proc (x, g) options operator, arrow;
clicollect(simplify(`&c`(`&c`(alpha(g),x),star(g)))) end;}{%
\[
\mathit{Pin\_action} := (x, \,g)\rightarrow \mathrm{clicollect}(
\mathrm{simplify}((\alpha (g)\,\mathrm{`\&c`}\,x)\,\mathrm{`\&c`}
\,\mathrm{star}(g)))
\]
}
\end{maplelatex}
\end{maplegroup}

\begin{maplegroup}
\begin{mapleinput}
\mapleinline{active}{1d}{
Spin_action:=(x,g)->clicollect(simplify(g &c x &c star(g)));#action of Spin(3)}{%
}
\end{mapleinput}

\mapleresult
\begin{maplelatex}
\mapleinline{inert}{2d}{Spin_action := proc (x, g) options operator, arrow;
clicollect(simplify(`&c`(`&c`(g,x),star(g)))) end;}{%
\[
\mathit{Spin\_action} := (x, \,g)\rightarrow \mathrm{simplify}((g
\,\mathrm{`\&c`}\,x)\,\mathrm{`\&c`}\,\mathrm{star}(g))
\]
}
\end{maplelatex}
\end{maplegroup}
\noindent
Let $v,v_1,v_2$ be three arbitrary vectors in $\BR^3$ with some undetermined coefficients expressed in a pseudo-orthonormal basis $\{\e_1,\e_2,\e_3\}:$
\begin{maplegroup}
\begin{mapleinput}
\mapleinline{active}{1d}{
v:=c1*e1+c2*e2+c3*e3:v1:=c11*e1+c12*e2+c13*e3:v2:=c21*e1+c22*e2+c23*e3:}{%
}
\end{mapleinput}
\end{maplegroup}
\noindent
Then the Euclidean norm in $\BR^3$ is:
\begin{maplegroup}
\begin{mapleinput}
\mapleinline{active}{1d}{Enorm(v);}{%
}
\end{mapleinput}

\mapleresult
\begin{maplelatex}
\mapleinline{inert}{2d}{c1^2+c2^2+c3^2;}{%
\[
\mathit{c1}^{2} + \mathit{c2}^{2} + \mathit{c3}^{2}
\]
}
\end{maplelatex}
\end{maplegroup}
\noindent
The action of $\Pin$ on $\cl_{0,3}$ is realized as the procedure {\tt Pin\_action} defined above. Let's verify Selig's claim (\cite{Selig}, page 153) that when 
$\bx,\bg$ are both in $V,$ then $\bg \bx \bg^{\ast}$ automatically belongs to~$V:$
\begin{maplegroup}
\begin{mapleinput}
\mapleinline{active}{1d}{Pin_action(v,v1);}{%
}
\end{mapleinput}
\mapleresult
\begin{maplelatex}
\mapleinline{inert}{2d}{(c11^2*c3-c13^2*c3-2*c13*c12*c2+c12^2*c3-2*c13*c11*c1)*e3+(-c11^2*c1-
2*c11*c13*c3+c12^2*c1+c13^2*c1-2*c11*c12*c2)*e1+(c11^2*c2-2*c12*c11*c1
-2*c12*c13*c3+c13^2*c2-c12^2*c2)*e2;}{%
\maplemultiline{
(\mathit{c11}^{2}\,\mathit{c3} - \mathit{c13}^{2}\,\mathit{c3} - 
2\,\mathit{c13}\,\mathit{c12}\,\mathit{c2} + \mathit{c12}^{2}\,
\mathit{c3} - 2\,\mathit{c13}\,\mathit{c11}\,\mathit{c1})\,
\mathit{e3} \\
\mbox{} + ( - \mathit{c11}^{2}\,\mathit{c1} - 2\,\mathit{c11}\,
\mathit{c13}\,\mathit{c3} + \mathit{c12}^{2}\,\mathit{c1} + 
\mathit{c13}^{2}\,\mathit{c1} - 2\,\mathit{c11}\,\mathit{c12}\,
\mathit{c2})\,\mathit{e1} \\
\mbox{} + (\mathit{c11}^{2}\,\mathit{c2} - 2\,\mathit{c12}\,
\mathit{c11}\,\mathit{c1} - 2\,\mathit{c12}\,\mathit{c13}\,
\mathit{c3} + \mathit{c13}^{2}\,\mathit{c2} - \mathit{c12}^{2}\,
\mathit{c2})\,\mathit{e2} }
}
\end{maplelatex}

\end{maplegroup}
\noindent
As we can see from the above, the output of {\tt Pin\_action(v,v1)} belongs to $V.$ In order to check that indeed the action of $\Pin$ in $V$ preserves the scalar product, we will first find all possible forms of $\bg \in \Pin(3).$ Recall that according to (\ref{eq:pindef}) any element $\bg \in \Pin(n)$ must satisfy two conditions: $(1)\; \bg \Star{\bg} = \Id$ and 
$(2)\; \alpha(\bg) \bv \Star{\bg} \in V$ for any $\bv \in V.$
\noindent
Suppose that $\bg$ is an arbitrary element in $\cl_{0,3,0}$ expressed in {\sc CLIFFORD} in terms of the Grassmann basis $\{\Id,\e_1,\e_2,\e_{12}\}$\footnote{Recall that $\be_{12}$ was defined above as an alias of $\be_1 \w \be_2$ with the command {\tt makealiases}.} 
\begin{maplegroup}
\begin{mapleinput}
\mapleinline{active}{1d}{
g:=add(x.i * clibas[i],i=1..nops(clibas)); #a general element in Cl(0,3,0)}{%
}
\end{mapleinput}

\mapleresult
\begin{maplelatex}
\mapleinline{inert}{2d}{g := x1*Id+x2*e1+x3*e2+x4*e3+x5*e12+x6*e13+x7*e23+x8*e123;}{%
\[
g := \mathit{x1}\,\mathit{Id} + \mathit{x2}\,\mathit{e1} + 
\mathit{x3}\,\mathit{e2} + \mathit{x4}\,\mathit{e3} + \mathit{x5}
\,\mathit{e12} + \mathit{x6}\,\mathit{e13} + \mathit{x7}\,
\mathit{e23} + \mathit{x8}\,\mathit{e123}
\]
}
\end{maplelatex}

\end{maplegroup}
\noindent
We will now attempt to find conditions that the coefficients $x_i,i=1,\ldots,8,$ must satisfy so  that $\bg \Star{\bg} = \Id.$ We will again use the command {\tt clisolve2}. In order to shorten its outputs, additional aliases $\kappa_j,j=1,\ldots,5,$ need to be defined 
(see Appendix~2).

\begin{maplegroup}
The first condition (1) gives:

\end{maplegroup}
\begin{maplegroup}
\begin{mapleinput}
\mapleinline{active}{1d}{sol:=clisolve2(cmul(g,star(g))-Id,[x.(1..8)]);}{%
}
\end{mapleinput}
\mapleresult
\begin{maplelatex}
\mapleinline{inert}{2d}{sol := [\{x1 = -(-x7*kappa5-x5*x4+x6*x3)/x8, x8 = x8, x4 = x4, x7 = x7, x6 = x6, x5 = x5, x2 = kappa5, x3 = x3\}, 
\{x8 = 0, x1 = kappa4/x7, x4 = x4, x7 = x7, x6 = x6, x5 = x5, x3 = x3, x2 = (-x5*x4+x6*x3)/x7\}, 
\{x7 = 0, x8 = 0, x1 = kappa3/x6, x3 = x4*x5/x6, x2 = x2, x4 = x4, x6 = x6, x5 = x5\}, 
\{x4 = 0, x7 = 0, x6 = 0, x8 = 0, x2 = x2, x5 = x5, x3 = x3, x1 = kappa2\}, 
\{x7 = 0, x6 = 0, x5 = 0, x8 = 0, x2 = x2, x4 = x4, x3 = x3, x1 = kappa1\}];}{%
\maplemultiline{
\mathit{sol} := [\{\mathit{x1}= - {\displaystyle \frac { - 
\mathit{x7}\,\kappa 5 - \mathit{x5}\,\mathit{x4} + \mathit{x6}\,
\mathit{x3}}{\mathit{x8}}} , \,\mathit{x8}=\mathit{x8}, \,
\mathit{x4}=\mathit{x4}, \,\mathit{x7}=\mathit{x7}, \,\mathit{x6}
=\mathit{x6}, \,\mathit{x5}=\mathit{x5}, \,\mathit{x2}=\kappa 5, 
 \\
\mathit{x3}=\mathit{x3}\}, 
\{ \mathit{x8}=0, \,\mathit{x1}={\displaystyle \frac {\kappa 4}{
\mathit{x7}}} , \,\mathit{x4}=\mathit{x4}, \,\mathit{x7}=\mathit{
x7}, \,\mathit{x6}=\mathit{x6}, \,\mathit{x5}=\mathit{x5}, \,
\mathit{x3}=\mathit{x3}, \,\mathit{x2}={\displaystyle \frac { - 
\mathit{x5}\,\mathit{x4} + \mathit{x6}\,\mathit{x3}}{\mathit{x7}}
} \},  \\
\{\mathit{x7}=0, \,\mathit{x8}=0, \,\mathit{x1}={\displaystyle 
\frac {\kappa 3}{\mathit{x6}}} , \,\mathit{x3}={\displaystyle 
\frac {\mathit{x4}\,\mathit{x5}}{\mathit{x6}}} , \,\mathit{x2}=
\mathit{x2}, \,\mathit{x4}=\mathit{x4}, \,\mathit{x6}=\mathit{x6}
, \,\mathit{x5}=\mathit{x5}\},  \\
\{\mathit{x4}=0, \,\mathit{x7}=0, \,\mathit{x6}=0, \,\mathit{x8}=
0, \,\mathit{x2}=\mathit{x2}, \,\mathit{x5}=\mathit{x5}, \,
\mathit{x3}=\mathit{x3}, \,\mathit{x1}=\kappa 2\},  \\
\{\mathit{x7}=0, \,\mathit{x6}=0, \,\mathit{x5}=0, \,\mathit{x8}=
0, \,\mathit{x2}=\mathit{x2}, \,\mathit{x4}=\mathit{x4}, \,
\mathit{x3}=\mathit{x3}, \,\mathit{x1}=\kappa 1\}] }
}
\end{maplelatex}
\end{maplegroup}
\noindent
Thus, there are five different possible solutions, three of which requiring respectively that $x_6,x_7$ and $x_8$ be non-zero. Let's substitute these solutions into $\bg.$
\begin{maplegroup}
\begin{mapleinput}
\mapleinline{active}{1d}{for i from 1 to nops(sol) do g.i:=subs(sol[i],g) od;}{%
}
\end{mapleinput}

\mapleresult
\begin{maplelatex}
\mapleinline{inert}{2d}{g1 :=
-(-x7*kappa5-x5*x4+x6*x3)*Id/x8+kappa5*e1+x3*e2+x4*e3+x5*e12+x6*e13+x7
*e23+x8*e123;}{%
\maplemultiline{
\mathit{g1} :=  - {\displaystyle \frac {( - \mathit{x7}\,\kappa 5
 - \mathit{x5}\,\mathit{x4} + \mathit{x6}\,\mathit{x3})\,\mathit{
Id}}{\mathit{x8}}}  + \kappa 5\,\mathit{e1} + \mathit{x3}\,
\mathit{e2} + \mathit{x4}\,\mathit{e3} + \mathit{x5}\,\mathit{e12
} + \mathit{x6}\,\mathit{e13} + \mathit{x7}\,\mathit{e23} \\
\mbox{} + \mathit{x8}\,\mathit{e123} }
}
\end{maplelatex}

\begin{maplelatex}
\mapleinline{inert}{2d}{g2 :=
kappa4*Id/x7+(-x5*x4+x6*x3)*e1/x7+x3*e2+x4*e3+x5*e12+x6*e13+x7*e23;}{%
\[
\mathit{g2} := {\displaystyle \frac {\kappa 4\,\mathit{Id}}{
\mathit{x7}}}  + {\displaystyle \frac {( - \mathit{x5}\,\mathit{
x4} + \mathit{x6}\,\mathit{x3})\,\mathit{e1}}{\mathit{x7}}}  + 
\mathit{x3}\,\mathit{e2} + \mathit{x4}\,\mathit{e3} + \mathit{x5}
\,\mathit{e12} + \mathit{x6}\,\mathit{e13} + \mathit{x7}\,
\mathit{e23}
\]
}
\end{maplelatex}

\begin{maplelatex}
\mapleinline{inert}{2d}{g3 := kappa3*Id/x6+x2*e1+x4*x5*e2/x6+x4*e3+x5*e12+x6*e13;}{%
\[
\mathit{g3} := {\displaystyle \frac {\kappa 3\,\mathit{Id}}{
\mathit{x6}}}  + \mathit{x2}\,\mathit{e1} + {\displaystyle 
\frac {\mathit{x4}\,\mathit{x5}\,\mathit{e2}}{\mathit{x6}}}  + 
\mathit{x4}\,\mathit{e3} + \mathit{x5}\,\mathit{e12} + \mathit{x6
}\,\mathit{e13}
\]
}
\end{maplelatex}

\begin{maplelatex}
\mapleinline{inert}{2d}{g4 := kappa2*Id+x2*e1+x3*e2+x5*e12;}{%
\[
\mathit{g4} := \kappa 2\,\mathit{Id} + \mathit{x2}\,\mathit{e1}
 + \mathit{x3}\,\mathit{e2} + \mathit{x5}\,\mathit{e12}
\]
}
\end{maplelatex}

\begin{maplelatex}
\mapleinline{inert}{2d}{g5 := kappa1*Id+x2*e1+x3*e2+x4*e3;}{%
\[
\mathit{g5} := \kappa 1\,\mathit{Id} + \mathit{x2}\,\mathit{e1}
 + \mathit{x3}\,\mathit{e2} + \mathit{x4}\,\mathit{e3}
\]
}
\end{maplelatex}

\end{maplegroup}
\noindent
The above are five different types of $\bg$ in $\cl_{0,3,0}$ satisfying $\bg \Star{\bg} = \Id.$
\begin{maplegroup}
\begin{mapleinput}
\mapleinline{active}{1d}{for i from 1 to nops(sol) do simplify(cmul(g.i,star(g.i))) od;}{%
}
\end{mapleinput}

\mapleresult
\begin{maplelatex}
\mapleinline{inert}{2d}{Id;Id;Id;Id;Id;}{%
\[
\mathit{Id},\,\mathit{Id},\,\mathit{Id},\,\mathit{Id},\,\mathit{Id}
\]
}
\end{maplelatex}
\end{maplegroup}
\noindent
We need to make sure now that each $g_i$ displayed above satisfies also the second condition $(2),$ namely, $\alpha(\bg) \bv \Star{\bg}$ is in $V$ for any $\bv \in V.$ We begin with the simplest element~$g_5.$ By computing the $\Pin$ group action on $\bv$ and requiring that the result be a $1$-vector, we get for $g_5:$
\begin{maplegroup}
\begin{mapleinput}
\mapleinline{active}{1d}{Pin_action(v,g5);}{%
}
\end{mapleinput}

\mapleresult
\begin{maplelatex}
\mapleinline{inert}{2d}{(2*kappa1*c3*x4+2*kappa1*c1*x2+2*kappa1*c2*x3)*Id+(-2*x3*x4*c2-2*x4*x
2*c1-2*x4^2*c3+c3)*e3+(-2*x4*x2*c3-2*x3*x2*c2-2*x2^2*c1+c1)*e1+(-2*x3^
2*c2-2*x4*x3*c3-2*x2*x3*c1+c2)*e2;}{%
\maplemultiline{
(2\,\kappa 1\,\mathit{c3}\,\mathit{x4} + 2\,\kappa 1\,\mathit{c1}
\,\mathit{x2} + 2\,\kappa 1\,\mathit{c2}\,\mathit{x3})\,\mathit{
Id} + ( - 2\,\mathit{x3}\,\mathit{x4}\,\mathit{c2} - 2\,\mathit{
x4}\,\mathit{x2}\,\mathit{c1} - 2\,\mathit{x4}^{2}\,\mathit{c3}
 + \mathit{c3})\,\mathit{e3} \\
\mbox{} + ( - 2\,\mathit{x4}\,\mathit{x2}\,\mathit{c3} - 2\,
\mathit{x3}\,\mathit{x2}\,\mathit{c2} - 2\,\mathit{x2}^{2}\,
\mathit{c1} + \mathit{c1})\,\mathit{e1} \\
\mbox{} + ( - 2\,\mathit{x3}^{2}\,\mathit{c2} - 2\,\mathit{x4}\,
\mathit{x3}\,\mathit{c3} - 2\,\mathit{x2}\,\mathit{x3}\,\mathit{
c1} + \mathit{c2})\,\mathit{e2} }
}
\end{maplelatex}

\end{maplegroup}
\noindent
It should be clear from the above that since the coefficient of $Id$ must be zero for any $c_1,c_2,c_3,$ either $x_2=x_3=x_4=0$ or $\kappa_1=0.$ Let $\varepsilon = \pm 1.$\footnote{In Maple one way to make $\varepsilon = \pm 1$ is to define  {\tt alias(eps=RootOf(\_Z\^{}2-1)):}. In the following, Maple outputs will contain the alias $eps.$} Thus, the former gives 
$\bg = \pm \Id,$
\begin{maplegroup}
\begin{mapleinput}
\mapleinline{active}{1d}{g.5.1:=subs(\{x2=0,x3=0,x4=0\},g5);}{%
}
\end{mapleinput}

\mapleresult
\begin{maplelatex}
\mapleinline{inert}{2d}{g51 := eps*Id;}{%
\[
\mathit{g51} := \mathit{eps}\,\mathit{Id}
\]
}
\end{maplelatex}

\end{maplegroup}
\noindent
while the latter gives
\begin{maplegroup}
\begin{mapleinput}
\mapleinline{active}{1d}{g.5.2:=subs(\{kappa1=0,x4=lambda1\},g5);}{%
}
\end{mapleinput}

\mapleresult
\begin{maplelatex}
\mapleinline{inert}{2d}{g52 := lambda1*e3+x2*e1+x3*e2;}{%
\[
\mathit{g52} := \lambda 1\,\mathit{e3} + \mathit{x2}\,\mathit{e1}
 + \mathit{x3}\,\mathit{e2}
\]
}
\end{maplelatex}

\end{maplegroup}
\noindent
where $\lambda_1=\pm \sqrt{1-x_2^2-x_3^2}.$\footnote{In Maple we define {\tt alias(lambda1=RootOf(\_Z\^{}2+x2\^{}2+x3\^{}2-1)):.}} We will collect all $\Pin$ group elements in a set {\tt Pin\_group}.
\begin{maplegroup}
\begin{mapleinput}
\mapleinline{active}{1d}{Pin_group:=\{g.5.1,g.5.2\};}{%
}
\end{mapleinput}

\mapleresult
\begin{maplelatex}
\mapleinline{inert}{2d}{Pin_group := \{\mathit{eps}*Id, lambda1*e3+x2*e1+x3*e2\};}{%
\[
\mathit{Pin\_group} := \{\mathit{eps}\,\mathit{Id},\,\lambda 1\,\mathit{e3} + \mathit{x2}\,
\mathit{e1} + \mathit{x3}\,\mathit{e2}\}
\]
}
\end{maplelatex}

\end{maplegroup}

Similarly, we consider $g_4.$ We assign the identity coefficient of the $\Pin$ action $\alpha(\bg) \bv \Star{\bg}$ to a variable $eq$ and find a solution to the resulting two equations that will be parameterized by~$c_1,c_2:$
\begin{maplegroup}
\begin{mapleinput}
\mapleinline{active}{1d}{a:=Pin_action(v,g4);}{%
}
\end{mapleinput}

\mapleresult
\begin{maplelatex}
\mapleinline{inert}{2d}{a :=
(2*x3*kappa2*c2-2*x5*x2*c2+2*x5*x3*c1+2*x2*kappa2*c1)*Id+c3*e3+(-2*x2^
2*c1-2*x5*kappa2*c2-2*x3*x2*c2+c1-2*x5^2*c1)*e1+(-2*x5^2*c2-2*x2*x3*c1
+2*kappa2*x5*c1-2*x3^2*c2+c2)*e2;}{%
\maplemultiline{
a := (2\,\mathit{x3}\,\kappa 2\,\mathit{c2} - 2\,\mathit{x5}\,
\mathit{x2}\,\mathit{c2} + 2\,\mathit{x5}\,\mathit{x3}\,\mathit{
c1} + 2\,\mathit{x2}\,\kappa 2\,\mathit{c1})\,\mathit{Id} + 
\mathit{c3}\,\mathit{e3} \\
\mbox{} + ( - 2\,\mathit{x2}^{2}\,\mathit{c1} - 2\,\mathit{x5}\,
\kappa 2\,\mathit{c2} - 2\,\mathit{x3}\,\mathit{x2}\,\mathit{c2}
 + \mathit{c1} - 2\,\mathit{x5}^{2}\,\mathit{c1})\,\mathit{e1}
 \\
\mbox{} + ( - 2\,\mathit{x5}^{2}\,\mathit{c2} - 2\,\mathit{x2}\,
\mathit{x3}\,\mathit{c1} + 2\,\kappa 2\,\mathit{x5}\,\mathit{c1}
 - 2\,\mathit{x3}^{2}\,\mathit{c2} + \mathit{c2})\,\mathit{e2} }
}
\end{maplelatex}

\end{maplegroup}
\begin{maplegroup}
\begin{mapleinput}
\mapleinline{active}{1d}{eq:=collect(coeff(a,Id),\{c1,c2\}):eq1:=coeff(eq,c1):eq2:=coeff(eq,c2):
sol:=[solve(\{eq1,eq2\},\{x2,x5,x3\})];}{%
}
\end{mapleinput}

\mapleresult
\begin{maplelatex}
\mapleinline{inert}{2d}{sol := [\{x5 = 0, x2 = lambda2, x3 = x3\}, \{x3 = 0, x2 = 0, x5 =
x5\}];}{%
\[
\mathit{sol} := [\{\mathit{x5}=0, \,\mathit{x2}=\lambda 2, \,
\mathit{x3}=\mathit{x3}\}, \,\{\mathit{x3}=0, \,\mathit{x2}=0, \,
\mathit{x5}=\mathit{x5}\}]
\]
}
\end{maplelatex}

\end{maplegroup}
\noindent
In the above, $\lambda_2=\pm \sqrt{1-x_3^2}.$ Likewise, we set 
$\lambda_3=\pm \sqrt{1-x_5^2}.$\footnote{Both are defined as aliases:
{\tt alias(lambda2=RootOf(\_Z\^{}2-1+x3\^{}2)):}   
{\tt alias(lambda3=RootOf(\_Z\^{}2-1+x5\^{}2)):}.} The two new elements we assign to $g_{41},g_{42}$ and add to {\tt Pin\_group}. 
\begin{maplegroup}
\begin{mapleinput}
\mapleinline{active}{1d}{for i from 1 to nops(sol) do g.4.i:=simplify(subs(sol[i],g4)) od:
Pin_group:=Pin_group union \{g41,g42\};}{%
}
\end{mapleinput}

\mapleresult
\begin{maplelatex}
\mapleinline{inert}{2d}{
Pin_group := \{ eps*Id, lambda2*e1+x3*e2, lambda3*Id+x5*e12, lambda1*e3+x2*e1+x3*e2\};}{%
\[
\mathit{Pin\_group} := \{\mathit{eps}\,\mathit{Id},\,
\lambda 1\,\mathit{e3} + \mathit{x2}\,\mathit{e1}+ \mathit{x3}\,\mathit{e2},\,
\lambda 2\,\mathit{e1} + \mathit{x3}\,\mathit{e2}, \,
\lambda 3\,\mathit{Id} + \mathit{x5}\,\mathit{e12}
\}
\]
}
\end{maplelatex}

\end{maplegroup}
In order to continue with $g_3$ displayed above, we must make the assumption $x_6 \neq 0$ known to Maple. Then, the action of $g_3$ on a vector can be computed.\footnote{We won't display it due to its length.}  
\begin{maplegroup}
\begin{mapleinput}
\mapleinline{active}{1d}{assume(x6>0,x6<0);a:=Pin_action(v,g3):}{%
}
\end{mapleinput}

\end{maplegroup}
\noindent
As in the previous two case, the quantity $a$ is spanned by $\{Id,e_1,e_2,e_3\}.$ We will isolate the coefficient of the identity element in $a,$ assign it to a variable $eq,$ and then determine for which values of $x_2,x_4,x_5,x_6$ it will be automatically zero for every choice of $c_1,c_2,c_3.$ This will require solving a set of three equations $\{eq_1,eq_2,eq_3\}$ for $x_2,x_4,x_5,x_6.$  Maple reminds us that $x_6 \neq 0$ by displaying it as $\mathit{x6\symbol{126}}:$
\begin{maplegroup}
\begin{mapleinput}
\mapleinline{active}{1d}{cliterms(a);}{%
}
\end{mapleinput}

\mapleresult
\begin{maplelatex}
\mapleinline{inert}{2d}{\{Id, e3, e1, e2\};}{%
\[
\{\mathit{Id}, \,\mathit{e3}, \,\mathit{e1}, \,\mathit{e2}\}
\]
}
\end{maplelatex}
\end{maplegroup}
\begin{maplegroup}
\begin{mapleinput}
\mapleinline{active}{1d}{eq:=collect(coeff(a,Id),\{c.(1..3)\});}{%
}
\end{mapleinput}

\mapleresult
\begin{maplelatex}
\mapleinline{inert}{2d}{eq :=
-2*(x6^3*x2-lambda4*x4*x6)*c3/(x6^2)-2*(-lambda4*x2*x6-x5^2*x4*x6-x4*x
6^3)*c1/(x6^2)-2*(-lambda4*x5*x4+x2*x5*x6^2)*c2/(x6^2);}{%
\maplemultiline{
\mathit{eq} :=  - 2\,{\displaystyle \frac {(\mathit{x6
\symbol{126}}^{3}\,\mathit{x2} - \lambda 4\,\mathit{x4}\,\mathit{
x6\symbol{126}})\,\mathit{c3}}{\mathit{x6\symbol{126}}^{2}}}  - 2
\,{\displaystyle \frac {( - \lambda 4\,\mathit{x2}\,\mathit{x6
\symbol{126}} - \mathit{x5}^{2}\,\mathit{x4}\,\mathit{x6
\symbol{126}} - \mathit{x4}\,\mathit{x6\symbol{126}}^{3})\,
\mathit{c1}}{\mathit{x6\symbol{126}}^{2}}}  \\
\mbox{} - 2\,{\displaystyle \frac {( - \lambda 4\,\mathit{x5}\,
\mathit{x4} + \mathit{x2}\,\mathit{x5}\,\mathit{x6\symbol{126}}^{
2})\,\mathit{c2}}{\mathit{x6\symbol{126}}^{2}}}  }
}
\end{maplelatex}

\end{maplegroup}
\begin{maplegroup}
\begin{mapleinput}
\mapleinline{active}{1d}{for i from 1 to 3 do eq.i:=coeff(eq,c.i) od:
sol:=solve(\{eq.(1..3)\},\{x6,x4,x5,x2\});}{%
}
\end{mapleinput}

\mapleresult
\begin{maplelatex}
\mapleinline{inert}{2d}{sol := \{x2 = 0, x6 = x6, x5 = x5, x4 = 0\};}{%
\[
\mathit{sol} := \{\mathit{x2}=0, \,\mathit{x6\symbol{126}}=
\mathit{x6\symbol{126}}, \,\mathit{x5}=\mathit{x5}, \,\mathit{x4}
=0\}
\]
}
\end{maplelatex}

\end{maplegroup}
\noindent
This time we only have one solution which we then substitute into $g_3.$ 
\begin{maplegroup}
\begin{mapleinput}
\mapleinline{active}{1d}{g31:=subs(sol,g3);}{%
}
\end{mapleinput}

\mapleresult
\begin{maplelatex}
\mapleinline{inert}{2d}{g31 := lambda5*Id/x6+x5*e12+x6*e13;}{%
\[
\mathit{g31} := {\displaystyle \frac {\lambda 5\,\mathit{Id}}{
\mathit{x6\symbol{126}}}}  + \mathit{x5}\,\mathit{e12} + \mathit{x6\symbol{126}}\,\mathit{e13}
\]
}
\end{maplelatex}
\end{maplegroup}
\noindent
where $\lambda_5$ is another alias (see Appendix 2).
\begin{maplegroup}
\begin{mapleinput}
\mapleinline{active}{1d}{Pin_group:=Pin_group union \{g31\};}{%
}
\end{mapleinput}

\mapleresult
\begin{maplelatex}
\mapleinline{inert}{2d}{Pin_group := \{lambda3*Id+x5*e12, lambda5*Id/x6+x5*e12+x6*e13,
lambda1*e3+x2*e1+x3*e2, lambda2*e1+x3*e2, eps*Id\};}{%
\maplemultiline{
\mathit{Pin\_group} := \{\lambda 3\,\mathit{Id} + \mathit{x5}\,
\mathit{e12}, \,{\displaystyle \frac {\lambda 5\,\mathit{Id}}{
\mathit{x6\symbol{126}}}}  + \mathit{x5}\,\mathit{e12} + \mathit{
x6\symbol{126}}\,\mathit{e13}, \,\lambda 1\,\mathit{e3} + 
\mathit{x2}\,\mathit{e1} + \mathit{x3}\,\mathit{e2},  \\
\lambda 2\,\mathit{e1} + \mathit{x3}\,\mathit{e2}, \,\mathit{eps}
\,\mathit{Id}\} }
}
\end{maplelatex}

\end{maplegroup}

By continuing in the similar fashion with the elements $g_2$ and $g_1,$ one can find all general types of $\Pin(3).$ Finally, all general elements of $\Pin(3)$ can be displayed:
\begin{maplegroup}
\begin{mapleinput}
\mapleinline{active}{1d}{'Pin_group'=Pin_group;}{%
}
\end{mapleinput}

\mapleresult
\begin{maplelatex}
\mapleinline{inert}{2d}{Pin_group = \{lambda3*Id+x5*e12, lambda5*Id/x6+x5*e12+x6*e13,
lambda7*Id/x7+x5*e12+x6*e13+x7*e23, lambda9*e1+x3*e2+x4*e3+x8*e123,
lambda1*e3+x2*e1+x3*e2, lambda2*e1+x3*e2, eps*Id\};}{%
\maplemultiline{
\mathit{Pin\_group} := \{\lambda 3\,\mathit{Id} + \mathit{x5}\,
\mathit{e12}, \,{\displaystyle \frac {\lambda 5\,\mathit{Id}}{
\mathit{x6\symbol{126}}}}  + \mathit{x5}\,\mathit{e12} + \mathit{
x6\symbol{126}}\,\mathit{e13},  \\
{\displaystyle \frac {\lambda 7\,\mathit{Id}}{\mathit{x7
\symbol{126}}}}  + \mathit{x5}\,\mathit{e12} + \mathit{x6}\,
\mathit{e13} + \mathit{x7\symbol{126}}\,\mathit{e23}, \,\lambda 9
\,\mathit{e1} + \mathit{x3}\,\mathit{e2} + \mathit{x4}\,\mathit{
e3} + \mathit{x8\symbol{126}}\,\mathit{e123},  \\
\lambda 1\,\mathit{e3} + \mathit{x2}\,\mathit{e1} + \mathit{x3}\,
\mathit{e2}, \,\lambda 2\,\mathit{e1} + \mathit{x3}\,\mathit{e2}
, \,\mathit{eps}\,\mathit{Id}\} }
}
\end{maplelatex}
\end{maplegroup}
\noindent
where $\lambda_7$ and $\lambda_9$ are displayed in the Appendix 2. It is a simple matter now to verify that all elements of $\Pin$ displayed in $Pin\_group$ satisfy both conditions $(1)$ and $(2)$ from the definition~(\ref{eq:pindef}).
\begin{maplegroup}
\begin{mapleinput}
\mapleinline{active}{1d}{
for g in Pin_group do evalb(simplify(cmul(g,star(g))=Id)) od;#Condition (1)}{%
}
\end{mapleinput}

\mapleresult
\begin{maplelatex}
\mapleinline{inert}{2d}{true;true;true;true;true;true;}{%
\[
\mathit{true},\,\mathit{true},\,\mathit{true},\,\mathit{true},\,\mathit{true},\,\mathit{true}
\]
}
\end{maplelatex}
\end{maplegroup}
\begin{maplegroup}
\begin{mapleinput}
\mapleinline{active}{1d}{for g in Pin_group do 
evalb(Pin_action(v,g)=vectorpart(Pin_action(v,g),1)) od;#Condition (2)}{%
}
\end{mapleinput}

\mapleresult
\begin{maplelatex}
\mapleinline{inert}{2d}{true;true;true;true;true;true}{%
\[
\mathit{true},\,\mathit{true},\,\mathit{true},\,\mathit{true},\,\mathit{true},\,\mathit{true}
\]
}
\end{maplelatex}
\end{maplegroup}

We are now in position to verify Selig's claim \cite{Selig}, page 153, that the scalar product on $V=\BR^3$ defined in $\cl_{d,p,q}$ as
$$
\bv_1 \cdot \bv_2 = \frac12(\bv_1 \Star{\bv_2} + \bv_2 \Star{\bv_1}),\quad \mbox{ for any } \bv_1,\bv_2 \in \BR^3,
$$
is preserved under the action of the $\Pin$ group. Let $v_1,v_2$ be two arbitrary $1$-vectors defined earlier. Procedure {\tt scalarprod} that gives the scalar product may be defined as follows:
\begin{maplegroup}
\begin{mapleinput}
\mapleinline{active}{1d}{
scalarprod:=(x,y)->scalarpart(1/2*(x &c star(y) + star(y) &c x)):
for g in Pin_group do
\phtab{x} simplify(scalarprod(Pin_action(v1,g),Pin_action(v2,g))-scalarprod(v1,v2)) 
od;}{%
}
\end{mapleinput}

\mapleresult
\begin{maplelatex}
\mapleinline{inert}{2d}{0;0;0;0;0;0;}{%
\[
0,\,0,\,0,\,0,\,0,\,0
\]
}
\end{maplelatex}
\end{maplegroup}
\noindent
Thus, $\Pin(3)$ preserves the scalar product in $\BR^3$ and, therefore, we have a homomorphism from $\Pin(3)$ to ${\rm O}(3)$ which is known to be a double-covering map. In the process, we have found all types of elements in $\Pin(3).$

\subsection{Group $\Spin(3)$}
\label{sec:spin3}

In this section we will perform a few computations with $\Spin(3).$ Once we have found general elements in $\Pin(3),$ it is much easier to find elements in $\Spin(3).$ Recall from (\ref{eq:spindef}) that 
$$
\Spin(3)=\{\bg \in \cl^{+}_{0,3,0}: \bg\Star{\bg}=\Id \mbox{ and } \bg\bx\Star{\bg} \in \BR^3  \mbox{ for all }  \bx \in \BR^3\}.
$$ 
Let's find $gSpin,$ a general element in $\Spin(3).$ Since $\Spin(3) \subset \cl^{+}_{0,3,0},$ we will begin with decomposing $gSpin$ over the even basis elements.
\begin{maplegroup}
\begin{mapleinput}
\mapleinline{active}{1d}{clibaseven:=cbasis(3,'even');}{%
}
\end{mapleinput}

\mapleresult
\begin{maplelatex}
\mapleinline{inert}{2d}{clibaseven := [Id, e12, e13, e23];}{%
\[
\mathit{clibaseven} := [\mathit{Id}, \,\mathit{e12}, \,\mathit{
e13}, \,\mathit{e23}]
\]
}
\end{maplelatex}

\end{maplegroup}
\begin{maplegroup}
\begin{mapleinput}
\mapleinline{active}{1d}{gSpin:=c0*Id+c3*e12+c2*e13+c1*e23;}{%
}
\end{mapleinput}

\mapleresult
\begin{maplelatex}
\mapleinline{inert}{2d}{gSpin := c0*Id+c3*e12+c2*e13+c1*e23;}{%
\[
\mathit{gSpin} := \mathit{c0}\,\mathit{Id} + \mathit{c3}\,
\mathit{e12} + \mathit{c2}\,\mathit{e13} + \mathit{c1}\,\mathit{
e23}
\]
}
\end{maplelatex}

\end{maplegroup}
\noindent
Notice that under the $\Spin$ group action defined as a procedure {\tt Spin\_action}
\begin{maplegroup}
\begin{mapleinput}
\mapleinline{active}{1d}{Spin_action:=(x,g)-> simplify(g &c x &c star(g));}{%
}
\end{mapleinput}

\mapleresult
\begin{maplelatex}
\mapleinline{inert}{2d}{Spin_action := proc (x, g) options operator, arrow;
simplify(`&c`(`&c`(g,x),star(g))) end;}{%
\[
\mathit{Spin\_action} := (x, \,g)\rightarrow \mathrm{simplify}((g
\,\mathrm{`\&c`}\,x)\,\mathrm{`\&c`}\,\mathrm{star}(g))
\]
}
\end{maplelatex}

\end{maplegroup}
\noindent
vectors are automatically mapped into vectors:
\begin{maplegroup}
\begin{mapleinput}
\mapleinline{active}{1d}{Spin_action(v,gSpin)-vectorpart(Spin_action(v,gSpin),1);}{%
}
\end{mapleinput}

\mapleresult
\begin{maplelatex}
\mapleinline{inert}{2d}{0;}{%
\[
0
\]
}
\end{maplelatex}

\end{maplegroup}
\noindent
We just need to make sure that $\bg \Star{\bg} = \Id$ for each $\bg \in \Spin(3).$ To simplify Maple output, we define $\kappa=\sqrt{1-c_1^2-c_2^2-c_3^2}$ as a Maple alias.  
\begin{maplegroup}
\begin{mapleinput}
\mapleinline{active}{1d}{alias(kappa=sqrt(-c1^2-c2^2-c3^2+1)):
sol:=clisolve2(cmul(gSpin,star(gSpin))-Id,[c.(0..3)]);}{%
}
\end{mapleinput}

\mapleresult
\begin{maplelatex}
\mapleinline{inert}{2d}{sol := [\{c1 = c1, c2 = c2, c3 = c3, c0 = kappa\}, \{c1 = c1, c2 =
c2, c3 = c3, c0 = -kappa\}];}{%
\[
\mathit{sol} := [\{\mathit{c1}=\mathit{c1}, \,\mathit{c2}=
\mathit{c2}, \,\mathit{c3}=\mathit{c3}, \,\mathit{c0}=\kappa \}, 
\,\{\mathit{c1}=\mathit{c1}, \,\mathit{c2}=\mathit{c2}, \,
\mathit{c3}=\mathit{c3}, \,\mathit{c0}= - \kappa \}]
\]
}
\end{maplelatex}

\end{maplegroup}
\begin{maplegroup}
\begin{mapleinput}
\mapleinline{active}{1d}{gSpin:=eps*kappa*Id+c3*e12+c2*e13+c1*e23;}{%
}
\end{mapleinput}

\mapleresult
\begin{maplelatex}
\mapleinline{inert}{2d}{gSpin := eps*kappa*Id+c3*e12+c2*e13+c1*e23;}{%
\[
\mathit{gSpin} := \mathit{eps}\,\kappa \,\mathit{Id} + \mathit{c3
}\,\mathit{e12} + \mathit{c2}\,\mathit{e13} + \mathit{c1}\,
\mathit{e23}
\]
}
\end{maplelatex}

\end{maplegroup}
\noindent
Thus, the most general element in $\Spin(3)$ is just 
$\bg = \varepsilon \kappa \Id + c_3 \e_{12} + c_2 \e_{13} + c_1 \e_{23}$ 
where $\varepsilon = \pm 1.$
\noindent
Notice, that the defining properties of $\bg$ are easily checked:
\begin{maplegroup}
\begin{mapleinput}
\mapleinline{active}{1d}{simplify(cmul(gSpin,star(gSpin)));}{%
}
\end{mapleinput}

\mapleresult
\begin{maplelatex}
\mapleinline{inert}{2d}{Id;}{%
\[
\mathit{Id} 
\]
}
\end{maplelatex}

\end{maplegroup}
\begin{maplegroup}
\begin{mapleinput}
\mapleinline{active}{1d}{evalb(Spin_action(v,gSpin)=vectorpart(Spin_action(v,gSpin),1));}{%
}
\end{mapleinput}

\mapleresult
\begin{maplelatex}
\mapleinline{inert}{2d}{true;}{%
\[
\mathit{true}
\]
}
\end{maplelatex}

\end{maplegroup}
\noindent
In fact, element $\bg \in \Spin(3)$ could be identified with a unit quaternion spanned over the basis $\{\Id, \e_{12},\e_{13},\e_{23}\}.$ Then, the $\Star{}$ conjugation becomes the quaternionic conjugation. It can be easily checked by hand or with {\sc CLIFFORD} that the basis (bi)vectors anticommute and square to $-\Id.$
\begin{maplegroup}
\begin{mapleinput}
\mapleinline{active}{1d}{quatbasis:=[e12,e13,e23];}{%
}
\end{mapleinput}

\mapleresult
\begin{maplelatex}
\mapleinline{inert}{2d}{quatbasis := [e12, e13, e23];}{%
\[
\mathit{quatbasis} := [\mathit{e12}, \,\mathit{e13}, \,\mathit{
e23}]
\]
}
\end{maplelatex}
\end{maplegroup}
\begin{maplegroup}
\begin{mapleinput}
\mapleinline{active}{1d}{M:=matrix(3,3,(i,j)->cmul(quatbasis[i],quatbasis[j]));}{%
}
\end{mapleinput}

\mapleresult
\begin{maplelatex}
\mapleinline{inert}{2d}{M := matrix([[-Id, e23, -e13], [-e23, -Id, e12], [e13, -e12,
-Id]]);}{%
\[
M :=  \left[ 
{\begin{array}{ccc}
 - \mathit{Id} & \mathit{e23} &  - \mathit{e13} \\
 - \mathit{e23} &  - \mathit{Id} & \mathit{e12} \\
\mathit{e13} &  - \mathit{e12} &  - \mathit{Id}
\end{array}}
 \right] 
\]
}
\end{maplelatex}

\end{maplegroup}
\noindent
We have unit quaternions on a unit sphere in $\BR^4$ isomorphic to $\Spin(3)$ while the even part of $\cl_{0,3}$ is isomorphic with the quaternionic division ring $\BH.$ $\Spin(3)$ acts on $\BR^3$ through the rotations. In Appendix 3 one can find a procedure {\tt rot} which takes as its first argument a vector $\bv$ in $\BR^3$ while as its second argument it takes a quaternion. For example, a counter-clockwise rotation in the $\{\be_1,\be_2\}$ is accomplished with a help of a unit quaternion $\cos(\frac{\theta}{2})+\sin(\frac{\theta}{2}) \, \be_{12}:$
\begin{maplegroup}
\begin{mapleinput}
\mapleinline{active}{1d}{rot(e1,cos(theta/2)+sin(theta/2)*e12);
rot(e2,cos(theta/2)+sin(theta/2)*e12);
rot(e3,cos(theta/2)+sin(theta/2)*e12);}{%
}
\end{mapleinput}

\mapleresult
\begin{maplelatex}
\mapleinline{inert}{2d}{
cos(theta)*e1+e2*sin(theta);-e1*sin(theta)+cos(theta)*e2;e3;}{%
\[
\mathrm{cos}(\theta )\,\mathit{e1} + \mathit{e2}\,\mathrm{sin}(\theta ),\;
 - \mathit{e1}\,\mathrm{sin}(\theta ) + \mathrm{cos}(\theta )\,\mathit{e2},\;
\mathit{e3}
\]
}
\end{maplelatex}
\end{maplegroup}
\noindent
Let's now take a general element from $\Spin(3)$ and act on all three unit basis vectors $\be_1,\be_2,\be_3.$ We can easily verify that the new elements $\be_{11},\be_{22},\be_{33}$ provide another orthonormal basis with the same orientation:
\begin{maplegroup}
\begin{mapleinput}
\mapleinline{active}{1d}{e11:=rot(e1,gSpin);e22:=rot(e2,gSpin);e33:=rot(e3,gSpin);}{%
}
\end{mapleinput}

\mapleresult
\begin{maplelatex}
\mapleinline{inert}{2d}{e11 :=
2*(eps*kappa*c2+c3*c1)*e3-(2*c2^2+2*c3^2-1)*e1+2*(eps*kappa*c3-c2*c1)*
e2;}{%
\[
\mathit{e11} := 2\,(\mathit{eps}\,\kappa \,\mathit{c2} + \mathit{
c3}\,\mathit{c1})\,\mathit{e3} - (2\,\mathit{c2}^{2} + 2\,
\mathit{c3}^{2} - 1)\,\mathit{e1} + 2\,(\mathit{eps}\,\kappa \,
\mathit{c3} - \mathit{c2}\,\mathit{c1})\,\mathit{e2}
\]
}
\end{maplelatex}

\begin{maplelatex}
\mapleinline{inert}{2d}{e22 :=
2*(eps*kappa*c1-c3*c2)*e3-2*(eps*kappa*c3+c2*c1)*e1-(2*c1^2-1+2*c3^2)*
e2;}{%
\[
\mathit{e22} := 2\,(\mathit{eps}\,\kappa \,\mathit{c1} - \mathit{
c3}\,\mathit{c2})\,\mathit{e3} - 2\,(\mathit{eps}\,\kappa \,
\mathit{c3} + \mathit{c2}\,\mathit{c1})\,\mathit{e1} - (2\,
\mathit{c1}^{2} - 1 + 2\,\mathit{c3}^{2})\,\mathit{e2}
\]
}
\end{maplelatex}

\begin{maplelatex}
\mapleinline{inert}{2d}{e33 :=
-(2*c2^2+2*c1^2-1)*e3-2*(-c3*c1+eps*kappa*c2)*e1-2*(eps*kappa*c1+c3*c2
)*e2;}{%
\[
\mathit{e33} :=  - (2\,\mathit{c2}^{2} + 2\,\mathit{c1}^{2} - 1)
\,\mathit{e3} - 2\,( - \mathit{c3}\,\mathit{c1} + \mathit{eps}\,
\kappa \,\mathit{c2})\,\mathit{e1} - 2\,(\mathit{eps}\,\kappa \,
\mathit{c1} + \mathit{c3}\,\mathit{c2})\,\mathit{e2}
\]
}
\end{maplelatex}

\end{maplegroup}
\begin{maplegroup}
\begin{mapleinput}
\mapleinline{active}{1d}{e1 &c e2 + e2 &c e1, e1 &c e3 + e3 &c e1, e2 &c e3 + e3 &c e2;}{%
}
\end{mapleinput}

\mapleresult
\begin{maplelatex}
\mapleinline{inert}{2d}{0, 0, 0;}{%
\[
0, \,0, \,0
\]
}
\end{maplelatex}

\end{maplegroup}
\begin{maplegroup}
\begin{mapleinput}
\mapleinline{active}{1d}{
e11 &c e22 + e22 &c e11, e11 &c e33 + e33 &c e11, e22 &c e33 + e33 &c e22;}{%
}
\end{mapleinput}

\mapleresult
\begin{maplelatex}
\mapleinline{inert}{2d}{0;0;0;}{%
\[
0,\,0,\,0
\]
}
\end{maplelatex}
\end{maplegroup}
\begin{maplegroup}
\begin{mapleinput}
\mapleinline{active}{1d}{e1 &w e2 &w e3,e11 &w e22 &w e33;}{%
}
\end{mapleinput}

\mapleresult
\begin{maplelatex}
\mapleinline{inert}{2d}{e123;e123;}{%
\[
\mathit{e123},\,\mathit{e123}
\]
}
\end{maplelatex}

\end{maplegroup}
\noindent
Length of a vector $\bv$ under the action $\Spin(3)$ is of course preserved:
\begin{maplegroup}
\begin{mapleinput}
\mapleinline{active}{1d}{Enorm(v),Enorm(Spin_action(v,gSpin));}{%
}
\end{mapleinput}

\mapleresult
\begin{maplelatex}
\mapleinline{inert}{2d}{c3^2+c1^2+c2^2;c3^2+c1^2+c2^2;}{%
\[
\mathit{c3}^{2} + \mathit{c1}^{2} + \mathit{c2}^{2},\;
\mathit{c3}^{2} + \mathit{c1}^{2} + \mathit{c2}^{2}
\]
}
\end{maplelatex}
\end{maplegroup}

\example{Example 1:}{Rotations in coordinate planes}

Let's define unit quaternions responsible for the rotations in the coordinate planes.  These are counter-clockwise rotations when looking down the rotation axis. We will define a pure-quaternion basis consisting of $\{\bq_i,\bq_j,\bq_k\}$ in place of traditionally used $\{\bi,\bj,\bk\}.$

\begin{maplegroup}
\begin{mapleinput}
\mapleinline{active}{1d}{qi:=e23:qj:=e13:qk:=e12:
q12:=cos(alpha/2)*Id+sin(alpha/2)*'qk';#rotation in the xy-plane}{%
}
\end{mapleinput}

\mapleresult
\begin{maplelatex}
\mapleinline{inert}{2d}{q12 := cos(1/2*alpha)*Id+sin(1/2*alpha)*qk;}{%
\[
\mathit{q12} := \mathrm{cos}({\displaystyle \frac {1}{2}} \,
\alpha )\,\mathit{Id} + \mathrm{sin}({\displaystyle \frac {1}{2}
} \,\alpha )\,\mathit{qk}
\]
}
\end{maplelatex}

\end{maplegroup}
\begin{maplegroup}
\begin{mapleinput}
\mapleinline{active}{1d}{q13:=cos(beta/2)*Id+sin(beta/2)*'qj';#rotation in the xz-plane}{%
}
\end{mapleinput}

\mapleresult
\begin{maplelatex}
\mapleinline{inert}{2d}{q13 := cos(1/2*beta)*Id+sin(1/2*beta)*qj;}{%
\[
\mathit{q13} := \mathrm{cos}({\displaystyle \frac {1}{2}} \,\beta
 )\,\mathit{Id} + \mathrm{sin}({\displaystyle \frac {1}{2}} \,
\beta )\,\mathit{qj}
\]
}
\end{maplelatex}

\end{maplegroup}
\begin{maplegroup}
\begin{mapleinput}
\mapleinline{active}{1d}{q23:=cos(gamma/2)*Id+sin(gamma/2)*'qi';#rotation in the yz-plane}{%
}
\end{mapleinput}

\mapleresult
\begin{maplelatex}
\mapleinline{inert}{2d}{q23 := cos(1/2*gamma)*Id+sin(1/2*gamma)*qi;}{%
\[
\mathit{q23} := \mathrm{cos}({\displaystyle \frac {1}{2}} \,
\gamma )\,\mathit{Id} + \mathrm{sin}({\displaystyle \frac {1}{2}
} \,\gamma )\,\mathit{qi}
\]
}
\end{maplelatex}

\end{maplegroup}
\noindent
Notice that to rotate by an angle $n\alpha$  it is enough to find the $n$-th Clifford power of the appropriate quaternion and then apply it to a vector.
\begin{maplegroup}
\begin{mapleinput}
\mapleinline{active}{1d}{q12 &c q12; #rotation by the angle 2*alpha}{%
}
\end{mapleinput}

\mapleresult
\begin{maplelatex}
\mapleinline{inert}{2d}{cos(alpha)*Id+e12*sin(alpha);}{%
\[
\mathrm{cos}(\alpha )\,\mathit{Id} + \mathit{e12}\,\mathrm{sin}(
\alpha )
\]
}
\end{maplelatex}

\end{maplegroup}
\begin{maplegroup}
\begin{mapleinput}
\mapleinline{active}{1d}{q12 &c q12 &c q12; #rotation by the angle 3*alpha}{%
}
\end{mapleinput}

\mapleresult
\begin{maplelatex}
\mapleinline{inert}{2d}{cos(3/2*alpha)*Id+sin(3/2*alpha)*e12;}{%
\[
\mathrm{cos}({\displaystyle \frac {3}{2}} \,\alpha )\,\mathit{Id}
 + \mathrm{sin}({\displaystyle \frac {3}{2}} \,\alpha )\,\mathit{
e12}
\]
}
\end{maplelatex}

\end{maplegroup}
\begin{maplegroup}
\begin{mapleinput}
\mapleinline{active}{1d}{q12 &c q12 &c q12 &c q12; #rotation by the angle 4*alpha}{%
}
\end{mapleinput}

\mapleresult
\begin{maplelatex}
\mapleinline{inert}{2d}{cos(2*alpha)*Id+sin(2*alpha)*e12;}{%
\[
\mathrm{cos}(2\,\alpha )\,\mathit{Id} + \mathrm{sin}(2\,\alpha )
\,\mathit{e12}
\]
}
\end{maplelatex}

\end{maplegroup}
\noindent
Let's see now how these basis rotations in the coordinate planes act on an arbitrary vector $\bv=a\e_1+b\e_2+c\e_3:$
\begin{maplegroup}
\begin{mapleinput}
\mapleinline{active}{1d}{v:=a*e1+b*e2+c*e3;}{%
} 
\end{mapleinput}

\mapleresult
\begin{maplelatex}
\mapleinline{inert}{2d}{v := a*e1+b*e2+c*e3;}{%
\[
v := a\,\mathit{e1} + b\,\mathit{e2} + c\,\mathit{e3}
\]
}
\end{maplelatex}

\end{maplegroup}
\noindent
The norm of $\bv$ is $\|\bv\|=\bv \Star{\bv}$ and it can be defined in {\sc CLIFFORD} as follows:

\begin{maplegroup}
\begin{mapleinput}
\mapleinline{active}{1d}{vlength:=sqrt(scalarpart(v &c star(v)));}{%
}
\end{mapleinput}

\mapleresult
\begin{maplelatex}
\mapleinline{inert}{2d}{vlength := sqrt(c^2+a^2+b^2);}{%
\[
\mathit{vlength} := \sqrt{c^{2} + a^{2} + b^{2}}
\]
}
\end{maplelatex}

\end{maplegroup}
\noindent
Certainly, rotations do not change length.  For example, let's rotate $\bv$ by  $q_{13} q_{12}:$
\begin{maplegroup}
\begin{mapleinput}
\mapleinline{active}{1d}{v123:=rot(v,q13 &c q12); #rotation q12 followed by q13}{%
}
\end{mapleinput}

\mapleresult
\begin{maplelatex}
\mapleinline{inert}{2d}{v123 :=
(-b*sin(beta)*sin(alpha)+a*sin(beta)*cos(alpha)+c*cos(beta))*e3+(-b*si
n(alpha)*cos(beta)-c*sin(beta)+a*cos(beta)*cos(alpha))*e1+(b*cos(alpha
)+a*sin(alpha))*e2;}{%
\maplemultiline{
\mathit{v123} := ( - b\,\mathrm{sin}(\beta )\,\mathrm{sin}(\alpha
 ) + a\,\mathrm{sin}(\beta )\,\mathrm{cos}(\alpha ) + c\,\mathrm{
cos}(\beta ))\,\mathit{e3} \\
\mbox{} + ( - b\,\mathrm{sin}(\alpha )\,\mathrm{cos}(\beta ) - c
\,\mathrm{sin}(\beta ) + a\,\mathrm{cos}(\beta )\,\mathrm{cos}(
\alpha ))\,\mathit{e1} + (b\,\mathrm{cos}(\alpha ) + a\,\mathrm{
sin}(\alpha ))\,\mathit{e2} }
}
\end{maplelatex}

\end{maplegroup}
\begin{maplegroup}
\begin{mapleinput}
\mapleinline{active}{1d}{vlength:=sqrt(scalarpart(v123 &c star(v123)));}{%
}
\end{mapleinput}

\mapleresult
\begin{maplelatex}
\mapleinline{inert}{2d}{vlength := sqrt(c^2+a^2+b^2);}{%
\[
\mathit{vlength} := \sqrt{c^{2} + a^{2} + b^{2}}
\]
}
\end{maplelatex}

\end{maplegroup}
\noindent
Thus, the length of $\bv_{123}$ is the same as the length of $\bv.$  However, rotations do not commute.  We will show that by applying quaternion $q_{12} q_{13}$ to $\bv$ and by comparing it with $q_{123}:$
\begin{maplegroup}
\begin{mapleinput}
\mapleinline{active}{1d}{v132:=rot(v,q12 &c q13); #rotation q13 followed by q12}{%
}
\end{mapleinput}

\mapleresult
\begin{maplelatex}
\mapleinline{inert}{2d}{v132 :=
(a*sin(beta)+c*cos(beta))*e3+(a*cos(beta)*cos(alpha)-c*sin(beta)*cos(a
lpha)-b*sin(alpha))*e1+(b*cos(alpha)-c*sin(beta)*sin(alpha)+a*sin(alph
a)*cos(beta))*e2;}{%
\maplemultiline{
\mathit{v132} := (a\,\mathrm{sin}(\beta ) + c\,\mathrm{cos}(\beta
 ))\,\mathit{e3} + (a\,\mathrm{cos}(\beta )\,\mathrm{cos}(\alpha 
) - c\,\mathrm{sin}(\beta )\,\mathrm{cos}(\alpha ) - b\,\mathrm{
sin}(\alpha ))\,\mathit{e1} \\
\mbox{} + (b\,\mathrm{cos}(\alpha ) - c\,\mathrm{sin}(\beta )\,
\mathrm{sin}(\alpha ) + a\,\mathrm{sin}(\alpha )\,\mathrm{cos}(
\beta ))\,\mathit{e2} }
}
\end{maplelatex}

\end{maplegroup}
\begin{maplegroup}
\begin{mapleinput}
\mapleinline{active}{1d}{clicollect(v123-v132);}{%
}
\end{mapleinput}

\mapleresult
\begin{maplelatex}
\mapleinline{inert}{2d}{(-b*sin(beta)*sin(alpha)+a*sin(beta)*cos(alpha)-a*sin(beta))*e3+(-b*s
in(alpha)*cos(beta)-c*sin(beta)+c*sin(beta)*cos(alpha)+b*sin(alpha))*e
1+(a*sin(alpha)+c*sin(beta)*sin(alpha)-a*sin(alpha)*cos(beta))*e2;}{%
\maplemultiline{
( - b\,\mathrm{sin}(\beta )\,\mathrm{sin}(\alpha ) + a\,\mathrm{
sin}(\beta )\,\mathrm{cos}(\alpha ) - a\,\mathrm{sin}(\beta ))\,
\mathit{e3} \\
\mbox{} + ( - b\,\mathrm{sin}(\alpha )\,\mathrm{cos}(\beta ) - c
\,\mathrm{sin}(\beta ) + c\,\mathrm{sin}(\beta )\,\mathrm{cos}(
\alpha ) + b\,\mathrm{sin}(\alpha ))\,\mathit{e1} \\
\mbox{} + (a\,\mathrm{sin}(\alpha ) + c\,\mathrm{sin}(\beta )\,
\mathrm{sin}(\alpha ) - a\,\mathrm{sin}(\alpha )\,\mathrm{cos}(
\beta ))\,\mathit{e2} }
}
\end{maplelatex}

\end{maplegroup}
\noindent
As it can be seen, $\bv_{123} \neq \bv_{132}.$

\example{Example 2:}{Counter-clockwise rotation by an angle $\alpha$ around the given axis}

In this example we will find a way to rotate a given vector $\bv$ by an angle $\alpha$ in a plane orthogonal to the given axis vector $\axis=a_1\be_1+a_2\be_2+a_3\be_3$ vector. This rotation will be counter-clockwise when looking down the axis towards to origin $(0,0,0)$ of the coordinate system. In order to derive symbolic formulas, we will assume that the symbolic vector $\axis$ has been normalized by defining $\lambda=\pm \sqrt{1-a_1^2-a_2^2}$ and $\axis = a_1\be_1+a_2\be_2+\lambda \be_3.$ However, it won't be necessary for $\axis$ to be of unit length when its components are numeric.
\begin{maplegroup}
\begin{mapleinput}
\mapleinline{active}{1d}{alias(lambda=RootOf(-a1^2-a2^2-_Z^2+1)):axis:=a1*e1+a2*e2+lambda*e3;}{%
}
\end{mapleinput}

\mapleresult
\begin{maplelatex}
\mapleinline{inert}{2d}{axis := a1*e1+a2*e2+lambda*e3;}{%
\[
\mathit{axis} := \mathit{a1}\,\mathit{e1} + \mathit{a2}\,\mathit{e2} + \lambda \,\mathit{e3}
\]
}
\end{maplelatex}

\end{maplegroup}
\begin{maplegroup}
\begin{mapleinput}
\mapleinline{active}{1d}{simplify(axis &c star(axis));}{%
}
\end{mapleinput}

\mapleresult
\begin{maplelatex}
\mapleinline{inert}{2d}{Id;}{%
\[
\mathit{Id}
\]
}
\end{maplelatex}

\end{maplegroup}
\noindent
Thus, in the symbolic case we will always have that $\|\axis\|=1.$ Notice that in order to represent a rotation around $\axis$ we need to find a dual unit quaternion which we will call $\qaxis.$ It will need to be defined in such a way as to give a desired orientation for the rotation. Since we have opted for counter-clockwise rotations, we define $\qaxis = - \axis \, \e_{123}$ where $\e_{123}$ is a unit pseudo-scalar in $\cl_{0,3}.$ In the following we will refer to the $\axis = a_1\be_1+a_2\be_2+a_3\be_3$ as the axis $(a_1,a_2,a_3).$
\begin{maplegroup}
\begin{mapleinput}
\mapleinline{active}{1d}{qaxis:=axis &c (-e123);}{%
}
\end{mapleinput}

\mapleresult
\begin{maplelatex}
\mapleinline{inert}{2d}{qaxis := lambda*e12+a1*e23-a2*e13;}{%
\[
\mathit{qaxis} := \lambda \,\mathit{e12} + \mathit{a1}\,\mathit{
e23} - \mathit{a2}\,\mathit{e13}
\]
}
\end{maplelatex}

\end{maplegroup}
\noindent
In Appendix 3 Reader can find a procedure {\tt qrot} which finds the dual quaternion $\qaxis.$ The first three arguments to {\tt qrot} are the (numeric or symbolic) components of the $\axis$ vector in the basis $\{\be_1,\be_2,\be_3\}$ while the fourth argument is the angle of rotation. For example, we can define various rotation quaternions:

\begin{maplegroup}
\begin{mapleinput}
\mapleinline{active}{1d}{q100:=qrot(1,0,0,theta);#rotation about the axis (1,0,0)}{%
}
\end{mapleinput}

\mapleresult
\begin{maplelatex}
\mapleinline{inert}{2d}{q100 := cos(1/2*theta)*Id+sin(1/2*theta)*e23;}{%
\[
\mathit{q100} := \mathrm{cos}({\displaystyle \frac {1}{2}} \,
\theta )\,\mathit{Id} + \mathrm{sin}({\displaystyle \frac {1}{2}
} \,\theta )\,\mathit{e23}
\]
}
\end{maplelatex}

\end{maplegroup}
\begin{maplegroup}
\begin{mapleinput}
\mapleinline{active}{1d}{q010:=qrot(0,1,0,theta);#rotation about the axis (0,1,0)}{%
}
\end{mapleinput}

\mapleresult
\begin{maplelatex}
\mapleinline{inert}{2d}{q010 := cos(1/2*theta)*Id-sin(1/2*theta)*e13;}{%
\[
\mathit{q010} := \mathrm{cos}({\displaystyle \frac {1}{2}} \,
\theta )\,\mathit{Id} - \mathrm{sin}({\displaystyle \frac {1}{2}
} \,\theta )\,\mathit{e13}
\]
}
\end{maplelatex}

\end{maplegroup}
\begin{maplegroup}
\begin{mapleinput}
\mapleinline{active}{1d}{q001:=qrot(0,0,1,theta);#rotation about the axis (0,0,1)}{%
}
\end{mapleinput}

\mapleresult
\begin{maplelatex}
\mapleinline{inert}{2d}{q001 := cos(1/2*theta)*Id+sin(1/2*theta)*e12;}{%
\[
\mathit{q001} := \mathrm{cos}({\displaystyle \frac {1}{2}} \,
\theta )\,\mathit{Id} + \mathrm{sin}({\displaystyle \frac {1}{2}
} \,\theta )\,\mathit{e12}
\]
}
\end{maplelatex}

\end{maplegroup}
\begin{maplegroup}
\begin{mapleinput}
\mapleinline{active}{1d}{q101:=qrot(1,0,1,theta);#rotation about the axis (1,0,1)}{%
}
\end{mapleinput}

\mapleresult
\begin{maplelatex}
\mapleinline{inert}{2d}{q101 :=
cos(1/2*theta)*Id+sin(1/2*theta)*(1/2*e12*sqrt(2)+1/2*sqrt(2)*e23);}{%
\[
\mathit{q101} := \mathrm{cos}({\displaystyle \frac {1}{2}} \,
\theta )\,\mathit{Id} + \mathrm{sin}({\displaystyle \frac {1}{2}
} \,\theta )\,({\displaystyle \frac {1}{2}} \,\mathit{e12}\,
\sqrt{2} + {\displaystyle \frac {1}{2}} \,\sqrt{2}\,\mathit{e23})
\]
}
\end{maplelatex}

\end{maplegroup}
\begin{maplegroup}
\begin{mapleinput}
\mapleinline{active}{1d}{q011:=qrot(0,1,1,theta);#rotation about the axis (0,1,1)}{%
}
\end{mapleinput}

\mapleresult
\begin{maplelatex}
\mapleinline{inert}{2d}{q011 :=
cos(1/2*theta)*Id+sin(1/2*theta)*(1/2*e12*sqrt(2)-1/2*sqrt(2)*e13);}{%
\[
\mathit{q011} := \mathrm{cos}({\displaystyle \frac {1}{2}} \,
\theta )\,\mathit{Id} + \mathrm{sin}({\displaystyle \frac {1}{2}
} \,\theta )\,({\displaystyle \frac {1}{2}} \,\mathit{e12}\,
\sqrt{2} - {\displaystyle \frac {1}{2}} \,\sqrt{2}\,\mathit{e13})
\]
}
\end{maplelatex}

\end{maplegroup}
\begin{maplegroup}
\begin{mapleinput}
\mapleinline{active}{1d}{q110:=qrot(1,1,0,theta);#rotation about the axis (1,1,0)}{%
}
\end{mapleinput}

\mapleresult
\begin{maplelatex}
\mapleinline{inert}{2d}{q110 :=
cos(1/2*theta)*Id+sin(1/2*theta)*(1/2*sqrt(2)*e23-1/2*sqrt(2)*e13);}{%
\[
\mathit{q110} := \mathrm{cos}({\displaystyle \frac {1}{2}} \,
\theta )\,\mathit{Id} + \mathrm{sin}({\displaystyle \frac {1}{2}
} \,\theta )\,({\displaystyle \frac {1}{2}} \,\sqrt{2}\,\mathit{
e23} - {\displaystyle \frac {1}{2}} \,\sqrt{2}\,\mathit{e13})
\]
}
\end{maplelatex}

\end{maplegroup}
\begin{maplegroup}
\begin{mapleinput}
\mapleinline{active}{1d}{q111:=qrot(1,1,1,theta);#rotation about the axis (1,1,1)}{%
}
\end{mapleinput}

\mapleresult
\begin{maplelatex}
\mapleinline{inert}{2d}{
q111 := cos(1/2*theta)*Id+sin(1/2*theta)*(1/3*e12*sqrt(3)-1/3*sqrt(3)*e13+1/3*sqrt(3)*e23);}{%
\[
\mathit{q111} := \mathrm{cos}({\displaystyle \frac {1}{2}} \,
\theta )\,\mathit{Id} + \mathrm{sin}({\displaystyle \frac {1}{2}
} \,\theta )\,({\displaystyle \frac {1}{3}} \,\mathit{e12}\,
\sqrt{3} - {\displaystyle \frac {1}{3}} \,\sqrt{3}\,\mathit{e13}
 + {\displaystyle \frac {1}{3}} \,\sqrt{3}\,\mathit{e23})
\]
}
\end{maplelatex}

\end{maplegroup}

Let's try to rotate first a basis vector $\e_1$ around various axes listed above by some angle $\alpha.$ In the next example, after we find a general formula for the components of the rotated e1, we will verify it for various angles.

\begin{maplegroup}
\begin{mapleinput}
\mapleinline{active}{1d}{v:=e1:vnew:=rot(v,q100); #rotation around the (100) axis}{%
}
\end{mapleinput}

\mapleresult
\begin{maplelatex}
\mapleinline{inert}{2d}{vnew := e1;}{%
\[
\mathit{vnew} := \mathit{e1}
\]
}
\end{maplelatex}

\end{maplegroup}
\begin{maplegroup}
\begin{mapleinput}
\mapleinline{active}{1d}{vnew:=rot(v,q010); #rotation around the (010) axis}{%
}
\end{mapleinput}

\mapleresult
\begin{maplelatex}
\mapleinline{inert}{2d}{vnew := -sin(theta)*e3+cos(theta)*e1;}{%
\[
\mathit{vnew} :=  - \mathrm{sin}(\theta )\,\mathit{e3} + \mathrm{
cos}(\theta )\,\mathit{e1}
\]
}
\end{maplelatex}

\end{maplegroup}
\begin{maplegroup}
\begin{mapleinput}
\mapleinline{active}{1d}{eval(subs(theta=Pi/2,vnew));}{%
}
\end{mapleinput}

\mapleresult
\begin{maplelatex}
\mapleinline{inert}{2d}{-e3;}{%
\[
 - \mathit{e3}
\]
}
\end{maplelatex}

\end{maplegroup}
\begin{maplegroup}
\begin{mapleinput}
\mapleinline{active}{1d}{vnew:=rot(v,q001); #rotation around the (001) axis}{%
}
\end{mapleinput}

\mapleresult
\begin{maplelatex}
\mapleinline{inert}{2d}{vnew := cos(theta)*e1+e2*sin(theta);}{%
\[
\mathit{vnew} := \mathrm{cos}(\theta )\,\mathit{e1} + \mathit{e2}
\,\mathrm{sin}(\theta )
\]
}
\end{maplelatex}

\end{maplegroup}
\begin{maplegroup}
\begin{mapleinput}
\mapleinline{active}{1d}{eval(subs(theta=Pi/2,vnew));}{%
}
\end{mapleinput}

\mapleresult
\begin{maplelatex}
\mapleinline{inert}{2d}{e2;}{%
\[
\mathit{e2}
\]
}
\end{maplelatex}

\end{maplegroup}
\begin{maplegroup}
\begin{mapleinput}
\mapleinline{active}{1d}{vnew:=rot(v,q101); #rotation around the (101) axis}{%
}
\end{mapleinput}

\mapleresult
\begin{maplelatex}
\mapleinline{inert}{2d}{vnew :=
-1/2*(-1+cos(theta))*e3+1/2*(cos(theta)+1)*e1+1/2*sqrt(2)*e2*sin(theta
);}{%
\[
\mathit{vnew} :=  - {\displaystyle \frac {1}{2}} \,( - 1 + 
\mathrm{cos}(\theta ))\,\mathit{e3} + {\displaystyle \frac {1}{2}
} \,(\mathrm{cos}(\theta ) + 1)\,\mathit{e1} + {\displaystyle 
\frac {1}{2}} \,\sqrt{2}\,\mathit{e2}\,\mathrm{sin}(\theta )
\]
}
\end{maplelatex}

\end{maplegroup}
\begin{maplegroup}
\begin{mapleinput}
\mapleinline{active}{1d}{eval(subs(theta=Pi/2,vnew));}{%
}
\end{mapleinput}

\mapleresult
\begin{maplelatex}
\mapleinline{inert}{2d}{1/2*e3+1/2*e1+1/2*e2*sqrt(2);}{%
\[
{\displaystyle \frac {1}{2}} \,\mathit{e3} + {\displaystyle 
\frac {1}{2}} \,\mathit{e1} + {\displaystyle \frac {1}{2}} \,
\mathit{e2}\,\sqrt{2}
\]
}
\end{maplelatex}

\end{maplegroup}
\begin{maplegroup}
\begin{mapleinput}
\mapleinline{active}{1d}{vnew:=rot(v,q011); #rotation around the (011) axis}{%
}
\end{mapleinput}

\mapleresult
\begin{maplelatex}
\mapleinline{inert}{2d}{vnew :=
-1/2*sqrt(2)*e3*sin(theta)+cos(theta)*e1+1/2*sqrt(2)*e2*sin(theta);}{%
\[
\mathit{vnew} :=  - {\displaystyle \frac {1}{2}} \,\sqrt{2}\,
\mathit{e3}\,\mathrm{sin}(\theta ) + \mathrm{cos}(\theta )\,
\mathit{e1} + {\displaystyle \frac {1}{2}} \,\sqrt{2}\,\mathit{e2
}\,\mathrm{sin}(\theta )
\]
}
\end{maplelatex}

\end{maplegroup}
\begin{maplegroup}
\begin{mapleinput}
\mapleinline{active}{1d}{eval(subs(theta=Pi/2,vnew));}{%
}
\end{mapleinput}

\mapleresult
\begin{maplelatex}
\mapleinline{inert}{2d}{-1/2*e3*sqrt(2)+1/2*e2*sqrt(2);}{%
\[
 - {\displaystyle \frac {1}{2}} \,\mathit{e3}\,\sqrt{2} + 
{\displaystyle \frac {1}{2}} \,\mathit{e2}\,\sqrt{2}
\]
}
\end{maplelatex}

\end{maplegroup}
\begin{maplegroup}
\begin{mapleinput}
\mapleinline{active}{1d}{vnew:=rot(v,q110); #rotation around the (110) axis}{%
}
\end{mapleinput}

\mapleresult
\begin{maplelatex}
\mapleinline{inert}{2d}{vnew :=
-1/2*sqrt(2)*e3*sin(theta)+1/2*(cos(theta)+1)*e1-1/2*(-1+cos(theta))*e
2;}{%
\[
\mathit{vnew} :=  - {\displaystyle \frac {1}{2}} \,\sqrt{2}\,
\mathit{e3}\,\mathrm{sin}(\theta ) + {\displaystyle \frac {1}{2}
} \,(\mathrm{cos}(\theta ) + 1)\,\mathit{e1} - {\displaystyle 
\frac {1}{2}} \,( - 1 + \mathrm{cos}(\theta ))\,\mathit{e2}
\]
}
\end{maplelatex}

\end{maplegroup}
\begin{maplegroup}
\begin{mapleinput}
\mapleinline{active}{1d}{eval(subs(theta=Pi/2,vnew));}{%
}
\end{mapleinput}

\mapleresult
\begin{maplelatex}
\mapleinline{inert}{2d}{-1/2*e3*sqrt(2)+1/2*e1+1/2*e2;}{%
\[
 - {\displaystyle \frac {1}{2}} \,\mathit{e3}\,\sqrt{2} + 
{\displaystyle \frac {1}{2}} \,\mathit{e1} + {\displaystyle 
\frac {1}{2}} \,\mathit{e2}
\]
}
\end{maplelatex}

\end{maplegroup}
\begin{maplegroup}
\begin{mapleinput}
\mapleinline{active}{1d}{vnew:=rot(v,q111); #rotation around the (111) axis}{%
}
\end{mapleinput}

\mapleresult
\begin{maplelatex}
\mapleinline{inert}{2d}{vnew :=
-1/3*(sqrt(3)*sin(theta)-1+cos(theta))*e3+1/3*(2*cos(theta)+1)*e1+1/3*
(sqrt(3)*sin(theta)+1-cos(theta))*e2;}{%
\maplemultiline{
\mathit{vnew} :=  - {\displaystyle \frac {1}{3}} \,(\sqrt{3}\,
\mathrm{sin}(\theta ) - 1 + \mathrm{cos}(\theta ))\,\mathit{e3}
 + {\displaystyle \frac {1}{3}} \,(2\,\mathrm{cos}(\theta ) + 1)
\,\mathit{e1} \\
\mbox{} + {\displaystyle \frac {1}{3}} \,(\sqrt{3}\,\mathrm{sin}(
\theta ) + 1 - \mathrm{cos}(\theta ))\,\mathit{e2} }
}
\end{maplelatex}

\end{maplegroup}
\begin{maplegroup}
\begin{mapleinput}
\mapleinline{active}{1d}{eval(subs(theta=Pi/2,vnew));}{%
}
\end{mapleinput}

\mapleresult
\begin{maplelatex}
\mapleinline{inert}{2d}{-1/3*(sqrt(3)-1)*e3+1/3*e1+1/3*(sqrt(3)+1)*e2;}{%
\[
 - {\displaystyle \frac {1}{3}} \,(\sqrt{3} - 1)\,\mathit{e3} + 
{\displaystyle \frac {1}{3}} \,\mathit{e1} + {\displaystyle 
\frac {1}{3}} \,(\sqrt{3} + 1)\,\mathit{e2}
\]
}
\end{maplelatex}

\end{maplegroup}
\begin{maplegroup}
\begin{mapleinput}
\mapleinline{active}{1d}{eval(subs(theta=Pi,vnew));}{%
}
\end{mapleinput}

\mapleresult
\begin{maplelatex}
\mapleinline{inert}{2d}{2/3*e3-1/3*e1+2/3*e2;}{%
\[
{\displaystyle \frac {2}{3}} \,\mathit{e3} - {\displaystyle 
\frac {1}{3}} \,\mathit{e1} + {\displaystyle \frac {2}{3}} \,
\mathit{e2}
\]
}
\end{maplelatex}

\end{maplegroup}

\example{Example 3:}{A general formula} 

Finally, we derive a general formula for a rotation of an arbitrary vector $\bv=v_1\be_1+v_2\be_2+v_3\be_3$ around an arbitrary $\axis=a_1\be_1+a_2\be_2+a_3\be_3$ and by an arbitrary angle $\alpha.$ In the purely symbolic case we assume that the axis is of unit length, that is, $a_3=\lambda.$

\begin{maplegroup}
\begin{mapleinput}
\mapleinline{active}{1d}{v:=v1*e1+v2*e2+v3*e3; #an arbitrary vector}{%
}
\end{mapleinput}

\mapleresult
\begin{maplelatex}
\mapleinline{inert}{2d}{v := v1*e1+v2*e2+v3*e3;}{%
\[
v := \mathit{v1}\,\mathit{e1} + \mathit{v2}\,\mathit{e2} + 
\mathit{v3}\,\mathit{e3}
\]
}
\end{maplelatex}

\end{maplegroup}
\begin{maplegroup}
\begin{mapleinput}
\mapleinline{active}{1d}{
qnew:=clicollect(rot(v,qrot(a1,a2,lambda,theta)));# new vector after rotation}{%
}
\end{mapleinput}

\mapleresult
\begin{maplelatex}
\mapleinline{inert}{2d}{qnew :=
(-lambda*v1*a1*cos(theta)-lambda*v2*a2*cos(theta)+a1^2*v3*cos(theta)-v
1*a2*sin(theta)+lambda*v1*a1+lambda*v2*a2+a2^2*v3*cos(theta)+v2*a1*sin
(theta)+v3-a2^2*v3-a1^2*v3)*e3+(v1*cos(theta)-a1*v2*a2*cos(theta)-lamb
da*v3*a1*cos(theta)+lambda*v3*a1+a1^2*v1-a1^2*v1*cos(theta)-v2*lambda*
sin(theta)+v3*a2*sin(theta)+a1*v2*a2)*e1+(a2^2*v2+a1*v1*a2-v3*a1*sin(t
heta)-a2^2*v2*cos(theta)+lambda*v3*a2+v1*lambda*sin(theta)-lambda*v3*a
2*cos(theta)+v2*cos(theta)-a1*v1*a2*cos(theta))*e2;}{%
\maplemultiline{
\mathit{qnew} := ( - \lambda \,\mathit{v1}\,\mathit{a1}\,\mathrm{
cos}(\theta ) - \lambda \,\mathit{v2}\,\mathit{a2}\,\mathrm{cos}(
\theta ) + \mathit{a1}^{2}\,\mathit{v3}\,\mathrm{cos}(\theta ) - 
\mathit{v1}\,\mathit{a2}\,\mathrm{sin}(\theta ) + \lambda \,
\mathit{v1}\,\mathit{a1} \\
\mbox{} + \lambda \,\mathit{v2}\,\mathit{a2} + \mathit{a2}^{2}\,
\mathit{v3}\,\mathrm{cos}(\theta ) + \mathit{v2}\,\mathit{a1}\,
\mathrm{sin}(\theta ) + \mathit{v3} - \mathit{a2}^{2}\,\mathit{v3
} - \mathit{a1}^{2}\,\mathit{v3})\mathit{e3}\mbox{} + (\mathit{v1
}\,\mathrm{cos}(\theta ) \\
\mbox{} - \mathit{a1}\,\mathit{v2}\,\mathit{a2}\,\mathrm{cos}(
\theta ) - \lambda \,\mathit{v3}\,\mathit{a1}\,\mathrm{cos}(
\theta ) + \lambda \,\mathit{v3}\,\mathit{a1} + \mathit{a1}^{2}\,
\mathit{v1} - \mathit{a1}^{2}\,\mathit{v1}\,\mathrm{cos}(\theta )
 \\
\mbox{} - \mathit{v2}\,\lambda \,\mathrm{sin}(\theta ) + \mathit{
v3}\,\mathit{a2}\,\mathrm{sin}(\theta ) + \mathit{a1}\,\mathit{v2
}\,\mathit{a2})\mathit{e1}\mbox{} + (\mathit{a2}^{2}\,\mathit{v2}
 + \mathit{a1}\,\mathit{v1}\,\mathit{a2} - \mathit{v3}\,\mathit{
a1}\,\mathrm{sin}(\theta ) \\
\mbox{} - \mathit{a2}^{2}\,\mathit{v2}\,\mathrm{cos}(\theta ) + 
\lambda \,\mathit{v3}\,\mathit{a2} + \mathit{v1}\,\lambda \,
\mathrm{sin}(\theta ) - \lambda \,\mathit{v3}\,\mathit{a2}\,
\mathrm{cos}(\theta ) + \mathit{v2}\,\mathrm{cos}(\theta ) \\
\mbox{} - \mathit{a1}\,\mathit{v1}\,\mathit{a2}\,\mathrm{cos}(
\theta ))\mathit{e2} }
}
\end{maplelatex}

\end{maplegroup}
\noindent
For example, let's rotate vector $\bv=\e_1+2 \e_2+3\e_3$ around the axis $(1,2,3)$ by
any angle $\alpha.$  Certainly, since the vector $\bv$ is on the axis of rotation, it should not change: 
\begin{maplegroup}
\begin{mapleinput}
\mapleinline{active}{1d}{clicollect(rot(e1+2*e2+3*e3,qrot(1,2,3,alpha)));}{%
}
\end{mapleinput}

\mapleresult
\begin{maplelatex}
\mapleinline{inert}{2d}{e1+2*e2+3*e3;}{%
\[
\mathit{e1} + 2\,\mathit{e2} + 3\,\mathit{e3}
\]
}
\end{maplelatex}

\end{maplegroup}
\noindent
Let's rotate $\bv=\e_1-2 \e_2+4 \e_3$ around the axis $(2,-3,4)$ by an angle $\alpha=\pi/4.$

\begin{maplegroup}
\begin{mapleinput}
\mapleinline{active}{1d}{clicollect(rot(e1-2*e2+4*e3,qrot(2,-3,4,Pi/4)));}{%
}
\end{mapleinput}

\mapleresult
\begin{maplelatex}
\mapleinline{inert}{2d}{(-1/58*sqrt(29)*sqrt(2)+10/29*sqrt(2)+96/29)*e3+(-2/29*sqrt(29)*sqrt(
2)+48/29-19/58*sqrt(2))*e1+(7/29*sqrt(2)-2/29*sqrt(29)*sqrt(2)-72/29)*
e2;}{%
\maplemultiline{
( - {\displaystyle \frac {1}{58}} \,\sqrt{29}\,\sqrt{2} + 
{\displaystyle \frac {10}{29}} \,\sqrt{2} + {\displaystyle 
\frac {96}{29}} )\,\mathit{e3} + ( - {\displaystyle \frac {2}{29}
} \,\sqrt{29}\,\sqrt{2} + {\displaystyle \frac {48}{29}}  - 
{\displaystyle \frac {19}{58}} \,\sqrt{2})\,\mathit{e1} \\
\mbox{} + ({\displaystyle \frac {7}{29}} \,\sqrt{2} - 
{\displaystyle \frac {2}{29}} \,\sqrt{29}\,\sqrt{2} - 
{\displaystyle \frac {72}{29}} )\,\mathit{e2} }
}
\end{maplelatex}

\end{maplegroup}

Thus, in this section we have shown how easy it is to derive vector rotation formulas from vector analysis using elements of $\Spin(3)$ considered as unit quaternions. It has been very helpful to be able to embed $\Spin(3)$ in $\cl_{0,3}.$

\subsection{Degenerate Clifford algebra and the proper rigid motions $SE(3)$}
\label{sec:e3}

In this final section we will use the ability of {\sc CLIFFORD} to perform computations in Clifford algebras of an arbitrary quadratic form including, of course, degenerate forms. We will consider the semi-direct product $\Spin(3) \rtimes \BR^3$ that double covers the group of proper rigid motions $SE(3).$ We will follow the notation used in \cite{Selig}, page 156, except that our basis vector that squares to $0$ will be $\be_4$ and not $\be.$ We begin by defining $B$ as a degenerate diagonal form $\diag(-1,-1,-1,0)$ of signature $(0,3,1).$ Recall from the previous section that procedure {\tt star} gives conjugation in $\cl_{0,3,1}=C(0,3,1).$
\begin{maplegroup}
\begin{mapleinput}
\mapleinline{active}{1d}{
dim:=4:n:=dim-1:eval(makealiases(dim)):B:=diag(-1$n,0); #Selig's C(0,3,1)}{%
}
\end{mapleinput}

\mapleresult
\begin{maplelatex}
\mapleinline{inert}{2d}{B := matrix([[-1, 0, 0, 0], [0, -1, 0, 0], [0, 0, -1, 0], [0, 0, 0,
0]]);}{%
\[
B :=  \left[ 
{\begin{array}{rrrr}
-1 & 0 & 0 & 0 \\
0 & -1 & 0 & 0 \\
0 & 0 & -1 & 0 \\
0 & 0 & 0 & 0
\end{array}}
 \right] 
\]
}
\end{maplelatex}
\end{maplegroup}
\noindent
Let the vector basis in $\BR^3$ be stored in {\tt vbasis} and let $\bt$ an arbitrary vector in $\BR^3.$ 
\begin{maplegroup}
\begin{mapleinput}
\mapleinline{active}{1d}{vbasis:=cbasis(n,1):t:=add(t.i*vbasis[i],i=1..n);}{%
}
\end{mapleinput}

\mapleresult
\begin{maplelatex}
\mapleinline{inert}{2d}{t := t1*e1+t2*e2+t3*e3;}{%
\[
t := \mathit{t1}\,\mathit{e1} + \mathit{t2}\,\mathit{e2} + 
\mathit{t3}\,\mathit{e3}
\]
}
\end{maplelatex}

\end{maplegroup}
\noindent
Elements of the form $\Id + \bt\,\e_4$ are invertible in $\cl_{0,3,1}$ since $u\, \e_4$ and 
$\e_4 \, u$ are nilpotent for any $u$ in $\cl_{0,3,1}.$ That is, the Jacobson radical  $J$ in $\cl_{0,3,1}$ is generated by $\be_4.$
\begin{maplegroup}
\begin{mapleinput}
\mapleinline{active}{1d}{cinv(1+t &c e4);#symbolic inverse of 1+t \& c e4}{%
}
\end{mapleinput}

\mapleresult
\begin{maplelatex}
\mapleinline{inert}{2d}{Id-t1*e14-t2*e24-t3*e34;}{%
\[
\mathit{Id} - \mathit{t1}\,\mathit{e14} - \mathit{t2}\,\mathit{
e24} - \mathit{t3}\,\mathit{e34}
\]
}
\end{maplelatex}

\end{maplegroup}

We will verify now statements made on page 156. Let $\kappa=\sqrt{-c_1^2-c_2^2-c_3^2+1}$ and $\varepsilon=\pm 1$ as before. Then the most general element $\bg$ in $\Spin(3)$ has the form:
\begin{maplegroup}
\begin{mapleinput}
\mapleinline{active}{1d}{alias(kappa=sqrt(-c1^2-c2^2-c3^2+1)):alias(eps=RootOf(_Z^2-1)):
gSpin:=eps*kappa*Id+c3*e12+c2*e13+c1*e23;}{%
}
\end{mapleinput}

\mapleresult
\begin{maplelatex}
\mapleinline{inert}{2d}{gSpin := eps*kappa*Id+c3*e12+c2*e13+c1*e23;}{%
\[
\mathit{gSpin} := \mathit{eps}\,\kappa \,\mathit{Id} + \mathit{c3
}\,\mathit{e12} + \mathit{c2}\,\mathit{e13} + \mathit{c1}\,
\mathit{e23}
\]
}
\end{maplelatex}

\end{maplegroup}
\noindent
We consider a subgroup $G$ of the group of units of $\cl_{0,3,1}$ of the form 
$\bg+\frac12 \bt \bg \e_4$ where $\bg$ belongs to $\Spin(3)$ and $\bt$ is a $1$-vector. Elements in $G$ will be given by a procedure {\tt ge}.
\begin{maplegroup}
\begin{mapleinput}
\mapleinline{active}{1d}{ge:=proc(g,t) RETURN(clicollect(simplify(g+1/2*t &c g &c e4))) end:}{%
}
\end{mapleinput}

\end{maplegroup}
\noindent
For the most general $\bg$ in $\Spin(3)$ and $\bt \in \BR^3,$ procedure {\tt ge} gives:
\begin{maplegroup}
\begin{mapleinput}
\mapleinline{active}{1d}{'ge(gSpin,t)'=ge(gSpin,t);}{%
}
\end{mapleinput}

\mapleresult
\begin{maplelatex}
\mapleinline{inert}{2d}{ge(gSpin,t)=(1/2*t2*kappa*eps-1/2*t1*c3+1/2*t3*c1)*e24+(1/2*t3*kappa*eps-1/2*t2*c1-1/2*t1*c2)*e34+(1/2*t2*c3+1/2*t3*c2+1/2*t1*kappa*eps)*e14+
(-1/2*t2*c2+1/2*t3*c3+1/2*t1*c1)*e1234+eps*kappa*Id+c3*e12+c2*e13+c1*e23;}{%
\maplemultiline{
\mathit{ge(gSpin,t)}=({\displaystyle \frac {1}{2}} \,\mathit{t2}\,\kappa \,\mathit{eps
} - {\displaystyle \frac {1}{2}} \,\mathit{t1}\,\mathit{c3} + 
{\displaystyle \frac {1}{2}} \,\mathit{t3}\,\mathit{c1})\,
\mathit{e24} + ({\displaystyle \frac {1}{2}} \,\mathit{t3}\,
\kappa \,\mathit{eps} - {\displaystyle \frac {1}{2}} \,\mathit{t2
}\,\mathit{c1} - {\displaystyle \frac {1}{2}} \,\mathit{t1}\,
\mathit{c2})\,\mathit{e34} \\
\mbox{} + ({\displaystyle \frac {1}{2}} \,\mathit{t2}\,\mathit{c3
} + {\displaystyle \frac {1}{2}} \,\mathit{t3}\,\mathit{c2} + 
{\displaystyle \frac {1}{2}} \,\mathit{t1}\,\kappa \,\mathit{eps}
)\,\mathit{e14} + ( - {\displaystyle \frac {1}{2}} \,\mathit{t2}
\,\mathit{c2} + {\displaystyle \frac {1}{2}} \,\mathit{t3}\,
\mathit{c3} + {\displaystyle \frac {1}{2}} \,\mathit{t1}\,
\mathit{c1})\,\mathit{e1234} \\
\mbox{} + \mathit{eps}\,\kappa \,\mathit{Id} + \mathit{c3}\,
\mathit{e12} + \mathit{c2}\,\mathit{e13} + \mathit{c1}\,\mathit{
e23} }
}
\end{maplelatex}

\end{maplegroup}
\noindent
First, let's verify Selig's statement that 
\begin{equation}
\Star{(\bg-\frac12 \bt \,\bg \,\e_4)} = (\Star{\bg}+\frac12 \Star{\bg}\,\bt\,\e_4).
\label{eq:st1}
\end{equation}
Notice that the left-hand-side in (\ref{eq:st1}) is just the conjugation of $ge(gSpin,-t)$ while the right-hand-side is equal to $ge(\Star{gSpin},-t)$ where $\Star{gSpin}$ denotes the conjugate of $gSpin.$
\begin{maplegroup}
\begin{mapleinput}
\mapleinline{active}{1d}{L:=clicollect(star(ge(gSpin,-t))):
R:=simplify(star(gSpin)+1/2*star(gSpin) &c t &c e4):
simplify(L-R);}{%
}
\end{mapleinput}

\mapleresult
\begin{maplelatex}
\mapleinline{inert}{2d}{0;}{%
\[
0
\]
}
\end{maplelatex}

\end{maplegroup}
\noindent
Next, we define the action of the group $G$ on the subspace of $\cl_{0,3,1}$ consisting of the elements of the form $1 + \bx \,\e_4$ as follows:
\begin{equation}
(\bg + \frac12 \bt\,\bg\,\e_4)(1 + \bx \,\e_4)\Star{(\bg - \frac12\bt\,\bg\e_4)} = 
1 + (\bg\,\bx \,\Star{\bg} + \bt)\,\e_4
\label{eq:st2}
\end{equation}
where $\bx,\bt \in \BR^3$ and $\bg \in \Spin(3).$ The above identity can be shown as follows. We define a procedure {\tt rigid} which will give this action on $\BR^3:$
$\bx \mapsto \bg\,\bx \,\Star{\bg} + \bt.$

\begin{maplegroup}
\begin{mapleinput}
\mapleinline{active}{1d}{rigid:=proc(x,g,t) local p;
if not evalb(x=vectorpart(x,1)) or not evalb(t=vectorpart(t,1)) then 
\phtab{} ERROR(`x and t must be vectors`)
fi:
if not type(g,evenelement) then ERROR(`g must be even`) fi;
RETURN(clicollect(simplify(cmul(g,x,star(g))+t)))
end:}{%
}
\end{mapleinput}

\end{maplegroup}
\noindent
This action will give us the rigid motion on $\BR^3.$ Thus, we can compute the right-hand-side of (\ref{eq:st2}) by using {\tt rigid} while the left-hand-side will be computed directly.
\begin{maplegroup}
\begin{mapleinput}
\mapleinline{active}{1d}{x:=add(x.i*vbasis[i],i=1..n);}{%
}
\end{mapleinput}

\mapleresult
\begin{maplelatex}
\mapleinline{inert}{2d}{x := x1*e1+x2*e2+x3*e3;}{%
\[
x := \mathit{x1}\,\mathit{e1} + \mathit{x2}\,\mathit{e2} + 
\mathit{x3}\,\mathit{e3}
\]
}
\end{maplelatex}

\end{maplegroup}
\begin{maplegroup}
\begin{mapleinput}
\mapleinline{active}{1d}{LHS:=simplify(ge(gSpin,t) &c (1 + x &c e4) &c star(ge(gSpin,-t))):
RHS:=simplify(Id+rigid(x,gSpin,t) &c e4):
simplify(LHS-RHS);}{%
}
\end{mapleinput}

\mapleresult
\begin{maplelatex}
\mapleinline{inert}{2d}{0;}{%
\[
0
\]
}
\end{maplelatex}
\end{maplegroup}
\noindent
Finally, we will verify directly that the action $\bx \mapsto \bg\,\bx \,\Star{\bg} + \bt$ is a rigid motion. If we denote by $\bx_p,\by_p$ the images of $\bx,\by$ under this action, we will need to show that $\|\bx_p - \by_p\| = \|\bx - \by\|$ in the Euclidean norm.  We can compute the norm $\|\bx - \by\|$ by taking $\sqrt{(\bx - \by)\Star{(\bx - \by)}}$ in $\cl_{0,3,1},$ or by using a procedure {\tt distance}. 
\begin{maplegroup}
\begin{mapleinput}
\mapleinline{active}{1d}{y:=add(y.i*vbasis[i],i=1..n):
distance:=proc(x,y) sqrt(simplify(scalarpart(cmul(x-y,star(x-y))))) end:
xp:=rigid(x,gSpin,t):yp:=rigid(y,gSpin,t):
evalb(distance(x,y)=distance(xp,yp));}{%
}
\end{mapleinput}

\mapleresult
\begin{maplelatex}
\mapleinline{inert}{2d}{true;}{%
\[
\mathit{true}
\]
}
\end{maplelatex}
\end{maplegroup}
\noindent
Thus the action of $G$ defined as $\bx \mapsto \bg\,\bx \,\Star{\bg} + \bt$ is a rigid motion in $\BR^3.$ In view of the presence of the radical in $\cl_{0,3,1},$ this group $G$ is in fact a semi-direct product of $\Spin(3)$ and $\BR^3,$ that is of rotations and translations. It is well known of course that $\Spin(3) \rtimes \BR^3$ doubly covers the group of proper rigid motions $SE(3).$ We leave it as an exercise for the Reader to check in {\sc CLIFFORD} that the composition of two rigid motions is a rigid motions. 

\section{Summary}
\label{sec:sum}

The main purpose of this paper has been to show a variety of computational problems that can be approached with the symbolic package {\sc CLIFFORD}. It is through an extensive experimentation with the package that the results reported in \cite{AblFauser99a} were found; we have seen some of the computations that have led to them in Section \ref{sec:hecke}. In view of their complexity it is unlikely that they could have been performed by hand. Relation between the Clifford products in $\cl(g)$ and $\cl(B)$ through the Helmstetter's formula appears more clear once it has been checked with {\sc CLIFFORD}. The SVD of a matrix as performed in Section \ref{sec:svd} is clearly feasible in the Clifford algebra language; however it is not clear if the approach presented there is the best. More study would need to be done here in order to possibly simplify the computations, perform them uniquely in the Clifford algebra language, and possibly better utilize the nilpotent-idempotent basis in the Clifford algebra rather than the Grassmann basis. In robotics applications  presented in Section \ref{sec:e3}, {\sc CLIFFORD} appears to be a very convenient tool to carry out practical computations in the low dimensional algebras such as $\cl_{0,3}$ and $\cl_{1,3,1}.$

\section{Acknowledgments}
\label{sec:ack}

The Author would like to acknowledge collaboration with and contribution from Bertfried Fauser, Universit\"{a}t Konstanz, Fakult\"{a}t f\"{u}r Physik, Fach M678, 78457 Konstanz, Germany, to the development of the additional procedures used in Section (\ref{sec:hecke}). 

\section{Appendix 1}
\label{app:first}

In addition to the main package {\sc CLIFFORD}, in Section \ref{sec:helm} we have used new procedures from a supplementary package {\sc suppl}. These procedures are: {\tt cinvg}, {\tt cmulg}, {\tt LCg}, {\tt makeF}, {\tt RCg}, {\tt revg}, and {\tt splitB}.

\begin{itemize}

\item Procedure {\tt cinvg} finds the Clifford inverse with respect to the symmetric part $g$ of the bilinear form $B,$ that is, it finds the inverse (if it exists) of $u$ in $\cl(g).$

\item Procedure {\tt cmulg} performs Clifford multiplication with respect to the symmetric part $g$ of the bilinear form $B:$ It gives the Clifford product $\gprod{u}{v}$ for any two elements $u$ and $v$ in $\cl(g).$

\item Procedures {\tt LCg} and {\tt RCg} give the left and right contraction in $\cl(g).$

\item Procedure {\tt makeF} computes element $F \in \bigw^2 V$ defined in (\ref{eq:defF}). It uses additional procedures {\tt pairs} and {\tt mysign} from the package {\sc suppl}.

\item Procedure {\tt revg} performs the reversion anti-automorphism $\tilde{~}$ in $\cl(g).$

\item Procedure {\tt splitB} splits a bilinear form $B$ in $V$ into its symmetric part $g$ and its antisymmetric part $A$ which are returned as a sequence $g, A.$ If $B$ is purely symbolic and a second optional parameter (of any type) is used, in addition to $g$ and $A$ the output sequence contains two lists of symbolic substitutions that relate entries of $g$ and $A$ to the entries of $B.$  If $B$ is not assigned, in order for {\tt splitB} to work it internally assigns a blank $9 \times 9$ matrix to $B$ and then it calls itself.

\end{itemize}

\noindent
In Section \ref{sec:hecke} we have used the following additional procedures from {\sc suppl}:
\begin{itemize}
\item Procedure {\tt cliexpand} expands the given Clifford number $u \in \cl(B)$ from the default Grassmann basis to a Clifford basis consisting of un-evaluated Clifford products of the generators $\{\e_1,\ldots,\e_n\}.$ Procedure {\tt clieval} converts from the Clifford basis back to the Grassmann basis.

\item Procedure {\tt reversion} gives the reversion in $\cl(B).$

\item Procedure {\tt clisolve2} solves equations of the type $u=0$ for the unknown parameters 
in~$u.$ Its code is shown in Appendix~3.

\item Procedure {\tt defB} is needed to define a bilinear form $B$ in an even-dimensional vector space $V.$

\item Procedure {\tt bexpand} expands any element in the Hecke algebra $H_{\BF}(n,q)$ in terms of the Hecke $b$-basis.

\item Procedure {\tt alpha2} provides a Hecke algebra automorphism.

\end{itemize}
\noindent
In Section \ref{sec:svd} we have used the following additional procedures from a supplementary package {\sc avsd}.
\begin{itemize}

\item Procedure {\tt phi} provides an isomorphism between a matrix algebra and a Clifford algebra.

\item Procedure {\tt radsimplify} simplifies radical expressions in matrices and vectors.

\item Procedure {\tt assignL} is needed to write output from a Maple procedure {\tt eigenvects} in a suitable form, it sorts eigenvectors according to the corresponding eigenvalues, and it uses the Gram-Schmidt orthogonalization process, if necessary, to return a complete list of orthogonal eigenvectors.

\item Procedure {\tt climinpoly} belongs to the main package {\sc CLIFFORD}. It computes a minimal polynomial of any element of a Clifford algebra.

\item Procedure {\tt makediag} makes a ``diagonal" $\Sigma$ matrix consisting of singular values.

\item Procedure {\tt embed} embeds the given non-square matrix or a matrix of smaller dimensions into a $2^k \times 2^k$ matrix of smallest $k$ such that it can be mapped into a Clifford algebra.

\end{itemize}

\section{Appendix 2}
\label{app:second}

We display a Maple code of some more important procedures used in Section \ref{sec:hecke}. 

\noindent
Procedure {\tt bexpand} expands any element in the Hecke algebra in terms of the Hecke $b$-basis:
\begin{mapleinput}
\mapleinline{active}{1d}{suppl[bexpand]:=proc(X) local a,setbs,setbe,i,eq,T,sys,sol;
option remember;
setbs:=[Id,'b1','b2','b12','b21','b121']:
setbe:=[Id,b1,b2,cmul(b1,b2),b21,cmul(b1,b21)]:       
eq:=CS(X-add(a[i]*setbe[i],i=1..6));
T:=cliterms(eq);
sys:=\{coeffs(eq,T)\};
sol:=solve(sys,\{seq(a[i],i=1..6)\});
subs(sol,add(a[i]*setbs[i],i=1..6));
end:}{%
}
\end{mapleinput}
\noindent
Procedure {\tt alpha2} gives an automorphism in the Hecke algebra:
\begin{mapleinput}
\mapleinline{active}{1d}{suppl[alpha2]:=proc(X) local a,setbb,setbe,setba,i,eq,T,sys,sol;
option remember;
setbe:=[Id,b1,b2,cmul(b1,b2),b21,cmul(b1,b21)]:
setba:=[Id,-1/q*reversion(b1),-1/q*reversion(b2),1/q^2*reversion(setbe[4]),
        1/q^2*reversion(setbe[5]),-1/q^3*reversion(setbe[6])]:        
eq:=CS(X-add(a[i]*setbe[i],i=1..6));
T:=cliterms(eq);
sys:=\{coeffs(eq,T)\};
sol:=solve(sys,\{seq(a[i],i=1..6)\});
subs(sol,add(a[i]*setba[i],i=1..6));
end:}{%
}
\end{mapleinput}
\noindent
Aliases needed in Section (\ref{sec:Garnir}.)

\begin{maplegroup}
\begin{mapleinput}
\mapleinline{active}{1d}{alias(H1=h[1]*q^3+2*h[6]*q-h[6]+h[5]*q^3-h[3]*q^2-2*h[5]*q^2-2*h[6]*q^2
\phtab{x} +h[5]*q+h[6]*q^3+h[4]*q^3+h[2]*q^3-h[2]*q^2+h[4]*q+h[3]*q^3-2*h[4]*q^2):
alias(H2=-2*h[6]*q+h[6]-h[4]*q+h[6]*q^2+h[2]*q^2+h[4]*q^2+h[5]*q^2-h[5]*q):
alias(H3=h[3]*q^2-h[5]*q-h[4]*q-2*h[6]*q+h[4]*q^2+h[6]*q^2+h[5]*q^2+h[6]):
alias(H4=-h[6]+h[5]*q+h[6]*q):
alias(H5=h[6]*q+h[4]*q-h[6]):
alias(H6=h[6]):}{%
}
\end{mapleinput}
\end{maplegroup}

\noindent
Aliases needed in Section (\ref{sec:Garnir}):

\begin{maplegroup}
\begin{mapleinput}
\mapleinline{active}{1d}{alias(w1=q^4*K[4]^2+3*q^3*K[4]^2-K[4]*q^3+4*q^2*K[4]^2-K[4]*q^2
\phtab{x} +3*q*K[4]^2-K[4]*q-q+K[4]^2-K[4]):
alias(w2=K[4]*q^3+2*K[4]*q^2-q^2+2*K[4]*q+K[4]):
alias(w3=-K[4]^2+K[4]+K[2]+K[4]*K[5]*q-q^3*K[4]*K[2]-q*K[2]*K[4]
\phtab{x} +K[6]*q^2*K[4]+K[4]*q^4*K[6]+K[6]*q^3*K[4]+K[4]*K[6]*q-q^2*K[2]*K[4]
\phtab{x} +q^2*K[4]*K[5]-K[2]*q-K[5]*q+K[4]*q^2-K[4]*q+q^4*K[4]^2-q^2*K[4]^2
\phtab{x} -q*K[4]^2-K[6]*q+K[6]*q^2-K[4]*K[2]);
alias(w4=K[6]*q^3+K[4]*q^3-K[6]*q^2-K[5]*q^2-K[6]*q-K[5]*q+K[6]+K[4]):
alias(w5=K[4]*q+K[2]*q-K[6]-K[5]):
alias(w6=K[6]*q^2+K[4]*q^2-2*K[6]*q-K[5]*q-K[4]*q+K[6]+K[4]):
alias(w7=K[4]*q^3+2*K[4]*q^2+q+2*K[4]*q+K[4]):
alias(w8=q^5*K[4]^2+3*q^4*K[4]^2+4*q^3*K[4]^2+K[4]*q^3+3*q^2*K[4]^2
\phtab{x} +K[4]*q^2+q*K[4]^2+K[4]*q-1+K[4]):
alias(w9=K[4]*q^3+K[2]*q^3-2*K[6]*q^2-K[5]*q^2-K[4]*q^2+K[5]*q
\phtab{x} +2*K[6]*q+K[4]*q-K[4]-K[6]):}{%
}
\end{mapleinput}
\end{maplegroup}

\noindent
Aliases needed in Section (\ref{sec:pin3}):

\begin{maplegroup}
\begin{mapleinput}
\mapleinline{active}{1d}{alias(kappa1=RootOf(_Z^2-1+x3^2+x2^2+x4^2)):
alias(kappa2=RootOf(_Z^2-1+x5^2+x3^2+x2^2)):
alias(kappa3=RootOf(-x6^2+x5^2*x6^2+x4^2*x5^2+x6^4+x2^2*x6^2+x4^2*x6^2+_Z^2)):
alias(kappa4=RootOf(-x7^2+x5^2*x7^2+x3^2*x7^2+x6^2*x7^2+x7^4+x4^2*x5^2
\phtab{x} -2*x4*x5*x6*x3+x6^2*x3^2+x4^2*x7^2+_Z^2)):
alias(kappa5=RootOf((x7^2+x8^2)*_Z^2+(2*x5*x4*x7-2*x3*x6*x7)*_Z-x8^2+
\phtab{x} +x6^2*x3^2+x8^4+x4^2*x5^2-2*x4*x5*x6*x3+x6^2*x8^2+x7^2*x8^2+x5^2*x8^2
\phtab{x} +x3^2*x8^2+x4^2*x8^2)):
alias(eps=RootOf(_Z^2-1)):
alias(lambda1=RootOf(\_Z\^{}2+x2\^{}2+x3\^{}2-1)):
alias(lambda2=RootOf(_Z^2-1+x3^2)):
alias(lambda3=RootOf(_Z^2-1+x5^2)):
alias(lambda4=RootOf(-x6^2+x5^2*x6^2+x4^2*x5^2+x6^4+x2^2*x6^2+x4^2*x6^2+_Z^2)):
alias(lambda5=RootOf(_Z^2-x6^2+x5^2*x6^2+x6^4)):
alias(lambda6=RootOf(-x7^2+x5^2*x7^2+x3^2*x7^2+x6^2*x7^2+x7^4+x4^2*x5^2
\phtab{x} -2*x4*x5*x6*x3+x6^2*x3^2+x4^2*x7^2+_Z^2)):
alias(lambda7=RootOf(_Z^2-x7^2+x5^2*x7^2+x6^2*x7^2+x7^4)):
alias(lambda8=RootOf((x7^2+x8^2)*_Z^2+(2*x5*x4*x7-2*x3*x6*x7)*_Z-x8^2
\phtab{x} +x6^2*x3^2+x8^4+x4^2*x5^2-2*x4*x5*x6*x3+x6^2*x8^2+x7^2*x8^2+x5^2*x8^2
\phtab{x} +x3^2*x8^2+x4^2*x8^2)):
alias(lambda9=RootOf(_Z^2*x8^2-x8^2+x8^4+x3^2*x8^2+x4^2*x8^2)):
}{%
}
\end{mapleinput}
\end{maplegroup}

\section{Appendix 3}
\label{app:three}
We display code for two procedures {\tt rot} and {\tt qrot} that were needed in Section \ref{sec:spin3}. Procedure {\tt rot} performs a rotation of a vector by a quaternion through a certain angle.
\begin{maplegroup}
\begin{mapleinput}
\mapleinline{active}{1d}{rot:=proc(v,q) local qs;
qs:=star(q):
RETURN(map(factor,clicollect(simplify(cmul(q,v,qs)))))
end:}{%
}
\end{mapleinput}
\end{maplegroup}
\noindent
Procedure {\tt qrot} finds a unit quaternion that is dual to the rotation axis vector ${\bf axis}.$
\begin{maplegroup}
\begin{mapleinput}
\mapleinline{active}{1d}{qrot:=proc(p1,p2,p3,theta) local bas,c,e,k,l,i,q,n;global qaxis;
if type(p1,name) or type(p2,name) then 
\phtab{xx} RETURN(cos(theta/2)*Id+sin(theta/2)*qaxis) fi;
if evalb(simplify(p1^2+p2^2+p3^2)=1) then 
\phtab{xx} q:=simplify(subs(\{a1=p1,a2=p2\},qaxis)) fi;
n:=sqrt(p1^2+p2^2+p3^2):
if n=0 then ERROR(`axis vector must be a non-zero vector`) elif
\phtab{xx} has(n,RootOf) then n:=max(allvalues(n)) fi;
q:=simplify(subs(\{a1=p1/n,a2=p2/n\},qaxis));
bas:=[Id,e12,e13,e23]:c:=[]:k:=0:
for i from 1 to 4 do l:=coeff(q,bas[i]);
\phtab{xxxxxxxxxxxxxxxxxx} if has(l,RootOf) then k:=i fi:
\phtab{xxxxxxxxxxxxxxxxxx} c:=[op(c),l] 
\phtab{xxxxxxxxxxxxxxx} od:
if k<>0 then e:=allvalues(c[k]);
\phtab{xxxx} if p3>0 then c:=subsop(k=max(e),c) elif
\phtab{xxxxxxx} p3<0 then c:=subsop(k=min(e),c) else ERROR(`p3=0`) fi;
fi;
q:=add(c[i]*bas[i],i=1..4);
RETURN(cos(theta/2)*Id+sin(theta/2)*q)
end:}{%
}
\end{mapleinput}
\end{maplegroup}
\noindent
Procedure {\tt clisolve2} was used extensively throughout this paper to solve linear equations in a Clifford algebra $\cl(B).$

\begin{mapleinput}
\mapleinline{active}{1d}{suppl[clisolve2]:=proc(eq,indet) local i,T,vars,sol,sys;
if type(indet,list) then vars:=convert(indet,set)
\phtab{xxxxxxxxxxxxxxxxx} else vars:=select(type,indets(indet),indexed)
fi;
T:=cliterms(eq);
sys:=\{coeffs(clicollect(simplify(eq)),T)\};
sol:=[solve(sys,vars)];
if type(indet,list) then RETURN(sol)
\phtab{xxxxxxxxxxxxxxxxx} else RETURN([seq(subs(sol[i],indet),i=1..nops(sol))]);
fi;
end:}{%
}
\end{mapleinput}

\begin{mapleinput}
\end{mapleinput}

\end{document}

%% file: Makro98.tex
\newcommand{\JJ}{\mathbin{\raisebox{0.25ex}{$\footnotesize
                       \rm\vphantom{I}%
                       \_\hskip -0.25em\_%
                       \vrule width 0.6pt$}}}           

\newcommand{\LL}{\mathbin{\raisebox{0.25ex}{$\footnotesize
                        \rm\vphantom{I}%
                        \vrule width 0.6pt \hskip -0.5pt%
                        \_\hskip -0.25em\_$}}}          

\newcommand{\w}{\wedge}
\newcommand{\bigw}{\bigwedge}

  \newcommand{\be}{{\bf e}}
                           \newcommand{\e}{{\bf e}} 

  \newcommand{\bg}{{\bf g}}
  
  \newcommand{\bi}{{\bf i}}
  \newcommand{\bj}{{\bf j}}
  \newcommand{\bk}{{\bf k}}

  \newcommand{\bq}{{\bf q}}

  \newcommand{\bt}{{\bf t}}
  
  \newcommand{\bv}{{\bf v}}
  
  \newcommand{\bx}{{\bf x}}
  \newcommand{\by}{{\bf y}}

\newcommand{\cl}{C \kern -0.1em \ell}     
\newcommand{\Spin}{{\bf Spin}}            \newcommand{\Pin}{{\bf Pin}}

\newcommand{\Mat}{{\rm Mat}}

\newcommand{\ut}[1]{{\setbox0=\hbox{$#1$}\mathsurround=0pt
       \rlap{\raisebox{-0.8\dp0}{\raisebox{-0.8ex}
       {\kern -0.15ex\hbox{$\tiny\sim$}\kern 0.15ex}}}#1}}
\newcommand{\uti}[1]{{\setbox0=\hbox{$#1$}\mathsurround=0pt
       \rlap{\raisebox{-0.8\dp0}{\raisebox{-0.8ex}
       {\kern -0.3ex\hbox{$\tiny\sim$}\kern 0.3ex}}}#1}}

\newdimen\arrayruleHwidth                 
     \makeatletter
     \def\Hline{\noalign{\ifnum0=`}\fi\hrule \@height \arrayruleHwidth
         \futurelet \@tempa\@xhline}
     \makeatother

%% file: CLIFF101.bbl
\begin{thebibliography}{99}

\bibitem{Abl96} R.~Ab\l amowicz; {\it Clifford algebra computations with Maple\/}, Proc. Clifford (Geometric) Algebras, Banff, Alberta Canada, 1995. Ed. W. E. Baylis, Birkh\"auser, Boston, 1996, pp. 463--501.

\bibitem{Structure} R.~Ab\l amowicz; {\it Structure of spin groups associated with degenerate Clifford algebras\/}, Journal of Mathematical Physics, Vol. {\bf 27}, No. {\bf 1}, January 1986, 
pp. 1--6.

\bibitem{Deformation} R.~Ab\l amowicz; {\it Deformation and contraction in Clifford algebras\/}, 
Journal of Mathematical Physics, Vol. {\bf 27}, No. {\bf 2}, January 1986, 
pp. 1--6.

\bibitem{Matrix} R.~Ab\l amowicz; {\it Matrix exponential via Clifford algebras\/}, Journal of Nonlinear Mathematical Physics, Vol. {\bf 27}, No. {\bf 2}, August 1998, pp. 423--427.

\bibitem{Symbolic} R.~Ab\l amowicz; {\it Spinor Representations of Clifford Algebras:
A Symbolic Approach\/}, CPC Thematic Issue ``Computer Algebra in Physics Research'', Physics Communications {\bf 115} (1998), pp. 510--535.

\bibitem{CLIFFORD} R.~Ab\l amowicz; {\sc CLIFFORD} - Maple V package for
Clifford algebra computations, ver. 4 (Copyright 1995-1999) and two supplementary packages {\sc suppl} and {\sc asvd} are available from \underline {{\protect \tt http://math.tntech.edu/rafal/cliff4/}}.

\bibitem{AblFauser99a}
R.~Ab\l amowicz and B.~Fauser, {\it Hecke algebra representations in ideals generated by q-Young Clifford idempotents\/}, Proceedings of the 5th International Conference on Clifford Algebras, Ixtapa, Mexico, 1999, Vol. {\bf 1} (submitted) and math.QA/9908062.

\bibitem{AblLoun85}
R.~Ab\l amowicz and P.~Lounesto, 
{\it Primitive idempotents and indecomposable left ideals in degenerate Clifford algebras\/}, Proc. of ``Clifford Algebras and Their Applications in Mathematical Physics (Canterbury, 1985), pp. 61--65, NATO Adv. Sci. Inst. Ser. C: Math. Phys. Sci., 183, Reidel, Dordrecht-Boston, Mass.

\bibitem{AblLoun96}
R.~Ab\l amowicz and P.~Lounesto, {\it On Clifford algebras of a bilinear form with an antisymmetric part\/}, in: R.~Ab\l amowicz, P.~Lounesto, and J.~M.~Parra, eds., {\em On Clifford Algebras with Numeric and Symbolic Computations\/} (Birkh\"{a}user, Boston, 1996), pp. 167--188.

\bibitem{LSI} M.~Berry and J.~Dongarra; {\it Atlanta Organizers Put Mathematics to Work for the Math Sciences Community\/}, SIAM News, Vol. {\bf 32}, No. {\bf 6}, July/August 1999, and references therein.

\bibitem{Bourbaki} N.~Bourbaki; {\it Algebra 1\/}, Chapters 1--3, Springer Verlag, Berlin, 1989.

\bibitem{Brooke7880} J.A.~Brooke; {\it A Galileian formulation of spin. I. Clifford algebras and Spin groups\/},  J. Math. Phys. {\bf 19}, no. {\bf 5}, 952--959 (1978); {\it A Galileian formulation of spin. I. Explicit realizations\/},  J. Math. Phys. {\bf 21}, no. {\bf 4}, pp. 617--621 (1980).

\bibitem{Brooke85} J.A.~Brooke; {\it Spin groups associated with degenerate orthogonal spaces\/}, Proc. of ``Clifford Algebras and Their Applications in Mathematical Physics (Canterbury, 1985), pp. 93--102, NATO Adv. Sci. Inst. Ser. C: Math. Phys. Sci., 183, Reidel, Dordrecht-Boston, Mass.

\bibitem{Crum90}
A.~Crumeyrolle, {\em Orthogonal and Symplectic Clifford Algebras: Spinor Structures\/} (Kluwer, Dordrecht, 1990).

\bibitem{CurtisReiner} C.W.~Curtis and I.~Reiner; {\it Methods of Representation Theory with Applications to Finite Groups and Orders\/}, Vol. 1, Wiley-Interscience, Now York, 1981.

\bibitem{Herstein} I.N.~Herstein; {\it Noncommutative Rings\/}, The Carus Mathematical Monographs, Number 15, Mathematical Association of America, Chicago, 1968. 

\bibitem{Fauser-hecke} B.~Fauser; {\it Hecke algebra representations within Clifford geometric algebras of multivectors\/}, J.~Phys.~A: Math. Gen. {\bf 32}, 1999, pp. 1919--1936.

\bibitem{Hamermesh} M.~Hamermesh; {\it Group Theory and Its application to Physical Problems\/}, Addison-Wesley, London, 1962.

\bibitem{Helm87}
J.~Helmstetter, {\it Mono\"{\i}des de Clifford et d\'{e}formations d'alg\`{e}bres de Clifford\/}, Journal of Algebra, Vol. {\bf 111} (1987), pp. 14--48.

\bibitem{KingWybourne} R.C.~King, B.G.~Wybourne; {\it Representations and traces of 
Hecke algebras $H_n(q)$ of type ${\bf A}_{n-1}$\/}, J. Math. Phys. {\bf 33}, 1992, pp. 4--14.

\bibitem{Clical} P.~Lounesto, R.~Mikkola, and V.~Vierros, `CLICAL User Manual', Helsinki
University of Technology, Institute of Mathematics, Research Reports {\bf A248}, Helsinki, 1987.

\bibitem{Lounesto} P.~Lounesto; {\it Clifford Algebras and Spinors\/}, Cambridge University Press, Cambridge, 1997.

\bibitem{Lounesto_private} P.~Lounesto; Private communication, 1997.

\bibitem{Macdonald} I.G.~Macdonald; {\it Symmetric Functions and Hall Polynomials\/}, Oxford University Press, Oxford, 1979.

\bibitem{Maciejowski} J.M.~Maciejowski; {\it Multivariable Feedback Design\/}, Addison-Wesley, Wokingham, England, 1989.

\bibitem{Port95}
I.R.~Porteous: {\em Clifford Algebras and the Classical Groups\/} (Cambridge University Press, Cambridge, 1995).

\bibitem{Selig} J.M.~Selig; {\it Geometrical Methods in Robotics\/}, Monographs in Computer Science, Springer-Verlag, New York, 1996.

\bibitem{Strang} G.~Strang; {\it Introduction to Linear Algebra\/}, Wellesley-Cambridge Press, Wellesley, 1998.

\end{thebibliography}
